\RequirePackage{fix-cm}
\documentclass[twocolumn]{svjour3}          

\smartqed  
\usepackage[misc,geometry]{ifsym}
\usepackage[utf8]{inputenc}
\usepackage{amsfonts, amssymb, amsmath, mathtools}
\usepackage{mathrsfs}
\usepackage{graphicx,graphics}
\usepackage{placeins}
\usepackage{epstopdf}
\usepackage{multirow}
\usepackage{subfig}
\usepackage{setspace}
\usepackage{anysize}
\usepackage{color}
\usepackage{calrsfs}
\usepackage{hyperref}
\usepackage{framed}
\usepackage{bm}
\usepackage{algorithm}
\usepackage{algpseudocode}

\usepackage{calrsfs}
\DeclareMathAlphabet{\pazocal}{OMS}{zplm}{m}{n}

\marginsize{2cm}{2cm}{2cm}{2cm} 
 \providecommand{\keywords}[1]{\textbf{\textit{Index
terms---}} #1}

\usepackage{etoolbox}
\newcommand*{\affaddr}[1]{#1} 
\newcommand*{\affmark}[1][*]{\textsuperscript{#1}}

\definecolor{edit}{RGB}{0, 61, 153}
\definecolor{remove}{RGB}{255, 78, 0}
\definecolor{Sampsa}{RGB}{0, 153, 61}

\begin{document}

\title{Conditionally Exponential Prior in 
Focal Near- and Far-Field EEG Source Localization via Randomized Multiresolution Scanning (RAMUS)
}

\author{Joonas Lahtinen\protect\affmark[a] \Letter         \and
        Alexandra Koulouri\affmark[a] \and
        Atena Rezaei\affmark[a] \and
        Sampsa Pursiainen\affmark[a]}


\institute{Joonas Lahtinen \at
              \email{joonas.j.lahtinen@tuni.fi}           
           \and
           Alexandra Koulouri \at
             \email{alexandra.koulouri@tuni.fi}
           \and
           Atena Rezaei \at
              \email{atena.rezaei@tuni.fi}
           \and
           Sampsa Pursiainen \at
              \email{sampsa.pursiainen@tuni.fi}
            \and 
           \affaddr{\affmark[a]Computing Sciences, Faculty of Information Sciences, Tampere University, Korkeakoulunkatu 1, 33014 Tampere}}

\date{Received: date / Accepted: date}

\maketitle

\begin{abstract}
In this paper, we focus on  the inverse problem of reconstructing  distributional brain activity with cortical and weakly detectable deep components in non-invasive Electro\-ence\-phalo\-graphy. We consider a recently introduced hybrid reconstruction strategy combining a hierarchical Bayesian model to incorporate {\em a  priori} information and the advanced randomized multiresolution scanning (RAMUS) source space decomposition approach to reduce modelling errors, respectively. In particular, we aim to  generalize the previously extensively used conditionally Gaussian prior (CGP) formalism  to achieve distributional reconstructions with higher focality. For this purpose, we  introduce as  a hierarchical prior, a  general exponential distribution, which we refer to as conditionally exponential prior (CEP). The first-degree CEP corresponds to focality enforcing Laplace prior, but  it also suffers from strong depth bias, when applied in numerical modelling, making the deep activity unrecoverable. We  sample over multiple resolution levels via RAMUS to reduce this bias as it is known to depend on the resolution of the source space. Moreover, we introduce a procedure based on the physiological {\em a priori} knowledge of the brain activity to obtain the shape and scale parameters of the gamma hyperprior that steer the CEP. The posterior estimates  are calculated using iterative statistical methods, expectation maximization and iterative alternating sequential algorithm, which we show to be algorithmically similar and to have a close resemblance to the iterative $\ell_1$ and $\ell_2$ re\-weighting methods. The performance of CEP is compared with the recent sampling-based dipole localization method {\em Sequential semi-analytic Monte Carlo estimation} (SE\-SA\-ME) in numerical experiments of simulated somatosensory evoked potentials related to the human median nerve stimulation. Our results obtained using synthetic  sources  suggest that a hybrid of the first-degree CEP and RAMUS can achieve an accuracy comparable to the second-degree case (CGP) while being more focal. Further, the proposed hybrid  is shown to be robust to noise effects and  compare well to the dipole reconstructions obtained with SESAME.

\end{abstract}

\section{Introduction}

This article aims to advance mathematical inverse meth\-odology in focal localization of brain activity at different depths based on non-invasive electro-ence\-phalo\-graphy (EEG) measurements \cite{niedermeyer2004}. In particular, we consider  the challenging task of reconstructing  simultaneous cortical and sub-cortical brain activity by comparing advanced reconstruction algorithms to find a solution. The EEG source localization task generally constitutes an ill-posed inverse problem \cite{hamalainen1993}, i.e.,  it does not have a unique solution and is sensitive to different modelling and measurement errors. Consequently, only a slight amount of noise in the measurement can significantly affect the reconstruction found by the inversion algorithm. This is especially the case for the far-field activity components, e.g., sub-cortical activity. The feasibility of depth-localization with non-invasive measurements has been suggested recently in studies concentrating on high-density measurements and filtering  \cite{seeber2019subcortical,pizzo2019deep}.

We approach the presented source localization task via the recently proposed \cite{rezaei2020randomized}  distributional  reconstruction strategy \cite{rezaei2020randomized} which utilizes a hierarchical Bayesian model (HBM) to incorporate  {\em a priori} information and, at the same time,  decomposes the source space into randomized sets at several different resolution levels to reduce modelling errors. This  hybrid technique  comprises the iterative alternating sequential (IAS) \cite{Calvetti2009} posterior optimization method  and the so-called randomized multiresolution scan (RAMUS), i.e., a Monte Carlo sampling approach which aims at  marginalizing the modelling erros over the source space decomposition. 
In this article, we introduce an extended version of this technique via the conditionally exponential prior (CEP), i.e., an exponential power distribution also referred to as a generalized normal distribution  \cite{nadarajah2005generalized}, which generalizes the well-known concept of the conditionally Gaussian prior (CGP) \cite{Calvetti2009}. CEP includes a variable {\em prior degree} parameter $q$ which defines the $\ell_q$-norm applied in the argument of the exponential prior probability density function. Bayesian formulation of the model allows us to use the methods of statistical analysis to find an estimate for the brain activity, especially, the expectation maximization (EM) and iterative alternating sequential (IAS) method which  are,  in the case of CEP, algorithmically similar to $\ell_q$-reweighted methods considered in \cite{Wipf&Nagarajan} for $q=1,2$. But as a difference, the hierarchical prior structure allows taking into account \cite{rezaei2020parametrizing} the physical and physiological properties of the underlying primary currents, here, in particular, the shape and scale parameter steering of the gamma hyperprior of CEP.

 Our focus is  in particular on  developing   new methodology that is applicable  in hybrid with the RAMUS technique to enhance the detectability of both cortical and sub-cortical activity \cite{rezaei2020randomized}. RAMUS differs from the other recently introduced  parcellation techniques \cite{parcellation-MattoutJ2005,parcellation-Chevalier-Alexis2020S} as it does not try to reduce variable dimension, i.e., source space size, or divide spatially the problem to sub-problems but rather gathers information from multiple randomized source spaces with different resolutions to  strengthen dipole localization at any depth. A similar method with a coarse-to-fine hierarchy has been shown to improve the detectability of sub-cortical activity in \cite{krishnaswamy2017sparsity}. However, this method limits its parcellation to cortical level and changes the whole inversion problem with its observations progressively, whereas RAMUS keeps observation untouched relying only on the sparsity  of the source space.

We show mathematically how the HBM approach, previously introduced in \cite{Calvetti2009}, can be extended to the case of the CEP and how the resulting statistical framework can be associated with the previously described CGP model and the reweighted posterior optimization methods that are applicable with RAMUS as they allow a hyperparameter progression over multiple resolution levels. In previous studies RAMUS, technique is applied with CGP \cite{rezaei2020randomized,rezaei2020parametrizing} concentrating on the detectability of the deep components while here we investigate  CEP as a potential improvement of the focality. As shown in \cite{rezaei2020randomized}, the numerical implementation of the method is an important factor contributing the eventual performance of the inverse model. Therefore, we investigate  EM and IAS as two alternative techniques for maximizing the posterior, and  compare the prior degrees $q=1$ and $q=2$  in reconstructing  different synthetic source configurations. In the former case, the reconstruction is found applying the Lasso algorithm  \cite{Tibshirani94regressionshrinkage}. Similar Laplace prior models have been previously studied in Bayesian framework with variance components as hyperparameters  in \cite{Lucka2,Lucka3}. 

In the numerical experiments, we aim  to find  the best combination of posterior maximizing algorithm and prior degree considering the estimation accuracy obtained for  cortical and sub-cortical sources. RAMUS is applied to reduce the modelling errors to improve the detection of deep components. As a forward modelling technique, we apply the finite element method (FEM)  \cite{miinalainen2019,pursiainen2016,demunck2012}  which is advantageous in the present application, since it, as a volumetric technique, allows decomposing the head model into multiple highly accurate cortical and sub-cortical compartments. The results are compared with  estimates obtained with the Sequential semi-analytic Monte Carlo estimation (SESAME) technique, which is  a recently introduced Monte Carlo based  algorithm \cite{sommariva2014} for Bayesian dipole  localization.

The synthetic measurements are produced according to the well-studied somatosensory evoked potentials (SEPs) occurring in the human median nerve stimulation \cite{niedermeyer2004,buchner1995somatotopy,buchner1994source,allison1991cortical,hari2017meg}. The median nerve SEPs have  known originators, i.e., locations of neural sources at given measurement time points. The early components occurring $\leq$ 20 ms post-stimulus involve sub-cortical far-field activity, i.e., activity far from the electrodes, which we here simulate and reconstruct numerically utilizing the CEP model together with a realistic head model segmentation. We consider especially the P14/N14, P16/N16 and P20/N20 components, i.e., the positive (P) and negative (N) 14, 16 and 20 ms post-stimulus peaks in the measured data with respect to the forehead potential. This setup not only fulfills the requirements of our numerical experiment,s but reflects to the source localization from real SEP data peaks \cite{RezaeiAtena2021Rsac}.  

The results obtained suggest that the first-degree CEP provides an advantageous  approach to detect focal near- and far-field activity when it is applied together with the RAMUS technique. Namely, while both  CEP prior degrees yield a similar source localization accuracy, the activity reconstructed via the first-degree CEP is overall more well-localized compared to CGP. This was found to be the case for simultaneous thalamic and cortical activity, approximating the simulated originators of P20/N20, as well as for a sub-thalamic dipolar or quadrupolar source configuration  corresponding approximately to  the originators of  P14/N14 and P16/N16 \cite{noel1996origin,buchner1995origin},  respectively. In comparison to SESAME, the proposed combination of RAMUS and CEP was found to be advantageous considering its robustness to high noise effects. In turn, SESAME can produce accurate localized activity in most cases, while it tends to show false activity on the cortical level in cases, where the cortical activity is absent (P14/N14 and P16/N16).

This article is organized as follows: In Section \ref{sec:Methods} we present a generalized exponential prior model to HBM based source localization and the EM and IAS algorithms as well as the main principles behind RAMUS and SESAME. Section \ref{sec:implementation} concentrates on the implementation of this methodology. The results are presented in Section \ref{sec:results} and discussed in Section \ref{sec:discussion}. Finally, Section \ref{sec:conclusion} concludes the study.

\section{Methods}\label{sec:Methods}
In this section, we introduce the hierarchical Bayesian framework for CEP model. Next, the statistical iterative maximum a posterior algorithms are derived for the model and, after that, we give some hypothetical and theoretical background for the randomized multiresolution scanning and the usage of multiple source space resolutions. Finally, we introduce Sequential semi-analytic Monte Carlo estimation and the assumptions behind it.

\subsection{Hierarchical Bayesian framework for a conditionally exponential prior}
\label{sec:hbm}

We consider the linear EEG observation model 
\begin{equation}
\label{eq:lead_field}
{\bf y}={\bf L} {\bf x}+{\bf e},
\end{equation}
where ${\bf x}\in \mathbb{R}^{3n}$, ${\bf y} \in\mathbb{R}^m$ and
$\pazocal{N}({\bf e};0,\sigma^2)$ with $\sigma$ a scale-invariant prior. From Bayes' rule we have $
\pi({\bf x} \mid {\bf y}) \propto \pi({\bf y} \mid {\bf x})\pi({\bf x})$. The leadfield matrix {\bf L} is obtained via the finite element method applied to Maxwell's equations in the quasi-static approximation. The unknown ${\bf x}$ represents the discretized primary current distribution ${\bf J}$ (vector field) of the neural activity in a three-dimensional source space with $n$ possible source locations. At a given position, ${\bf J}$ is described by a three-component vector or, equivalently, three entries of ${\bf x}$, i.e., a vector $
\left( x_{3\mu-2},\, x_{3\mu-1},\,  x_{3\mu} \right)$, where $\mu=1\ldots n$. Given the likelihood function 
$\pi({\bf y} \mid{\bf x})\propto$\\ 
\noindent $\exp{\left(-\frac{1}{2\sigma^2}\|{\bf L} {\bf x}-{\bf y}\|_2^2\right)}$ and a subjectively selected prior $\pi({\bf x})$, the posterior $\pi({\bf x} \mid {\bf y})$ is assumed to contain all the information about the underlying source activity. We associate each \textcolor{edit}{component} $x_i$ with 
an exponential power distribution 
\begin{equation}\label{eq:exponetialPrior}
\pi(x_i ,1/\gamma_{i},q)\propto\gamma_{i}^{1/q}\exp{\left(-{\gamma_{i}|x_i|^{q}}\right)}
\end{equation}
determined by the hyperparameter $\gamma_i$ and the degree $q$ of the prior, which is selected to be either one or two in this study. Our choice of hyperparameter differs from a common choice, separable prior variance variables, that are reciprocal to $\gamma_i$ \cite{hyperparam_SatoMasa-aki2004HBef,hyperparam_FristonKarl2007,hyperparam_Wipf&Nagarajan2008,Calvetti2009}. With this choice we gain two desirable properties: analytic expression for posterior maximizing hyperparameters and stability by avoiding tendency to infinity. Selecting between  $q=1$ and $q=2$  allows one to steer the focality of the reconstruction, since the exponential distribution can be justified to be more  heavy-tailed in the former case (i.e. Laplace distribution). By introducing an extra level of hierarchy in this prior
(\ref{eq:exponetialPrior}), we aim to obtain a minimization problem that
allows reconstructing more focal, intensity unbiased or deeper sources, as suggested in \cite{Calvetti2009,rezaei2020randomized}, compared to the minimization problem that employs a fixed $\gamma_{i}$ together with the  Laplace prior, when $q=1$, or the standard  Gaussian prior, when $q=2$. Consequently, the current analysis belongs to the
hierarchical Bayesian adaptive framework introduced, e.g., in \cite{Murphy2012}. In particular, with $\gamma_{i}$ being a random variable following
a Gamma hyperprior distribution, i.e., $\gamma_{i}\sim\mathrm{Ga}(\kappa_{},\theta_{})$
for $i=1,\ldots,3n$, where its density function is  
\begin{equation}
    \pi(\gamma_i)=\frac{\theta^\kappa\gamma_i^{\kappa-1}e^{-\theta \gamma_i}}{\Gamma(\kappa)}
\end{equation}
for $\kappa,\: \theta>0$. $\pi(\gamma_{i}|x_i)$ is also a Gamma
distribution, that is, 
$\gamma_{i}|x_i\sim\mathrm{Ga}(\gamma_{i}|x_i;\kappa_{}+1/q,\theta_{}+|x_i|^{q})$ by  conjugacy. The full posterior obeying a CEP 

\noindent $\pi({\bf x} \mid {\bf \gamma})$ is given by
\begin{equation}
\begin{split}
    \pi({\bf x},{\bm \gamma} \mid {\bf y}) &\propto \pi({\bf y}\mid {\bf x})\;\pi({\bf x}\mid{\bm \gamma})\;\pi({\bm \gamma})\\
    &=\pi({\bf y}\mid {\bf x}) \; \prod_{i=1}^{3n} \pi(x_i \mid\gamma_{i}) \;\pi(\gamma_{i}), 
\end{split}
\end{equation}
where $\pi({\bf x}\mid{\bm \gamma})$ and $\pi({\bm \gamma})$ correspond to the CEP following from (\ref{eq:exponetialPrior}) and the hyperprior, respectively. In the following subsections, we use two different approaches to estimate the mode of the marginal posterior $\pi({\bf x} \mid {\bf y})$ given the CEP. The first one relies on EM and the second one on IAS  \cite{calvetti2007,Calvetti2009,calvetti2018}. 

\subsection{EM for the hierarchical adaptive framework}
The EM-based maximum a posteriori (MAP) estimate  is given by the system
{\footnotesize
\begin{equation}
\label{eq:gamma_ex}
\begin{split}
\hat{\bf x}^{(j+1)}&= \,\underset{{\bf x}}{\hbox{arg max}} \Big\{ -\frac{1}{2\sigma^2}\|{\bf L}{\bf x}-{\bf y}\|_2^2 \\
&+\mathbb{E}_{\pi({\bm \gamma}|\hat{{\bf x}}^{(j)})}[\log
\pi({\bf x} \mid {\bm \gamma})]\Big\}, 
\end{split}
\end{equation}}
where the expectation of $\log
\pi({\bf x} \mid {\bm \gamma})$  with respect to the conditional probability density   $\pi({\bm \gamma}|\hat{\bf x}^{(j)})$ and $\pi(\gamma_{i}|\hat{x}_i^{(j)})$  given the MAP estimate $\hat{\bf x}_i^{(j)}$ for $i=1,\ldots,3n$ is 
{\footnotesize
\begin{equation}
\begin{split}
\label{eq:gamma_ex_2}
\mathbb{E}_{\pi({\bm \gamma}|\hat{\bf x}^{(j)})}[\log
\pi({\bf x} | {\bm \gamma})]&=\sum_{i=1}^{3n}\mathbb{E}_{\pi(\gamma_{i}|\hat{x}_i^{(j)})}[\log
\pi (x_i|\gamma_{i})]\\ &=\sum_{i=1}^{3n} \int_{0}^\infty
\pi(\gamma_{i}|\hat{x}_i^{(j)}) \log\pi
(x_i|\gamma_{i})\;d\gamma_{i} \\ &=- \! \sum_{i=1}^{3n} \! |x_i|^{q \! }\int_{0}^\infty \!
\gamma_{i} \, \pi(\gamma_{i}|\hat{x}_i^{(j)})
\; d\gamma_{i} \! + \! C,  
\end{split}
\end{equation}}
with  $\mathbb{E}_{\pi(\gamma_{i}|\hat{x}_i^{(j)})}[\gamma_{i}]=
\int_{0}^\infty
\gamma_{i}\pi(\gamma_{i}|\hat{x}_i^{(j)})
\; d\gamma_{i}= \frac{\kappa_{}+1/q}{\theta_{}+|\hat{x}_i^{(j)}|^q}$. Therefore, we have the following optimization problem:
\begin{equation}\label{eq:AdaptiveMinimization}
\begin{split}
\bar{\gamma}_i^{(j)} =& \frac{\kappa_{}+1/q}{\theta_{}+|\hat{x}_i^{(j)}|^q}\mbox{ for } i=1,\ldots,3n\\
\hat{\bf x}^{(j+1)}=&\,\underset{{\bf x}}{\hbox{arg min}} \left\lbrace \frac{1}{2\sigma^2}\|{\bf L} {\bf x}-{\bf y}\|_2^2+\sum_{i=1}^{3n}
\bar{\gamma}_i^{(j)} |x_i|^q \right\rbrace.
\end{split}
\end{equation}
When  $q=1$, this resembles the Lasso problem which,  expressed through (\ref{eq:EMminimization}),  corresponds to  a Laplace prior (\ref{eq:LaplacePrior}) with a fixed $\gamma_{i}$. 
Notice that EM (\ref{eq:AdaptiveMinimization})  finds the MAP estimate using as prior the marginal  distribution of $x_i$, i.e., $\pi(x_i)=\int_{\gamma_i}\pi(x_i|\gamma_i)\;\pi(\gamma_i)\;d\gamma_i$,
\begin{equation}\label{eq:t-distribution}
\pi(x_i;\kappa_{},\theta_{},q)\propto\left(\frac{|x_i|^q}{\theta_{}}+1\right)^{-(\kappa_{}+1/q)},
\end{equation}
where $\kappa_{}$ and $\theta_{}$ are the shape and scale
parameter of  $\mathrm{Ga}(\gamma_{i};\kappa_{},\theta_{})$, respectively.
\subsection{Iterative Alternating Sequential algorithm}
In IAS, we aim at estimating the MAP of the pair $({\bf x},\bm{\gamma})$ by solving the optimization problem $
(\hat{\bf x},\hat{\bm \gamma})=\mathrm{arg}\max_{{\bf x},{\bm \gamma}}\pi({\bf x},{\bm \gamma} \mid {\bf y})$.  A common procedure \cite{ohagan2004} is to evaluate it via alternating optimization with respect to  ${\bf x}$ and ${\bm \gamma}$ in the similar manner to \cite{DaubechiesIngridreweight,CandesEmmanuelJ2008ESbR}. In particular, the MAP estimates can be extracted by solving alternatingly and recursively the following two optimization problems
\begin{equation}\label{IAS_problem}
\begin{split}
\hat{\bm \gamma}^{(j)}=&\underset{{\bm \gamma}}{\hbox{arg max}}\log \pi({\bm \gamma} \mid {\bf y},\hat{\bf x}^{(j)}), \\
\hat{\bf x}^{(j+1)}=&\underset{{\bf x}}{\hbox{arg max}}\log \pi({\bf x} \mid {\bf y},\hat{\bm \gamma}^{(j)})  . 
\end{split}
\end{equation}
To express explicitly the two previous optimization problems, we write
the full posterior $\pi({\bf x},{\bm \gamma} \mid {\bf y})$ which is 
\begin{equation}
\begin{split}
    \pi({\bf x},{\bm \gamma} \mid {\bf y}) &\propto \exp\bigg(-\frac{1}{2\sigma^2}\|{\bf L} {\bf x}-{\bf y}\|_2^2\\
    &-\sum_{i=1}^{3n} \big( \gamma_{i}(|x_i|^q+\theta_{}) - (1/q+\kappa_{}-1)\log\gamma_{i} \big)
\bigg). 
\end{split}
\end{equation}
It follows that the optimization problem  (\ref{IAS_problem})  can be written as 
\begin{equation}\label{eq:IAS}
\begin{split}
\hat{\bf \gamma}^{(j)}_i=& \frac{\kappa_{}+1/q-1}{|\hat{x}^{(j)}_i|^q+\theta_{}}\mbox{ for } i=1,\ldots,3n,\\
\hat{\bf x}^{(j+1)}=&\underset{{\bf x}}{\hbox{arg min}} \left\lbrace \frac{1}{2\sigma^2}\|{\bf L}{\bf x}-{\bf y}\|_2^2+\sum_{i=1}^{3n}\hat{\gamma}_i^{(j)}|x_i|^q \right\rbrace.
\end{split}
\end{equation}

\subsection{Difference and similarity between EM and IAS approaches}

One can observe that the difference between the problems of  (\ref{eq:AdaptiveMinimization}) and  (\ref{eq:IAS}) is step updating the $\gamma$ parameter. In EM algorithm, the update of $\gamma_{i}$ is based on the estimation of the expectation of $\gamma_{i}$ (see \ref{eq:gamma_ex}) whereas in IAS, the update comes from the mode of $\pi(\gamma_{i}|x_i)$.
Since we have that $\gamma_{i}|x_i\sim\mathrm{Ga}(\gamma_{i};\kappa_{}+1/q,\theta_{}+|x_i|^q)$, we saw that the two updates can be explicitly expressed. In this article, we  investigate  how the prior degree $q=1$ or $q=2$ affects the source reconstructions and  how the mode or expectation of $\gamma_{i}$ as updating rules influence the performance of an EEG source localization solver. 
With these prior degrees the algorithms resembles the separable $\ell_1$ and $\ell_2$ reweighting algorithms that are examined in \cite{Wipf&Nagarajan}.

\subsection{Randomized multiresolution scanning (RAMUS)}

The RAMUS technique was applied in finding a MAP estimate in order to maximize the robustness of the source localization outcome for various source depths \cite{rezaei2020randomized}. That is, each MAP estimate was found for a large number of subsets of the source space and the final estimate was found as the average of all these subset-based estimates. In RAMUS, each  subset contains a given number of randomly and uniformly distributed source positions. These subsets are divided to resolution levels according to the  source position count. The estimates are evaluated in an ascending order with respect to the resolution, i.e., progressing from coarse (sparse) to fine resolutions. Furthermore, the hyperparameter estimate obtained at one resolution level is used as the initial guess for the next one. The presence of coarse resolutions can be shown to be essential, specially, regarding the distinguishability of the deep activity \cite{krishnaswamy2017sparsity}, whereas finer resolutions can provide an enhanced accuracy for the detection of the near-field (cortical) activity. The average estimate found for a randomized set of source sub-spaces with a given number of sources provides an enhanced robustness for each resolution, as it diminishes discretization and modelling errors \cite{rezaei2020randomized}. We consider RAMUS as a modelling error reduction technique based on the postulates presented in the following subsections. 

\subsubsection{Marginalization of source space related modelling errors} 
We assume  that the estimated unknown variable $\hat{\bf x}$  depends on uncertainties and modelling errors of the forward model ($\mathrm{FM}$), inverse algorithms ($\mathrm{IA}$), and the statistical inverse model itself ($\mathrm{IM}$). Consequently, the modelling error is of the form 
\begin{equation}
    \bm{\varepsilon}({\hbox{\small FM}},\hbox{\small IA},\hbox{\small IM})=\hat{{\bf y}}-{\bf L}(\hbox{\small FM}) \, \hat{{\bf x}}(\hbox{\small FM},\hbox{\small IA},\hbox{\small IM}),
\end{equation}
where $\hat{{\bf y}}$ is noiseless data. As the statistical inverse model is not affected by the discretization of the leadfield, we assume the inverse error to be identically and independently distributed with respect to the selection of the  source space. Consequently, a sample mean estimator obtained through a  Monte Carlo sampling process can be interpreted to  marginalize the error over the source space configurations that are uniformly distributed \cite{rezaei2020randomized}.

\subsubsection{A coarse-to-fine optimization process}
RAMUS uses the hyperparameter $\bm{\gamma}$ as a surrogate \cite{rezaei2020parametrizing}  model to reduce the optimization bias towards superficial brain regions, which otherwise occurs  with high-resolution source spaces \cite{krishnaswamy2017sparsity};  the  hyperparameter  obtained after processing one resolution level constitutes the initial guess for the hyperparameter on the next level to maintain the activity of the deep structures  found with low-resolutions when the optimizer proceeds towards the finest resolution level \cite{rezaei2020randomized}. In this study, each source set in the coarsest resolution level includes 10 source positions  and the source count is multiplied by a factor of 10 (sparsity factor) when moving up by one  resolution level in the multiresolution hierarchy. The results were averaged over 100 different multiresolution decompositions.

\subsubsection{Sparse source distinction}
The use of coarse resolution levels below the number of measurement sensors is justified by the well-posedness of the corresponding  Diriclet boundary value (forward) problem \cite{evans10}. In particular, the diffusion operator following from the finite element discretization of the Maxwell's equations under the quasi-static approximation is symmetric, continuos (bounded) and elliptic (coersive). By the Lax-Milgram theorem, given a  Lipschitz domain $\Omega$, positive-valued Lebesque integrable $L_\infty(\Omega)$ function $\sigma$, and a diverge conforming \cite{pursiainen2016} primary current distribution with a normalized and  square-integrable $L_2(\Omega)$ divergence $f$, $\| f \| =1$, the weak form $\int_\Omega \sigma \nabla u \cdot \nabla v \, d\Omega = \int_\Omega f v \, \hbox{d} \Omega$ has a unique solution $u$ which belongs to the Sobolev space $H^1(\Omega)$ or to its finite dimensional subspace $S$ such that the weak form is satisfied for any $v$ in $H^1(\Omega)$ or in $S$, respectively. Due to this well-posedness, any two boundary datasets $g_1$ and $g_2$ which tend infitesimally closer and closer  $g_2 - g_1 \to 0$ correspond to solutions $u_1$ and $u_2$ with a similar property  $u_1 - u_2 \to 0$ and sources whose difference tends to zero weakly, i.e.,   $ \big|\int_\Omega (f_2 - f_1) v  \, \hbox{d} \Omega \big| \to 0$ for any $v$ in $H^1(\Omega)$  or in $S$. It follows that each source in an $N$ dimensional subspace $S$, sources $f_1, f_2, \ldots, f_K$  can  be distinguished ($|\int_\Omega (f_i - f_j) v \, \hbox{d}\Omega| \geq \varepsilon$ for some $\varepsilon > 0$ for all non-zero $v$ in $S$ and any pair $i \neq j$) based on their boundary datasets $g_1, g_2, \ldots, g_K$, if $S$ is spanned by a set of functions $v_1,v_2, \ldots, v_N$ where any $v_i$, $i = 1, 2, \ldots, N$ overlaps with maximally one source $f_j$, i.e., if the distribution of the sources is sparse enough. With a finite number of measurements, the maximum number of uniquely  distinguishable sources is limited by the number of sensors, i.e., the number of linearly independent data vectors.

\subsection{Sequential semi-analytic Monte Carlo estimation (SESAME)}
When a posterior distribution of the source localization problem cannot be expressed analytically, one possibility to approximate such a distribution is to use a class of particle filters called {\em sequential Monte Carlo sampler} algorithms, where one constructs a sequence of artificial distributions that starts from an analytical distribution and converges smoothly to the true distribution. SESAME by Sorrentino et al.\ \cite{SorrentinoAlberto2014Bmmo,sommariva2014} utilizes this idea  via   the following artificial posterior:
\begin{equation}
    \pi_i({\bf x}\mid {\bf y})=\pi_0({\bf x})\pi_i({\bf y}\mid {\bf x})^{f(i)},
\end{equation}
in which the power $f(i)$ of the likelihood is given by
\begin{equation}
    f(i)=\sum_{k=1}^i\delta_k,\quad \sum_{k=1}^{n}\delta_k=1
\end{equation}
and $n$ is the total number of iterations, sufficient to bring the estimate satisfyingly close to the true posterior, and increments $\delta_k$ vary uniformly in the fixed interval $\left[10^{-5},10^{-1}\right]$. The prior density $\pi_0$ is based on  the following assumptions: \cite{SorrentinoAlberto2014Bmmo}
\begin{enumerate}
\item the number of dipoles is Poisson distributed;
\item the distribution of dipole locations and orientations is uniform;
\item dipole amplitude is random variable $10^{3U}$, where $U$ is uniform distributed in the interval $(0,1)$;
\item dipole location, amplitude and angle are independent from each other.
\end{enumerate}
The final dipole estimates are obtained via importance sampling, i.e., as the weighted mean values from the conditional distributions known as {\em importance}\\
\noindent {\em weights}. Unlike CEP, SESAME does not provide a sparse solution, but a finite number of estimated dipole positions and orientations that best explain the measurement data. Therefore SE\-SA\-ME needs the whole source space and leadfield to work with, and applying a multiresolution decomposition to it does not necessarily  improve the estimation by design.

\section{Implementation}
\label{sec:implementation}

The EM, IAS and SESAME source localization approaches were implemented as additional packages (plugins) of the Matlab based {\em Zeffiro interface} (ZI) code package\footnote{\url{https://github.com/sampsapursiainen/zeffiro_interface}}  \cite{he2019zeffiro} which allows for using the detailed multi-compartment head models obtained via high resolution magnetic resonance imaging (MRI) data. ZI's plugin has its own user interface and an access to all the variables, parameters and handles. ZI utilizes the volumetric FEM \cite{demunck2012} in the forward modelling stage which is why both cortical and sub-cortical compartments can be modelled accurately without limiting their number. The SESAME  plugin combined ZI's forward simulation routine also used in EM and IAS solver with the openly available core procedure of SE\-SA\-ME\footnote{\url{https://github.com/i-am-sorri/SESAME_core}}  \cite{SorrentinoAlberto2014Bmmo,sommariva2014}. The numerical calculations were performed using Dell Precision 5820 Tower Workstation with Intel Core i9-10900X X-series Processor CPUs, Quadro RTX 4000 GPUs and 128 GB RAM. 

\subsection{Discretization and source model}

In ZI, a given head model is discretized using a tetrahedral finite element (FE) mesh. One-millimeter mesh resolution \cite{rullmann2009eeg} is applied in order to achieve a high enough accuracy with respect to the strongly folded tissue structure of the brain and thin layers of the skull. The source space is modelled via the divergence conforming approach  \cite{miinalainen2019} in which the primary current density of the brain activity is formed by a superposition of dipole-like sources belonging to the Hilbert space H(div) of vector fields with a square integrable  $L_2(\Omega)$  divergence \cite{braess2001}. The H(div) model is advantageous, since it enables the accurate modelling of both well-localized dipolar and realistic  distributed sources  and,  especially, since analytical dipoles are inapplicable (singular) as sources of the  FEM forward simulation. 

The H(div) sources are  distributed evenly  in the active compartments of the head model. The mesh-based source orientations following from the FEM discretization are interpolated into  Cartesian directions using the position based optimization approach  \cite{baumgartner2010dipole} in the case of the 10-point stencil \cite{pursiainen2016} in which the source contained by a center tetrahedron is modelled by six edge-wise and four face-intersecting sources associated with the edges and faces of that tetrahedron. Each source position is associated with the three Cartesian source orientations. The source distribution of the cerebral cortex was assumed to be parallel to its local surface normal  due to the normally oriented axons of the cortex \cite{creutzfeldt1962influence}. This normal constraint was implemented by projecting the Cartesian source field to the nearest surface normal direction in the grey matter compartment. The source orientations of the other compartments were unconstrained.

\subsubsection{Accuracy measures}
\label{sec:accuracy_measures}

The source localization accuracy with respect to a given dipolar source was investigated by comparing its position, orientation, and amplitude with the corresponding integral means obtained for the reconstructed distribution. Each integral mean was calculated in a region of interest (ROI) with 30 mm radius and  center at the source position. The following measures were evaluated: (1) the Euclidean distance and (2) angle (degree) difference, and (3) the Briggsian logarithm of the amplitude ratio between the reconstructed and actual source, i.e.,  base 10 logarithm.

\subsubsection{Focality measures}
\label{sec:focality_measures}
We consider two different focality measures described in the following. (1) Hard thresholding  is calculated by dividing the number of source positions, where the reconstruction is larger than 75 \% of its maximum within the ROI,  by the number of sources inside the ROI. (2) The earth mover's distance (EMD), in units of mm, measures the minimum amount of work that needs to be done to transfer a given mass distribution to another location and shape. The measure is originally defined as an analytical distance function between probability distributions in metric spaces in \cite{Kantorovich,Vaserstein}. The linear optimization form and the name of the distance originate from  Rubner et al. \cite{RubnerY1998Amfd}. We use EMD to compare distribution-like estimates that are spread over the whole domain to a set of dipoles at point locations. That is, we move the reconstruction mass to the set of true dipole locations in a way that the workload is minimized. In order to use EMD as a focality measure,  we add a limit condition that the reconstruction mass beyond a certain distance away from the true source location is not moved. This is necessary since, otherwise, the base level noise for which a distributional source localization estimate exists in every source location will  dictate the EMD results.  We choose this moving limit to be 45 mm to have an approximately uniform set of source points between cortical and sub-cortical true sources. SESAME's  dipole estimation accuracy is  measured  using the EMD, while the distributional measures are not applicable  \cite{lucka2012}.

\subsubsection{Spherical model for source localization experiments}
\label{sec:spherical_model}

The quantitative performance of the EM and IAS techniques and SESAME were analyzed using the isotropical Ary model, which consists of three concentric spherical compartments modelling the brain, skull, and skin. The radii of these layers are 82, 86, and 92 mm and their electrical conductivities are 0.33, 0.0042, and 0.33 S/m, respectively. A spherical model is used to minimize the effect of tissue structure on the source localization accuracy. The Ary model was discretized using one-millimeter accuracy and the source space within the brain compartment consisted of 10000 source positions. This relatively large source space size was selected as it gives an appropriate forward model accuracy and allows running each inverse method examined in this study in a few minutes. While the distributional EM, IAS methods would allow a greater source count without slowing down significantly, which was observed to be the case for the dipole search of SESAME. 

\subsubsection{MRI-based model for source localization experiments}
\label{sec:realistic}

In the qualitative analysis, we used a multi-compart\-ment segmentation generated using open T1-weighted MRI data obtained from a healthy subject. Using this data, a surface segmentation was generated using the Free\-Surfer software suite \footnote{\url{https://surfer.nmr.mgh.harvard.edu/}}.  The number of individual source positions was selected to be 100,000. The tissue conductivities suggested in \cite{dannhauer2010} were applied, i.e., 0.14, 0.33, 0.0064, 1.79, and 0.33 S/m for the white and  grey matter, skull, cerebrospinal fluid (CSF), and skin, respectively. In addition to these compartments, a set of sub-cortical compartments were included in the model. The conductivity of those was assumed to be 0.33 S/m. 

     \begin{figure}[t!]
         \centering
         \includegraphics[height=4.5cm]{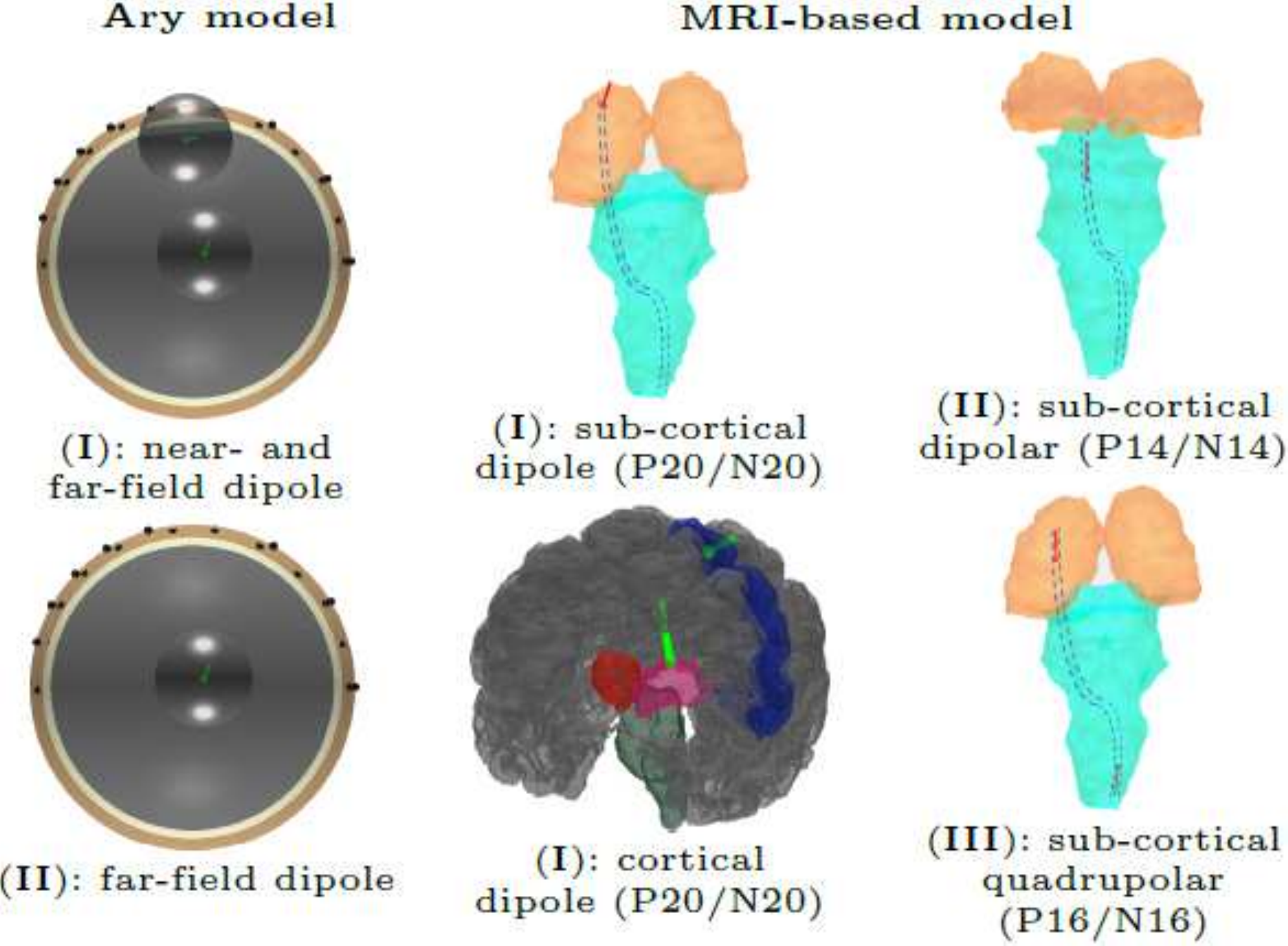}
         \caption{Visualization of near- (cortical) and far-field (sub-cortical) source location  in  the  spherical Ary  (1st from the left) and MRI-based head model (2nd and 3rd from the left). Each source position and orientation is depicted by a pointer (green and red) which, in the case of Ary model, is surrounded by a sphere showing the extent of the corresponding ROI.  Configuration ({\bf I}) corresponding to P20/N20 component of the median nerve SEP includes a superficial and deep dipole which in the MRI-based model are located in the ventral posterolateral (VPL) \cite{gotz2014thalamocortical,haueisen2007identifying} thalamus (top row, middle) and in the Brodmann area 3b \cite{buchner1995somatotopy,hari2017meg} of the primary somatosensory cortex (bottom row, blue area on the middle picture), respectively. The cortical source inherits its normal orientation with respect to the cortical surface (here the white matter surface) from that of the cortical neurons. The thalamic source is oriented along the dorsal column–medial lemniscus pathway  (blue dash)  which is a bundle of basically vertical neuron fibers conducting the SEP from the median nerve through the brainstem and thalamus to the primary somatosensory cortex. A single upward-pointing deep dipole constitutes configuration ({\bf II}) based on the P14/N14 component \cite{noel1996origin}. In the MRI-based model, it is located in the upper part of the brainstem, especially, in  the  medial lemniscus pathway \cite{noel1996origin}. Configuration  ({\bf III}) is quadrupolar, i.e.,  a combination of two oppositely  oriented dipoles, a ventrolateral  thalamic dipole with an upward orientation and a dipole in  the cuneate nucleus with an opposite orientation, creating the positive and negative pole of P16/N16   \cite{hsieh1995interaction,buchner1995origin}, respectively. }
         \label{fig:source_configurations}
           \label{fig:P14_Sourcepos}
     \end{figure}
     
\subsubsection{Ary model experiments and iteration numbers}
To analyze the accuracy and focality of the EM and IAS source localization estimates, we generated a sample of 100 different noise vector realizations utilizing the model of zero mean Gaussian noise with 5 \%  standard deviation relative to the largest absolute value in the simulated noiseless data vector. It is mathematically written as $\max \left|y_j\right|$ for $j=1,...,m$. Applying these noise realizations, a sample of EM and IAS estimates was obtained for the source configuration ({\bf I})---({\bf III}) evaluating the accuracy and focality measures for each estimate. To study the sensitivity of the methods to noise in more detail, we performed another experiment, where the reconstruction accuracy measures were  calculated for seven different noise level standard deviations: 3, 5, 7, 9 11, 13, and 15 percents (signal-to-noise ratio (SNR) 30, 26, 23, 19, 18, and 16 decibels, respectively) considering the data entry with the largest amplitude. A total of 25 estimates were evaluated for every noise level.

\begin{figure}[h!]
    \centering
    \includegraphics[width=5.5cm]{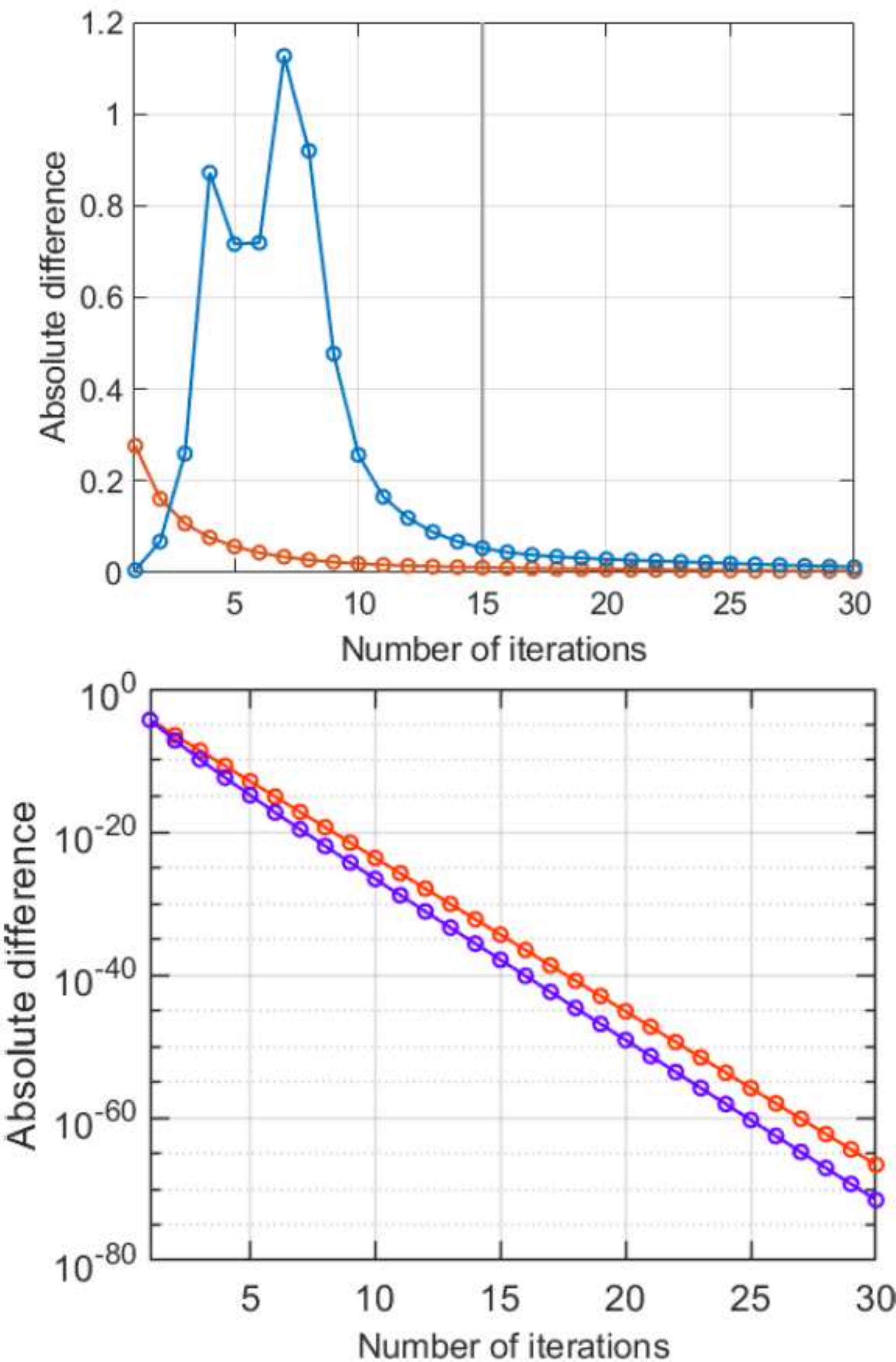}
    \caption{Convergence curves for IAS (red and orange) and EM (blue and violet) algorithms calculated as $\ell_2$-norm of difference of reconstruction vectors of sequential iteration steps from spherical Ary model ({\bf top}) and MRI-based head model ({\bf bottom}) with configuration ({\bf I}) starting by a zero vector. We find 15 iteration steps  appropriate because both  algorithms are then stabilized. The calculation time of one iteration step is is 12.5 s.}
    \label{fig:L1covergence}
\end{figure}

Following the findings of  \cite{rezaei2020randomized,Calvetti2009}, the number of iterations applied to find a reconstruction was chosen to be 10.  The number of the Lasso fixed-point iteration steps applied in finding $\hat{\bf x}^{(j+1)}$ in (\ref{eq:EMminimization}) and (\ref{eq:IAS}) with $q = 1$ was chosen to be 15 based on the convergence of $\ell_2$-norm of the reconstruction vector, see Figure \ref{fig:L1covergence}.

The number of Monte Carlo samples is an iteration-like parameter for SESAME that can be adjusted. The only limiting factor for sample size is time and its effects on reconstruction accuracy have always uncertainty in it. We decided to use (a) a default value 100 samples and (b) 700 samples that falls between the calculation times of CEP for both prior degrees. With CEP, the calculation time for degree 2 with 10 EM/IAS iterations is $18.20\pm 0.01$ s and for degree 1 with 10 EM/IAS iterations and 15 Lasso iterations it is 3 min $1.4\pm 0.1$ s. SESAME does not have constant iteration cycles because resampling is a random process and depends on the measurement data. The calculation time with 700 samples is 1 min $20.5\pm 0.5$ s in average. SESAME was also  tried with 1000 samples, but the results did not differ from the case of 700 samples. Longer processes were omitted as too time consuming compared to EM/IAS.

\subsection{Physiology-based parameter choice}

\label{sec:parameter_choice}
 To select the shape and scale parameter of the hyperprior optimally, we choose the expectation $\mathbb{E}_{\pi(x_i \, ; \, \kappa, \theta, q)[x_i]}$ of the marginal prior (\ref{eq:t-distribution}) as suggested in \cite{rezaei2020parametrizing}, i.e.,  so that it corresponds to the expected entry-wise deviation  of the reconstruction vector, to match  the random fluctuations predicted by the hyperprior approximately with the noise-induced fluctuations of the reconstruction. Consequently, the actual brain activity to be found will appropriately correspond to  the tail part of the hyperprior, i.e., the brain activity constitutes a data outlier compared to the measurement noise.  Due to the linear forward model, the relative noise level of the reconstruction may be assumed to be roughly that of the measurement noise. Thus, a random fluctuation may be assumed to have an amplitude of the relative measurement noise standard deviation (here 3 oe 5 \%) multiplied by a typical dipolar primary current amplitude in the brain, e.g., 10 nAm (1E-8 Am) \cite{hamalainen1993}. By taking into account that the largest absolute value in the measurement data is set to one microvolt which is a typical  EEG measurement amplitude \cite{niedermeyer2004} and the leadfield matrix is presented in SI-units, we can conclude the relative noise level of reconstruction to be $\alpha \cdot \hbox{1E-8} \cdot \hbox{1E6}$ in micro units, where $\alpha$ is the standard deviation of the  Gaussian noise associated with the likelihood. If we choose 3 \% as the likelihood, it follows that the expected deviation is $\hbox{3E-4}$. The shape parameter is chosen to be $\kappa = 4.4$ for which the ratio between $\mathbb{E}_{\pi(x_i \, ; \, \kappa, \theta, q)}[x_i]$ and $\theta^{1/q}$ is equal for both prior degrees $q = 1$ and $q =2$, i.e., $ \left. \frac{\mathbb{E}_{\pi(x_i \, ; \, \kappa, \theta, q)}[x_i]}{\theta^{1/q}} \right|_{q = 1} = \left. \frac{\mathbb{E}_{\pi(x_i \, ; \, \kappa, \theta, q)}[x_i]}{\theta^{1/q}} \right|_{q = 2} \approx 0.3$. In this way, we try to reduce the effect of the parameters in  the comparison between prior degrees $q$.
Consequently, the value of $\mathbb{E}_{\pi(x_i \, ; \, \kappa, \theta, q)}[x_i]$ follows from the scale parameter which is set to be $\theta = \hbox{1E-3}$ and $\theta = \hbox{1E-6}$ for $q = 1$ and $q =2$, respectively, in order to obtain the correspondence to the noise as described above. 

\begin{figure}
    \centering
    \begin{minipage}{8cm} \centering
    \includegraphics[width=5cm]{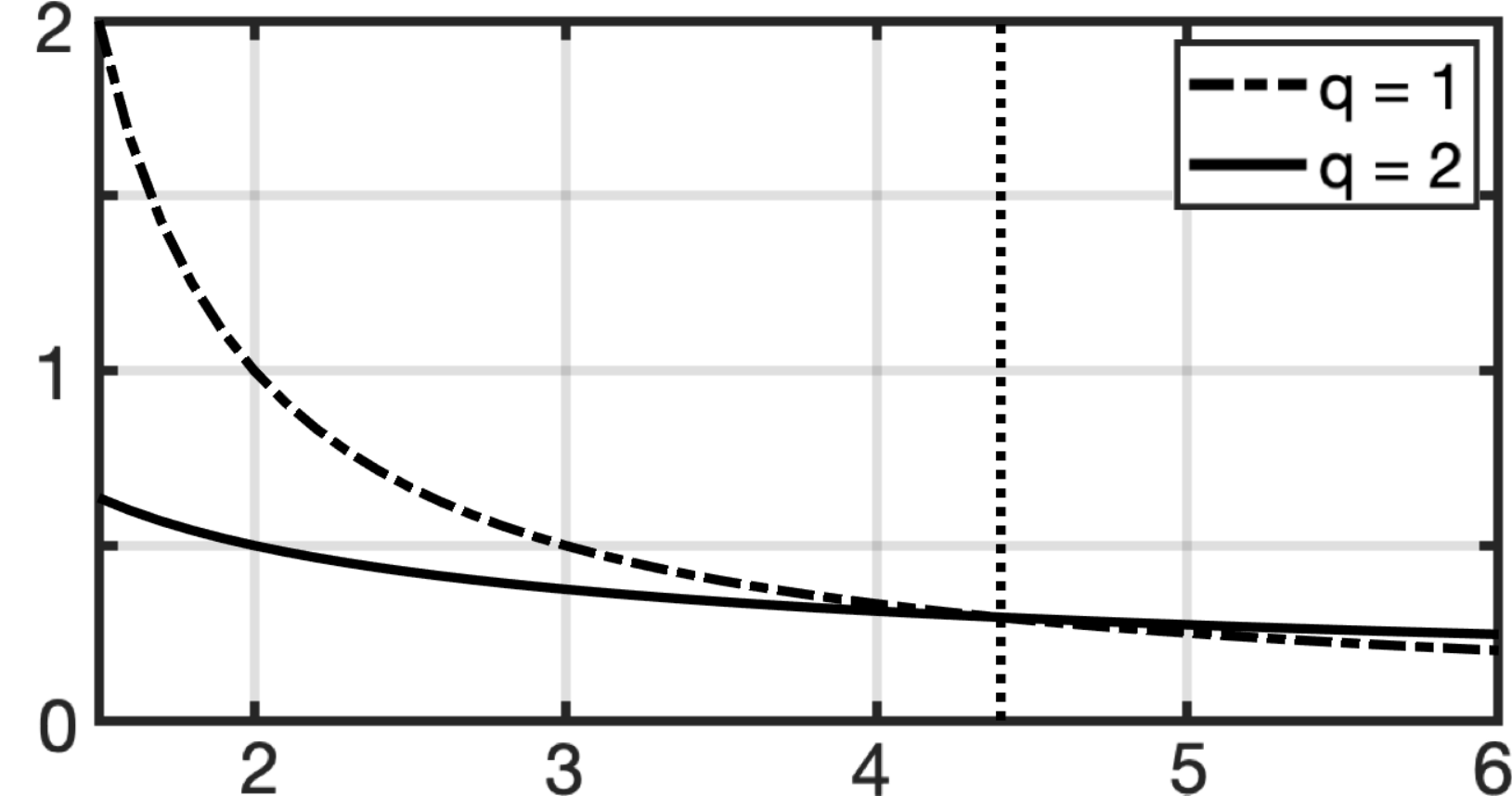} \\ \vskip0.2cm
    Expectation $\mathbb{E}_{\pi(x_i \, ; \, \kappa, \theta, q)} [x_i]/\theta^{1/q}$
    \end{minipage} \\
    \begin{minipage}{8cm}
    \centering
    \includegraphics[width=5cm]{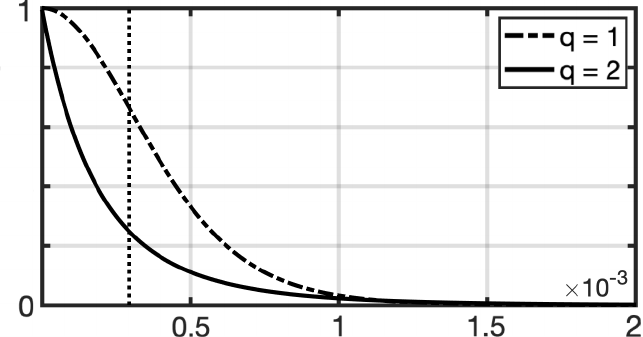} \\ \vskip0.2cm
    Marginal prior $\pi(x_i \, ; \, \kappa, \theta, q)$ with $\kappa = 4.4$
    \end{minipage}
    \caption{{\bf Left:} The ratio between the expectation $\mathbb{E}_{\pi(x_i \, ; \, \kappa, \theta, q)} [x_i]$ of the marginal  prior and $\theta^{1/q}$ as a function of the shape parameter value $\kappa$ from $\kappa = 1.5$ to $\kappa = 6$. The vertical line corresponds to the value  $\kappa = 4.4$ for which the cases $q = 1$ and $q = 2$ match. {\bf Right:} The marginal prior (normalized to one) with the expectation set to 3E-4, approximating typical deviations due to noise. The shape parameter is $\kappa = 4.4$ and the scale parameter $\theta = \hbox{1E-3}$, when $q = 1$, and $\theta = \hbox{1E-6}$, when $q =2$, resulting into $\mathbb{E}_{\pi(x_i \, ; \, \kappa, \theta, q)} [x_i] = \hbox{3E-4}$ in both cases. With these parameter choices the expectation can be assumed to appropriately coincide with the amplitude of the  expected random fluctuations of the reconstruction and the intensity of the actual source to be found is located in the tail part, i.e., it is an outlier with respect to the noise level (Section  \ref{sec:parameter_choice}).  }
    \label{fig:marginal_prior}
\end{figure}

\subsubsection{Synthetic data of somatosensory evoked potentials}\label{Sect:Sythetic_data}
\begin{table*}[ht!]
\caption{Collected information of the synthetic originators used in numerical experiments.}\label{tb:synthDipoles}
\centering
\begin{tabular}{|l|c|c|r|} 
 \hline
 Name & Time (ms) & Location & Type \\ \hline 
 P14/N14 & 14 & Pons & Dipolar (single dipole)\\ 
 P16/N16 & 16 & Thalamus / Cuneate nucleus & Quadrupolar (two opposite dipoles)\\ 
  P20/N20 & 20 & Brodmann area 3b / Thalamus &  Dipolar (cortical and sub-cortical dipoles)\\ 
 \hline
\end{tabular}
\end{table*}

As an example case in the numerical experiments, we consider the detection of  synthetic SEPs modelled according to the SEPs occurring in human median nerve stimulation. SEPs occur as a response to electrical pulses that stimulate the median nerve in the wrist area.  We modelled the originators of three SEP components of which the P20/N20 component originates 20 ms post-stimulus at 3b Brodmann area, and involves simultaneous sub-cortical activity at the ventral posterolateral (VPL) thalamus \cite{haueisen2007identifying,gotz2014thalamocortical,noel1996origin,RezaeiAtena2021Rsac}. P20/N20 is preceded by far-field components P14/N14 and P16/\-N16 occurring at 14 and 16 ms post-stimulus, respectively. P14/N14  originates in  the brainstem, where  the spike afferent volley travels through medial lemniscus pathway  \cite{noel1996origin} .  The quadrupolar P16/N16  includes a positive thalamic and a negative sub-thalamic  originator  \cite{hsieh1995interaction,buchner1995origin}. The first one of these is located at the ventrolateral thalamus and the second one is at the cuneate nucleus.  

The following three configurations  ({\bf I})--({\bf III}) of dipolar sources were applied to model SEP components.  Configuration ({\bf I}) consists of  two  sources modelling the simultaneous near- (cortical) and far-field (thalamic) activity corresponding to P20/N20. The amplitude of the cortical source is assumed to be 70 \% of the thalamic one, modelling a situation, where the amplitude of the cortical activity is intensified before reaching its maximum. ({\bf II}) is formed by a single source in the medial lemniscus pathway approximating the sub-cortical activity of P14/N14.  ({\bf III}) includes a quadrupolar configuration with  an upward  component in the ventral thalamus and a downward component  in the cuneate nucleus of the brainstem  \cite{hsieh1995interaction}. Since the thalamic component is likely to be peaked slightly after 16 ms \cite{hsieh1995interaction}, its amplitude is assumed to be 77 \% compared to the amplitude of the sub-thalamic component. The present quadrupolar setup is not to be mixed with a quadrupolar afferent volley which is difficult to be detected as is. Instead the two dipolar components correspond to two subsequent  quadrupolar  spikes with one of the dipolar components visible due to a local  discontinuity of the conductivity distribution \cite{buchner1995origin}. The details of the originators are collected in Table \ref{tb:synthDipoles}.
A spherical and an MRI-based model were applied in the source localization tests. Configuration ({\bf III}) was applied only to the MRI-based model, which allows distinguishing the thalamic and brainstem areas. In the first one of these, the source positions  were  selected according to   \cite{buchner1995origin,buchner1995somatotopy}   and, in the second one, they were placed in the aforementioned originator areas (Figure  \ref{fig:source_configurations}). The noise vector ${\bf e}$ in (\ref{eq:lead_field}) was assumed to be zero-mean Gaussian random variable with diagonal covariance and 5 \% relative standard deviation compared to the amplitude of the noiseless signal.

          \begin{figure*}[h!]
    \centering
    \begin{tiny}
     \vskip0.07cm \hrulefill \vskip0.07cm  
   {\bf Configuration (I): Expectation maximization} \\ \vskip0.07cm   
\centering
 \includegraphics[height=5.80cm]{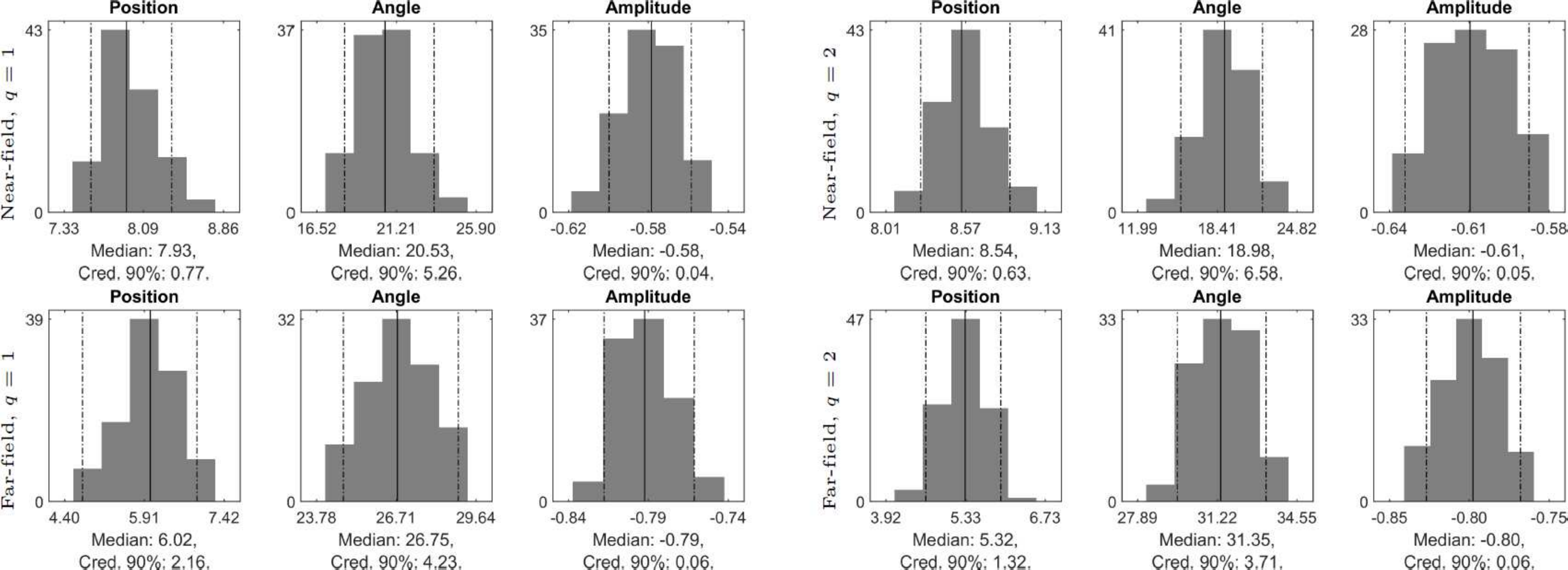} 
\\  \vskip0.07cm \hrulefill \vskip0.07cm  
{\bf Configuration (I): Iterative alternating sequential} \\ \vskip0.07cm
\centering
    \includegraphics[height=5.80cm]{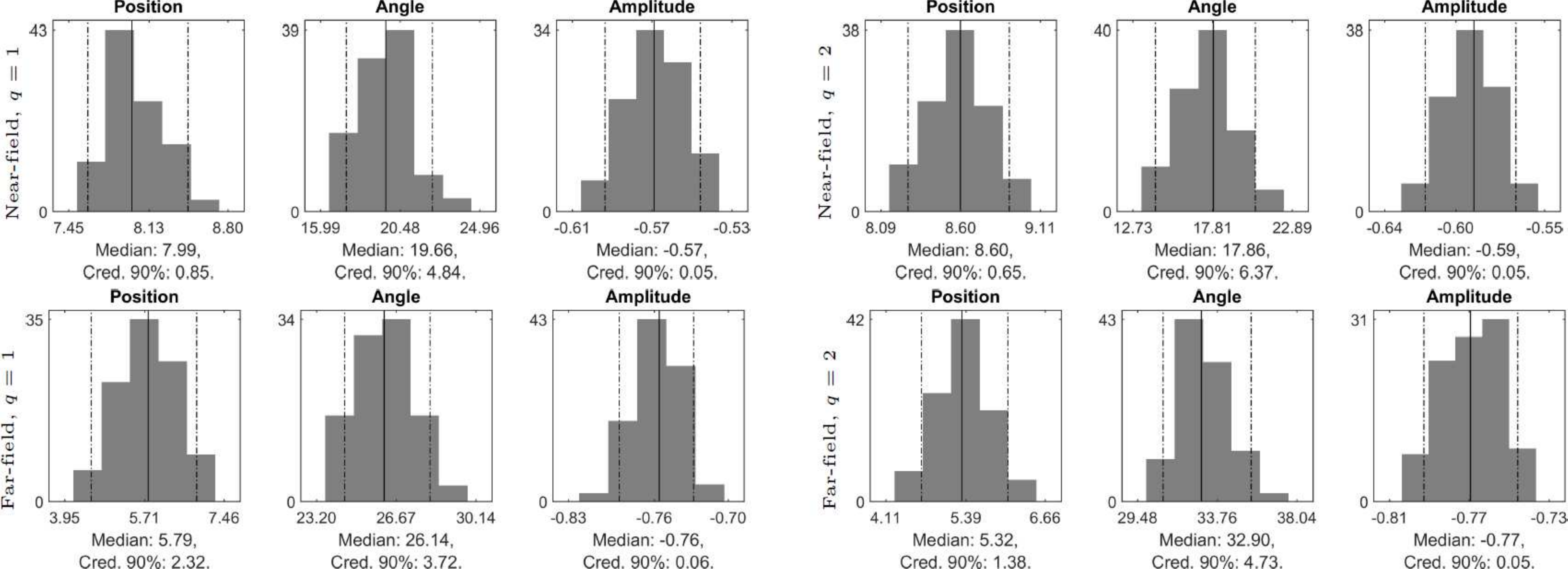} 
\\
 \vskip0.07cm \hrulefill  \vskip0.07cm  
  {\bf Configuration (II): Expectation maximization} \\ \vskip0.07cm
\centering
    \includegraphics[height=2.80cm]{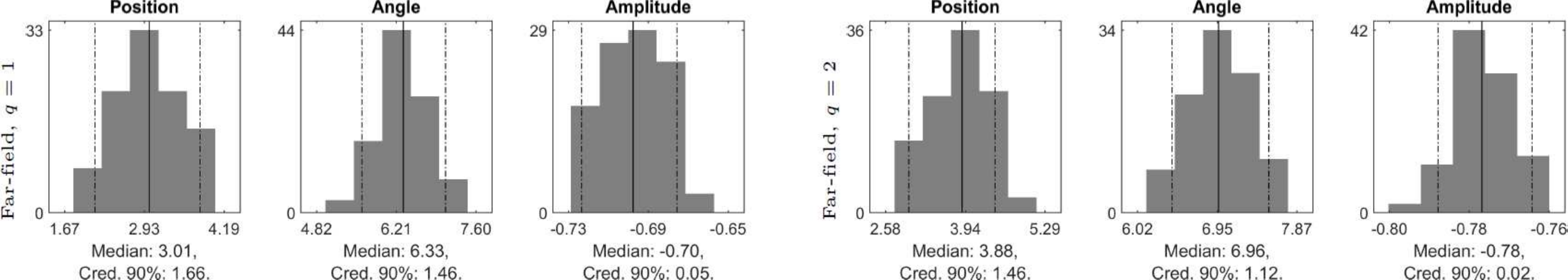}  
\\
    \vskip0.07cm \hrulefill  \vskip0.07cm  
{\bf Configuration (II): Iterative alternating sequential} \\ \vskip0.07cm
\centering
    \includegraphics[height=2.80cm]{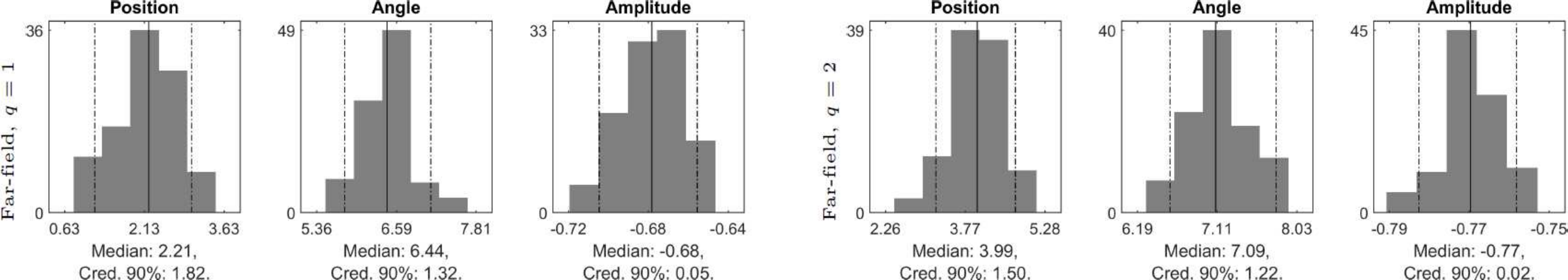} 
 \\
         \vskip0.07cm \hrulefill  \vskip0.07cm  
         \end{tiny}
    \caption{The source localization accuracy measures (position, angle, and amplitude)  evaluated in the spherical Ary model (Section \ref{sec:spherical_model}) with 5 \% noise for source configuration ({\bf I}) and ({\bf II}) from Section \ref{Sect:Sythetic_data}} applying the EM and IAS algorithm. The histograms show the results obtained for 100 different reconstructions, each corresponding to a different realization of the measurement noise. The measures concern the difference between the actual source and the mass centre of the reconstructed distribution in the corresponding ROI. The units of the position, angle and amplitude are, respectively, in mm, degrees and  $\log_{10} (A_r / A_s)$, where $A_r$ is the amplitude of the reconstructed and $A_s$ of the actual source.
    \label{fig:results_accuracy5}
\end{figure*}

\begin{figure*}
    \centering
    \begin{minipage}{0.45\textwidth}
    \centering
    \begin{tiny}
     \vskip0.07cm \rule{0.9\textwidth}{0.4pt}  \vskip0.07cm  
   {\bf Configuration (I): Expectation maximization} \\ \vskip0.07cm
\centering
    \includegraphics[height=5.70cm]{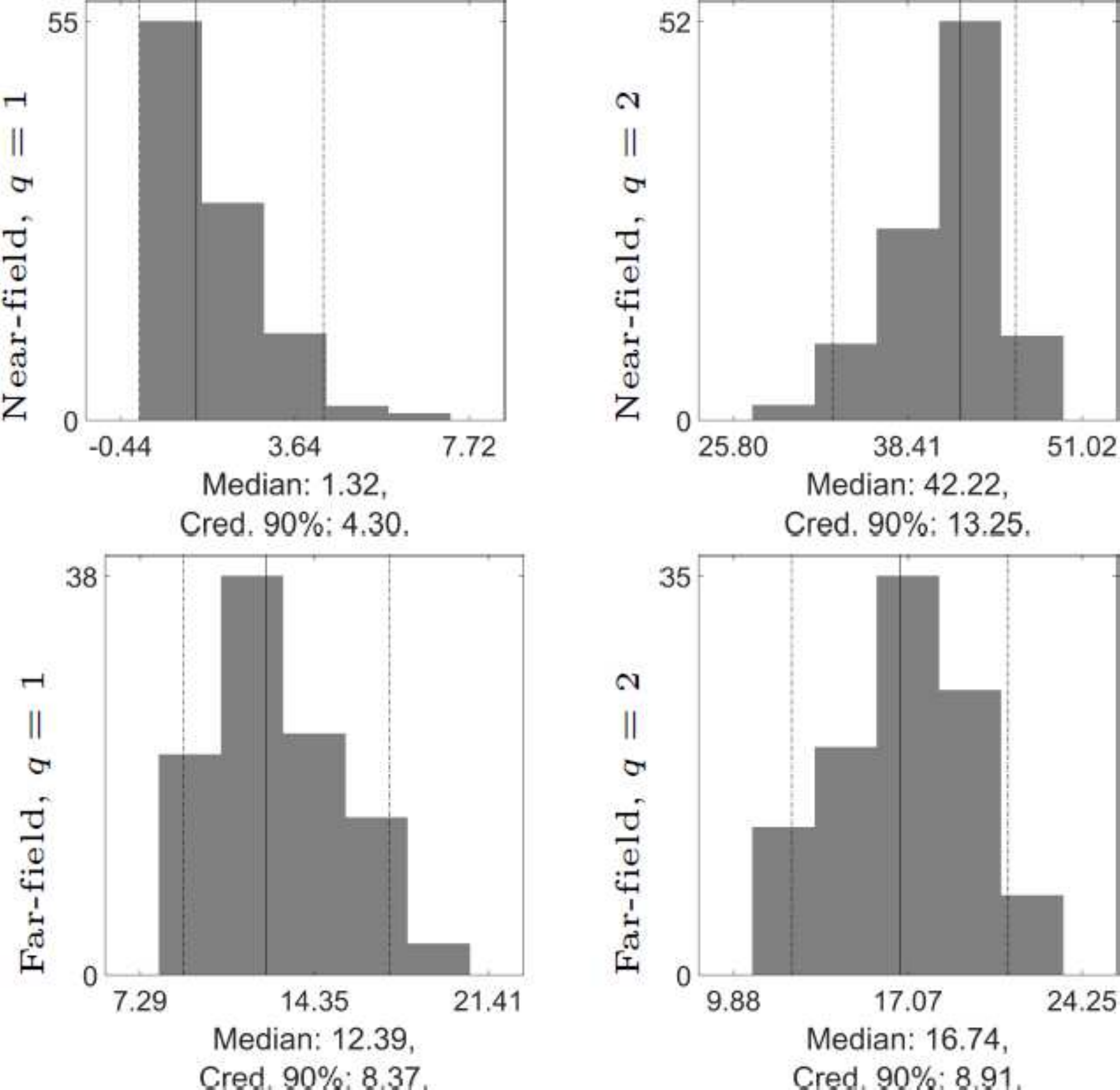} 
 \\
 \vskip0.07cm \rule{0.9\textwidth}{0.4pt} \vskip0.07cm
{\bf Configuration (I): Iterative alternating sequential} \\ \vskip0.07cm
\centering
    \includegraphics[height=5.70cm]{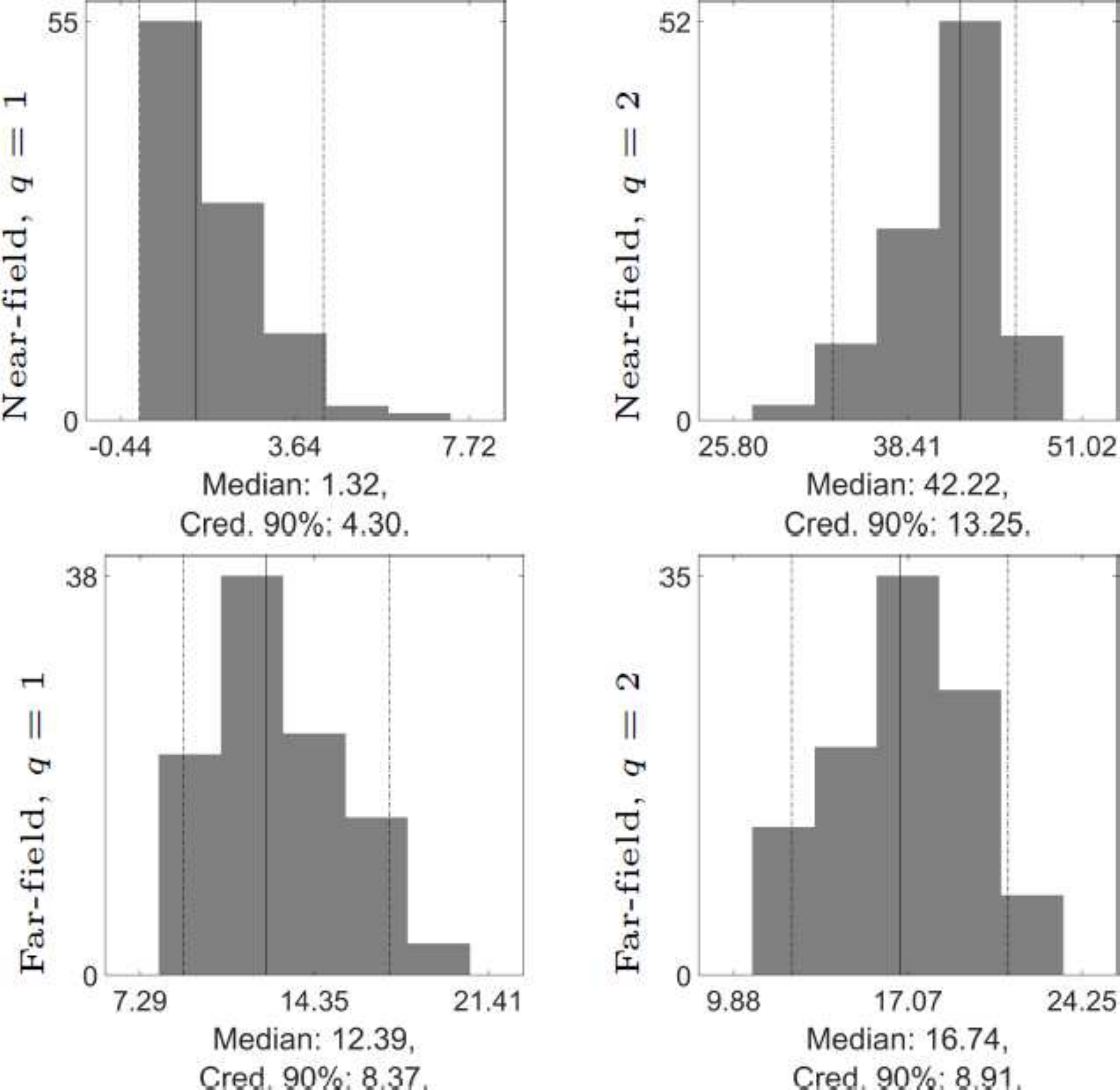}
 \\
 \vskip0.07cm \rule{0.9\textwidth}{0.4pt} \vskip0.07cm
  {\bf Configuration (II): Expectation maximization} \\ \vskip0.07cm
\centering
    \includegraphics[height=2.80cm]{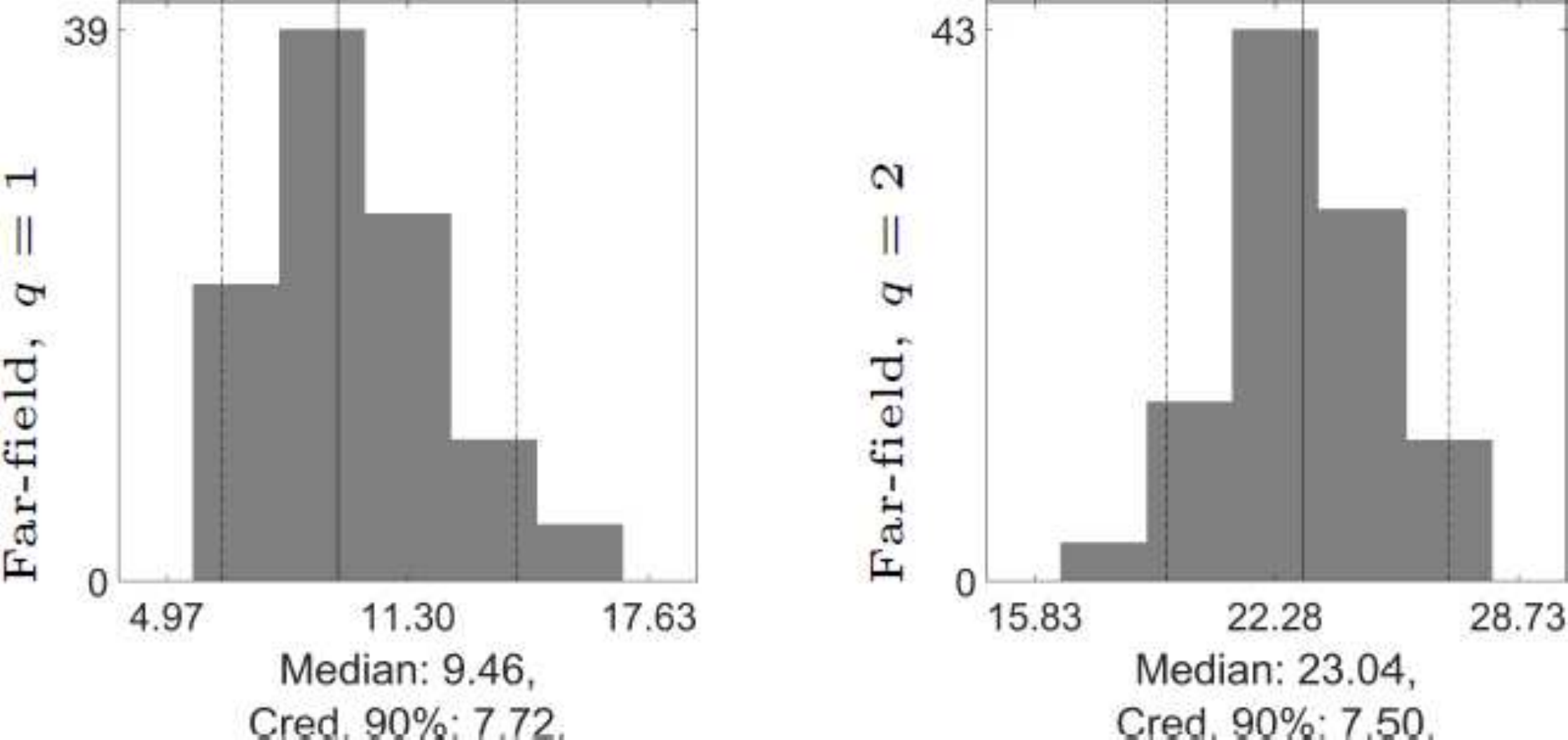} 
   \\
        \vskip0.07cm \rule{0.9\textwidth}{0.4pt} \vskip0.07cm
{\bf Configuration (II): Iterative alternating sequential} \\ \vskip0.07cm
\centering
    \includegraphics[height=2.80cm]{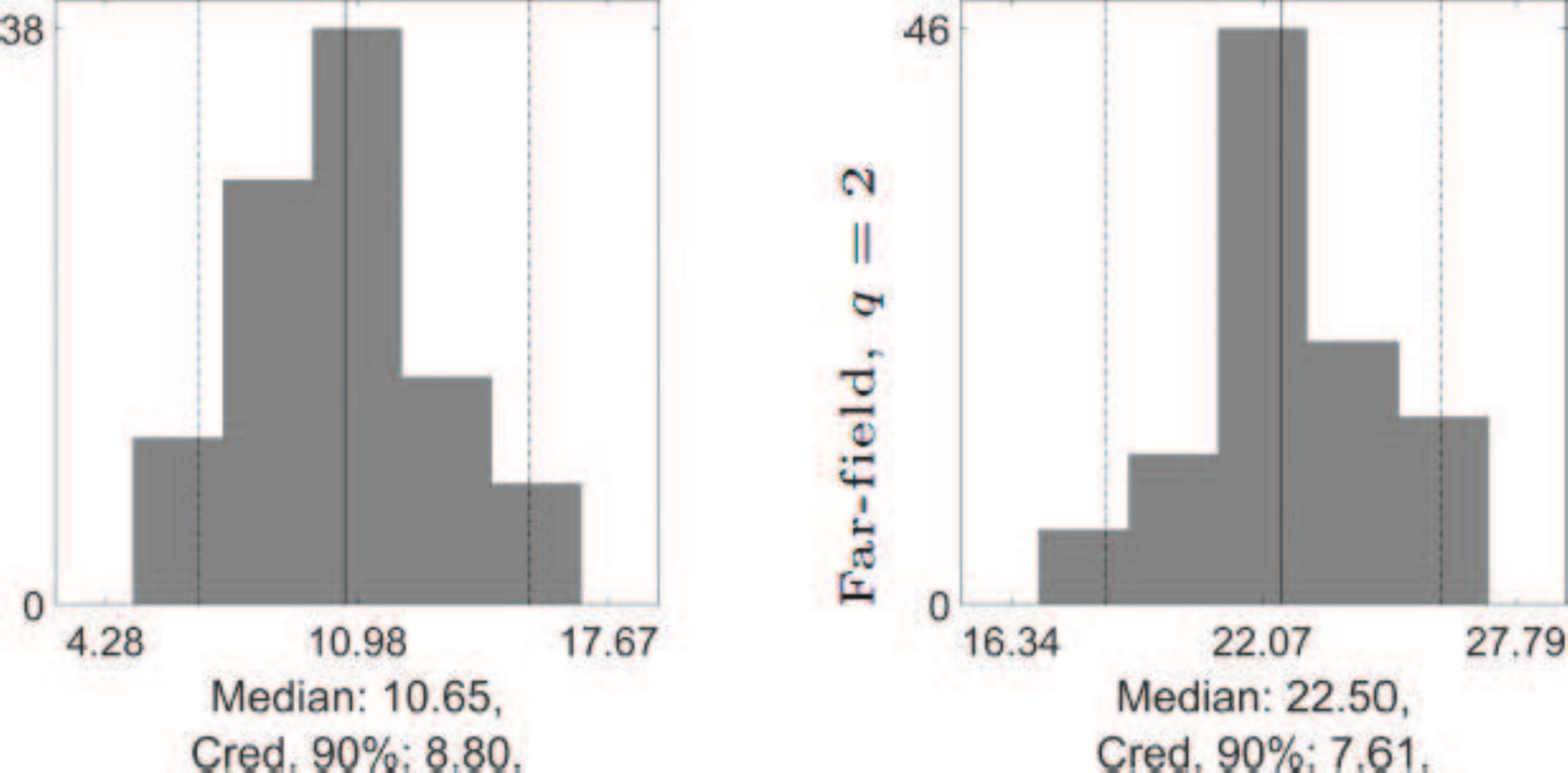} 
   \\
\vskip0.07cm \rule{0.9\textwidth}{0.4pt} \vskip0.07cm  
        \end{tiny}
    \captionof{figure}{The focality measure hard threshold evaluated in the spherical Ary model (Section \ref{sec:spherical_model}) with 5 \% noise for source configuration ({\bf I}) and ({\bf II}) applying the EM and IAS algorithm. The histograms show the results obtained for 100 different reconstructions, each corresponding to a different realization of the measurement noise. The hard threshold measure concerns the area where the intensity of the reconstruction is at least 75 \% of its maximum in the corresponding ROI.}
    \label{fig:results_focality5}
\end{minipage}\hskip0.5cm\begin{minipage}{0.5\textwidth}
    \centering
    \begin{tiny}
     \vskip0.07cm \rule{0.9\textwidth}{0.4pt}  \vskip0.07cm  
   {\bf Configuration (I): Expectation maximization} \\ \vskip0.07cm
\centering
    \includegraphics[height=2.80cm]{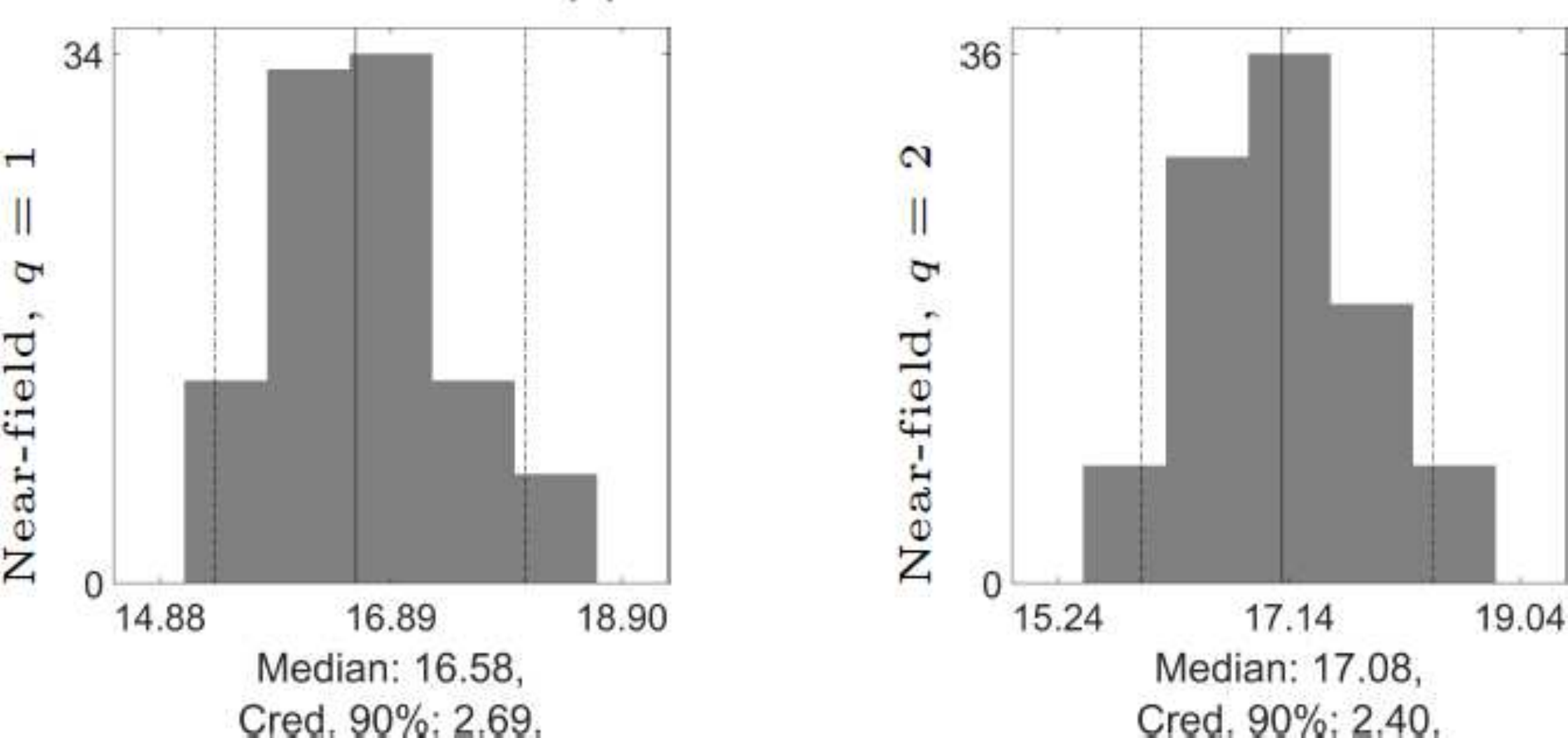} 
\\
 \vskip0.07cm \rule{0.9\textwidth}{0.4pt} \vskip0.07cm
{\bf Configuration (I): Iterative alternating sequential} \\ \vskip0.07cm
\centering
    \includegraphics[height=2.80cm]{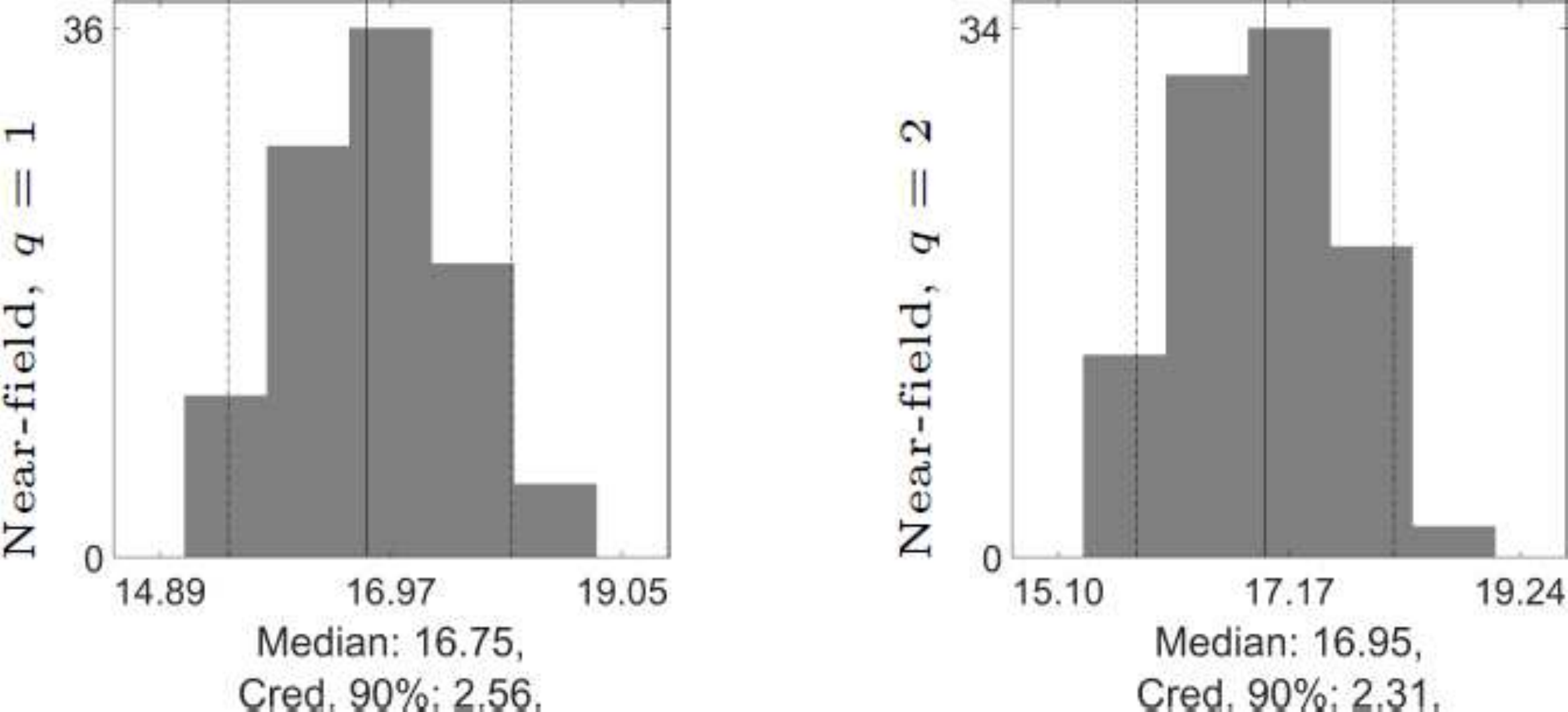} 
 \\
       \vskip0.07cm \rule{0.9\textwidth}{0.4pt} \vskip0.07cm
{\bf Configuration (I): SESAME} \\ \vskip0.07cm
\centering
    \includegraphics[height=2.80cm]{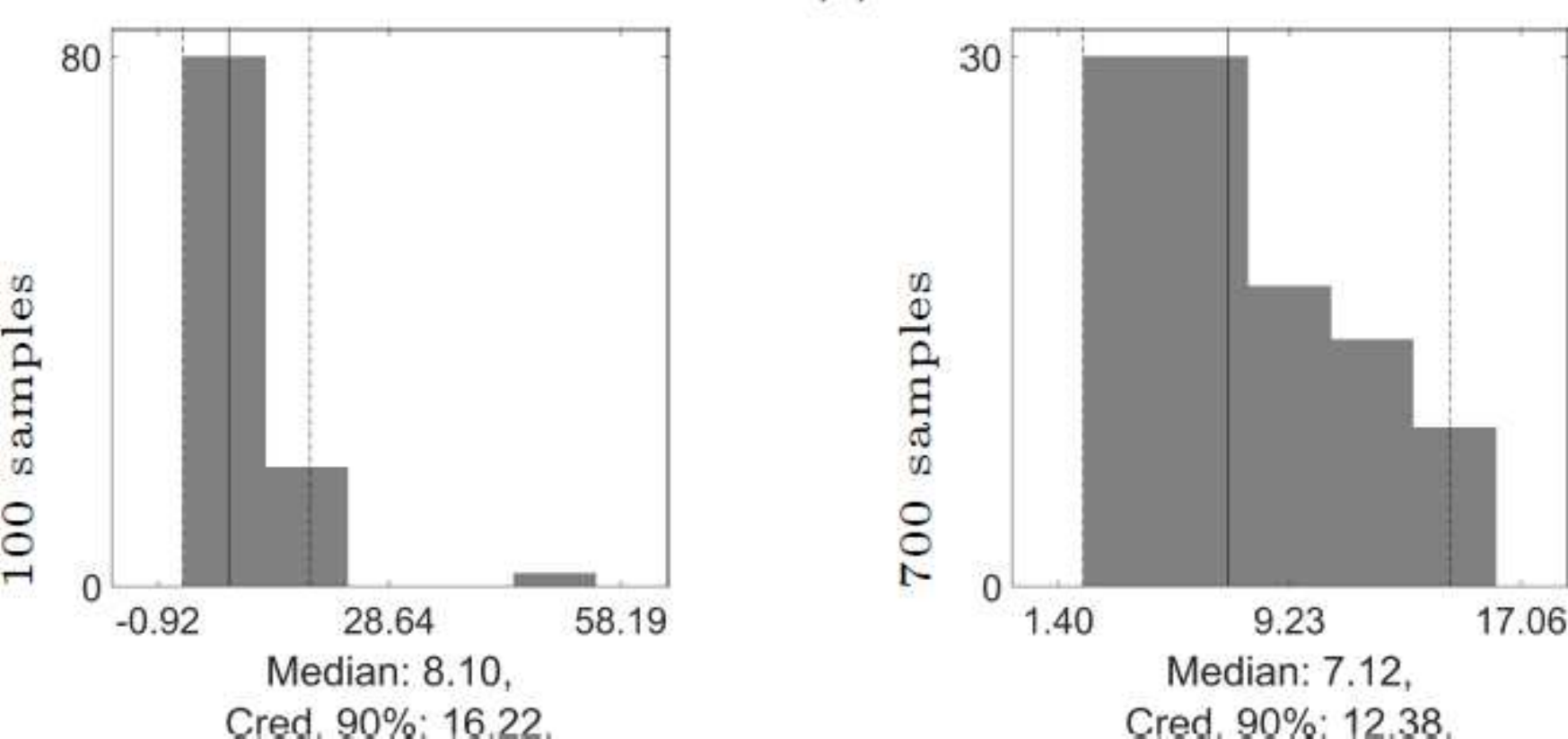} 
 \\
 \vskip0.07cm \rule{0.9\textwidth}{0.4pt} \vskip0.07cm
  {\bf Configuration (II): Expectation maximization} \\ \vskip0.07cm
\centering
    \includegraphics[height=2.80cm]{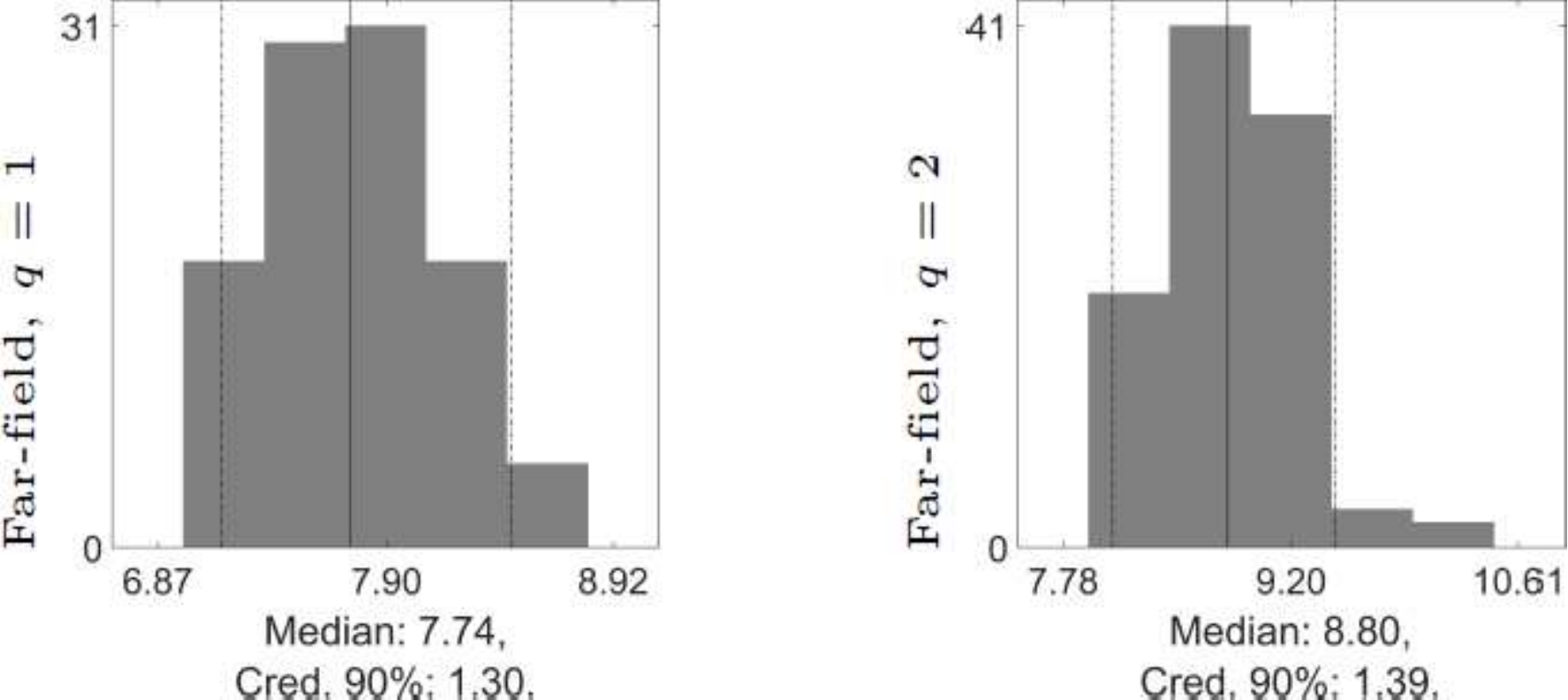}
 \\
        \vskip0.07cm \rule{0.9\textwidth}{0.4pt} \vskip0.07cm
{\bf Configuration (II): Iterative alternating sequential} \\ \vskip0.07cm
\centering
    \includegraphics[height=2.80cm]{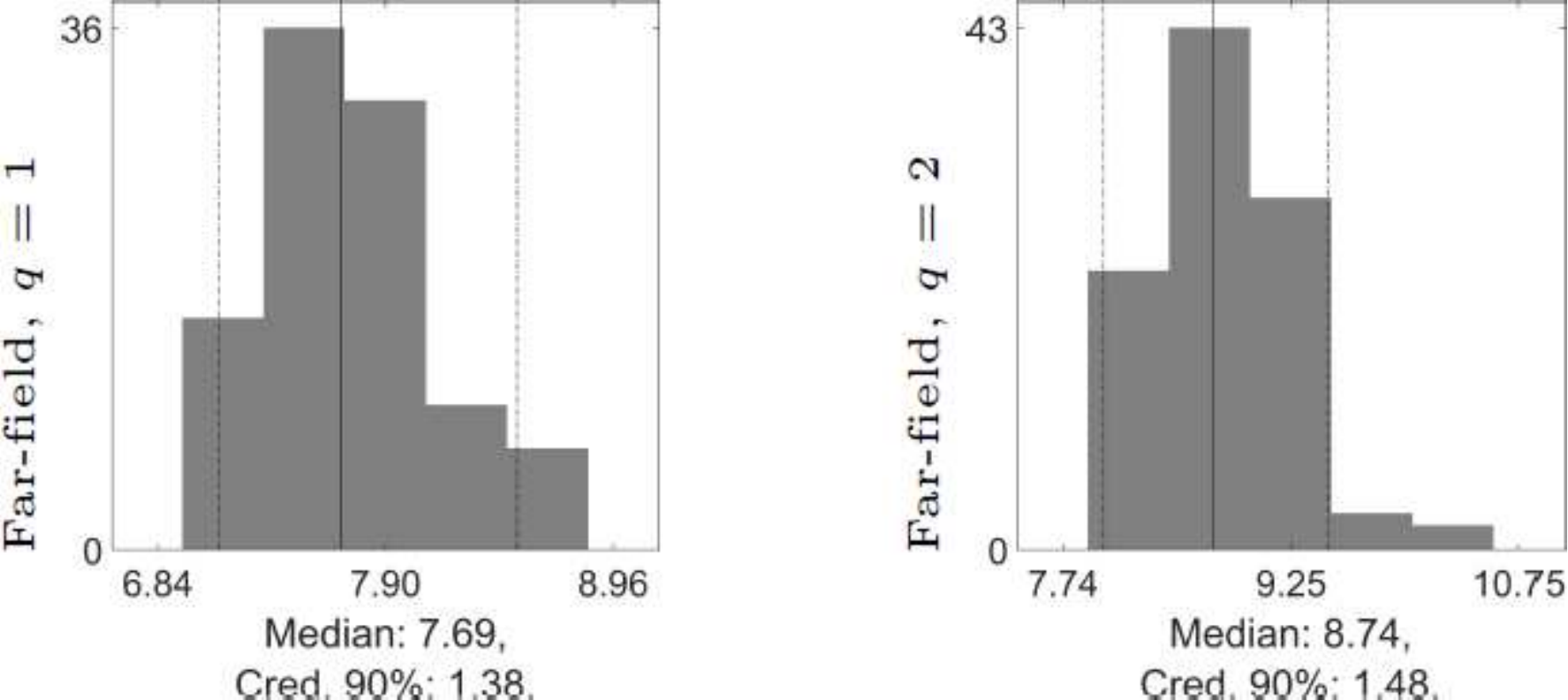} 
\\
       \vskip0.07cm \rule{0.9\textwidth}{0.4pt} \vskip0.07cm
{\bf Configuration (II): SESAME} \\ \vskip0.07cm
\centering
    \includegraphics[height=2.80cm]{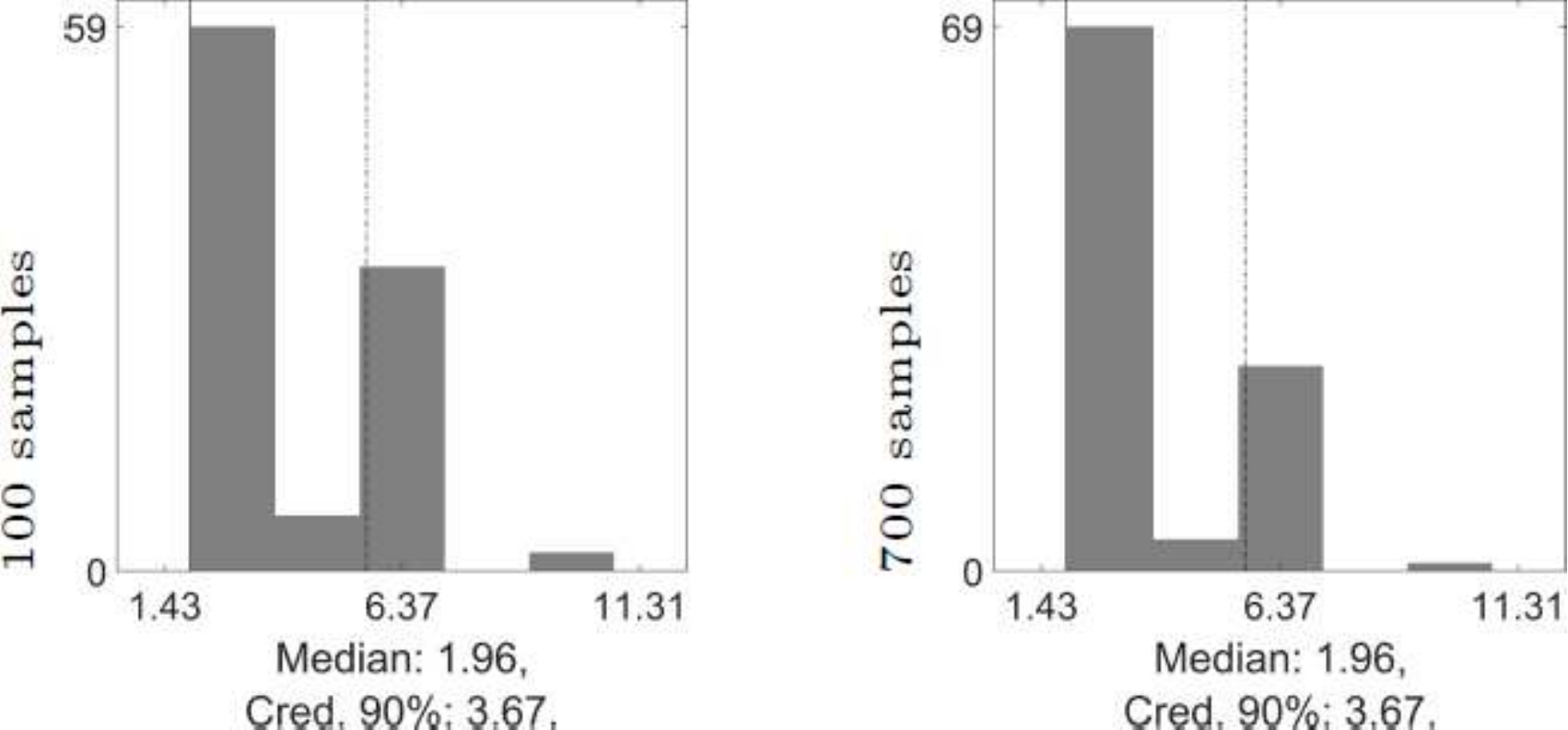} 
\\
 \vskip0.07cm \rule{0.9\textwidth}{0.4pt} \vskip0.07cm
        \end{tiny}
    \captionof{figure}{The earth mover's distance (EMD) with 45 mm moving limit evaluated in the spherical Ary model (Section \ref{sec:spherical_model}) with 5 \% noise for source configuration ({\bf I}) and ({\bf II}) applying the EM and IAS algorithm.}
    \label{fig:results_EarthMover5}
    \end{minipage}
\end{figure*}

\begin{figure*}[h!]
    \centering
    \begin{tiny}
     \vskip0.07cm \hrulefill \vskip0.07cm  
   {\bf Configuration (I): SESAME} \\ \vskip0.07cm   
\centering
 \includegraphics[height=5.8cm]{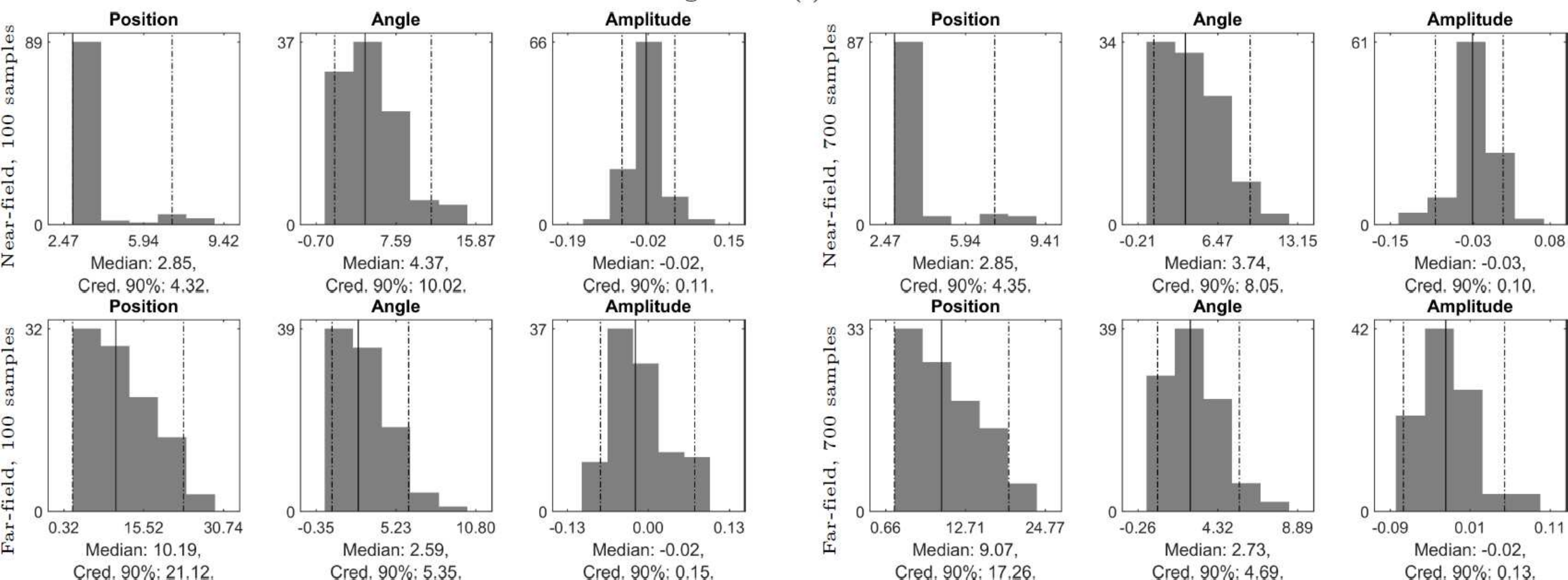} 
 \vskip0.07cm \hrulefill \vskip0.07cm  
  {\bf Configuration (II): SESAME} \\ \vskip0.07cm
\centering
    \includegraphics[height=2.8cm]{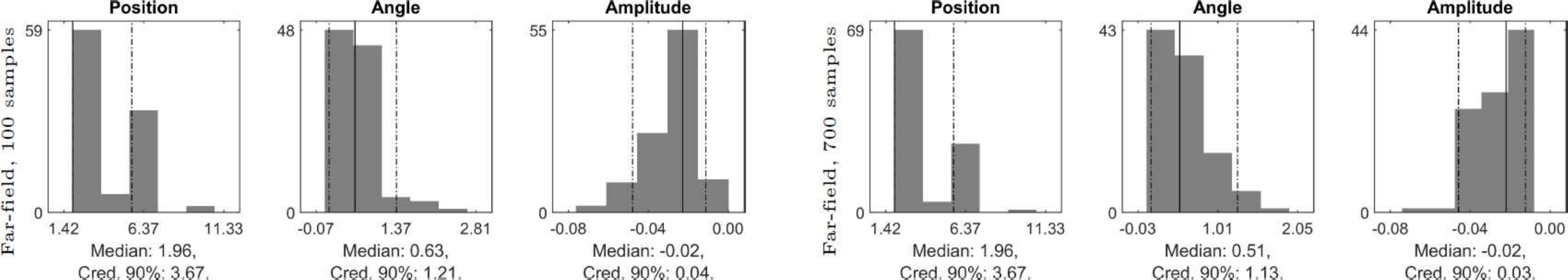}  
\\
    \vskip0.07cm \hrulefill  \vskip0.07cm
        \end{tiny}
    \caption{The source localization accuracy measures  (position, angle, amplitude, and maximum)  evaluated in the spherical Ary model (Section \ref{sec:spherical_model}) with 5 \% noise for source configuration ({\bf I}) and ({\bf II}) applying the SESAME algorithm with 100 and 700 samples. Results are obtained with the same scheme as in the case of CEP.}
    \label{fig:results_SESAME_accuracy}
\end{figure*}
     
      \begin{figure}[h!]
    \centering
    \begin{tiny}

     {\bf  Configuration (I): Expectation maximization} \\ \vskip0.07cm
     
   \centering
    \includegraphics[height=3.55cm]{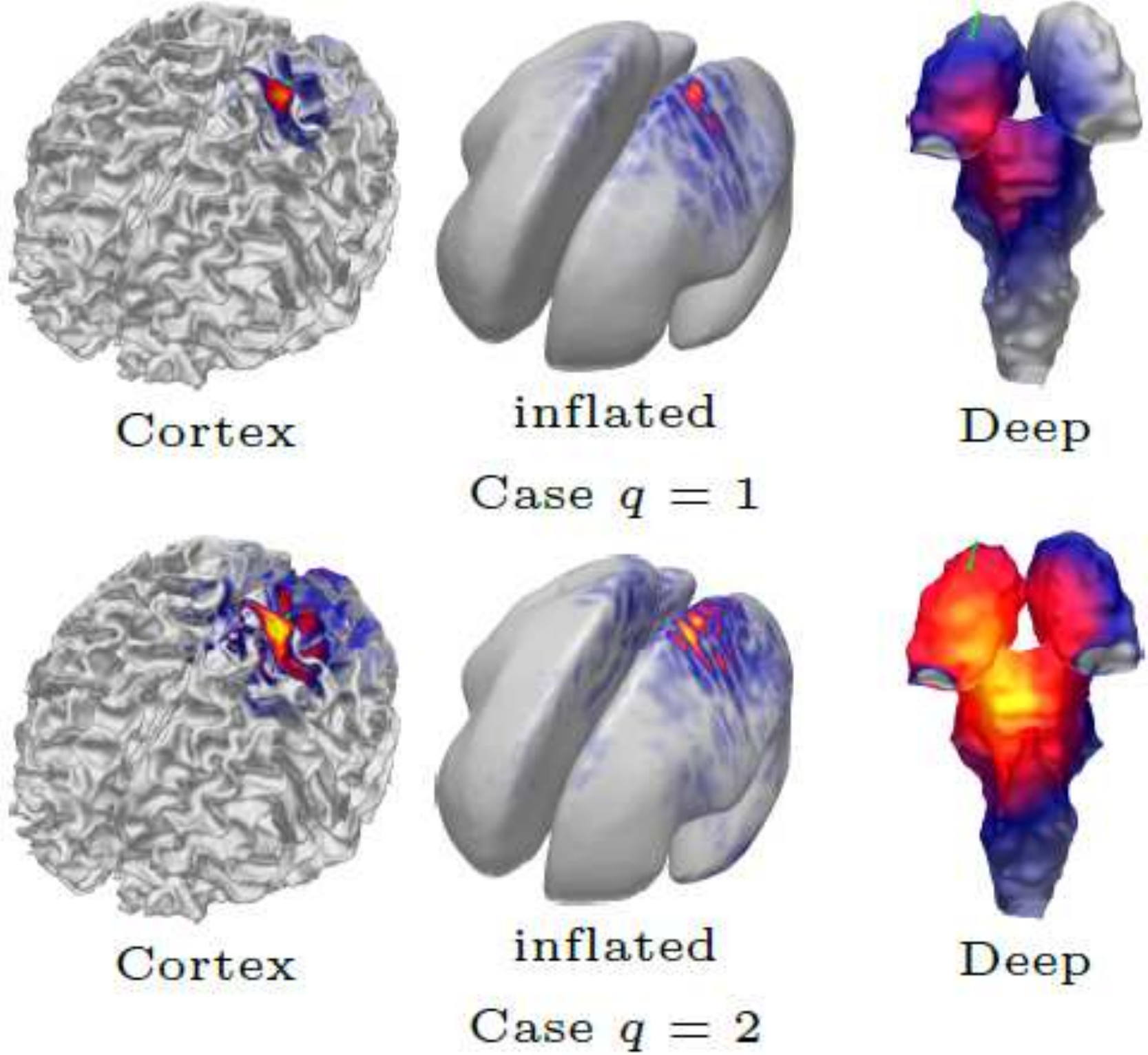}      \\ \vskip0.07cm         
      {\bf Configuration (I): Iterative alternating sequential } \\ \vskip0.07cm
   \centering
    \includegraphics[height=3.55cm]{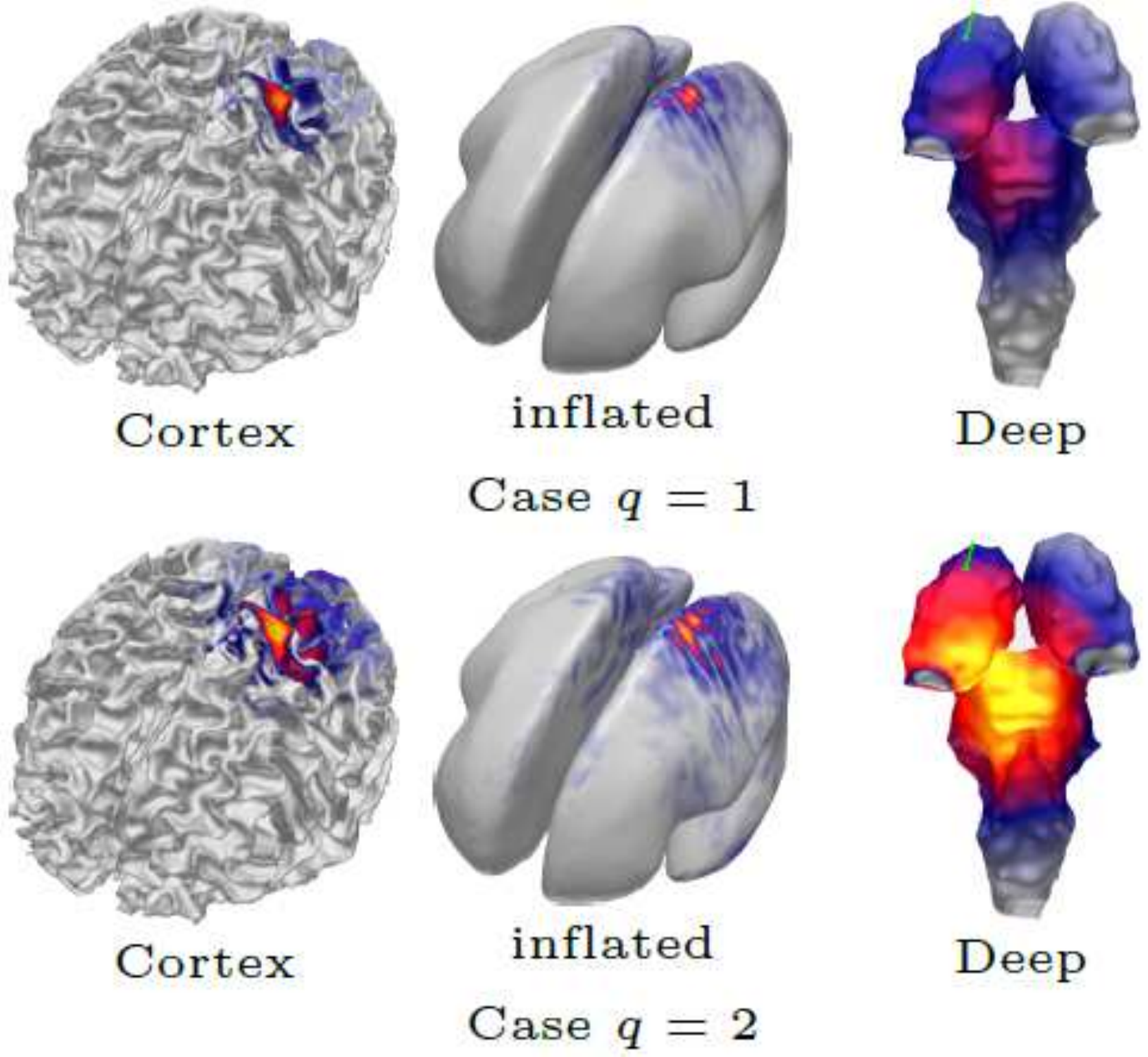} \\ \vskip0.07cm
     {\bf Configuration (I): SESAME } \\ \vskip0.07cm
      
   \centering
    \includegraphics[height=1.5cm]{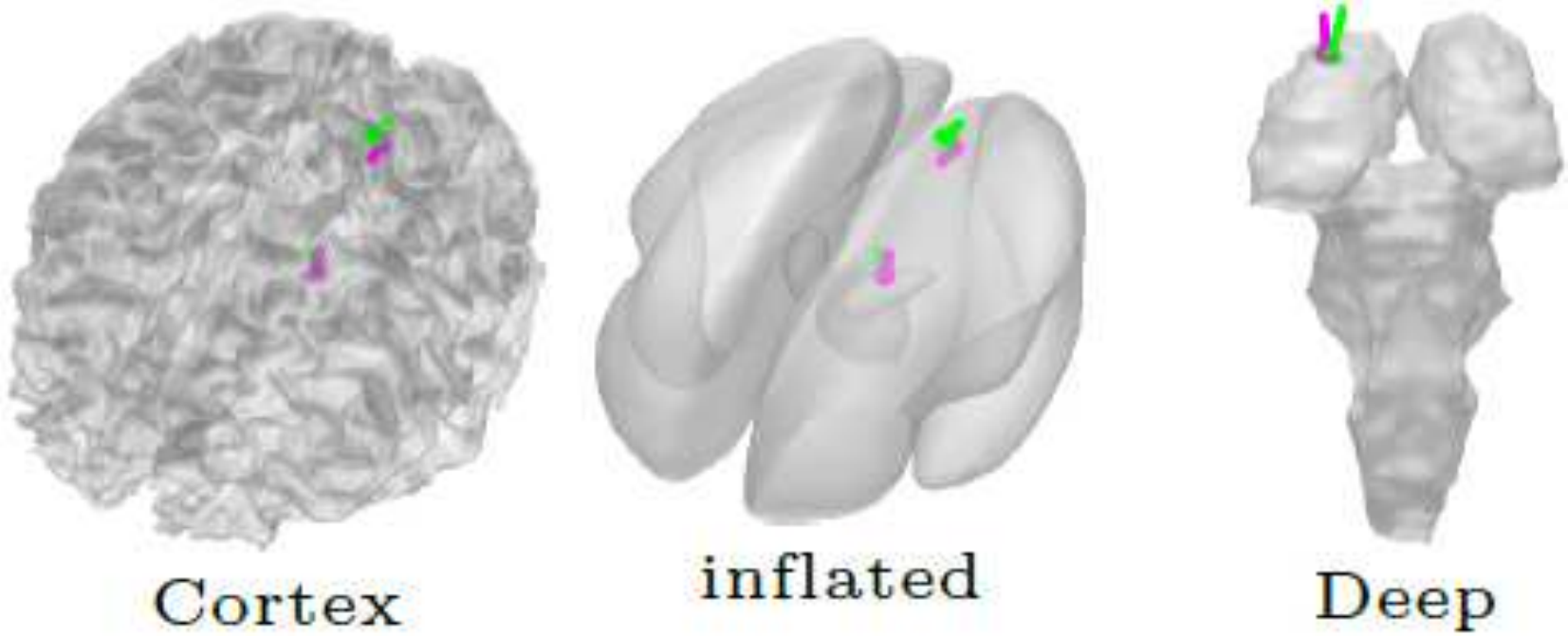}\\ \vskip0.07cm
    \end{tiny}
   
    \caption{The reconstructions obtained with the MRI-based head multi-compartment model with 3 \% noise for source  configuration ({\bf I})  including two sources, one placed in the left 3b Brodmann area of the central sulcus, pointing inwards in the direction of the local surface normal vector,  and a vertical source placed in the ventral posterolateral part of the left thalamus (Figure \ref{fig:source_configurations}).  The actual source position is visualized by a green pointer in each image. The first-degree conditionally exponential prior (CEP), i.e., the case $q = 1$, leads to a more focal reconstruction than the second-degree CEP ($q = 2$). The dipole estimations of SESAME are indicated by magenta pointers.}
    \label{fig:mri_based_1}
     \end{figure}
       \begin{figure}[h!]
    \centering
    \begin{tiny} 
    {\bf  Configuration (II): Expectation maximization} \\ \vskip0.07cm
     \centering
    \includegraphics[height=3.55cm]{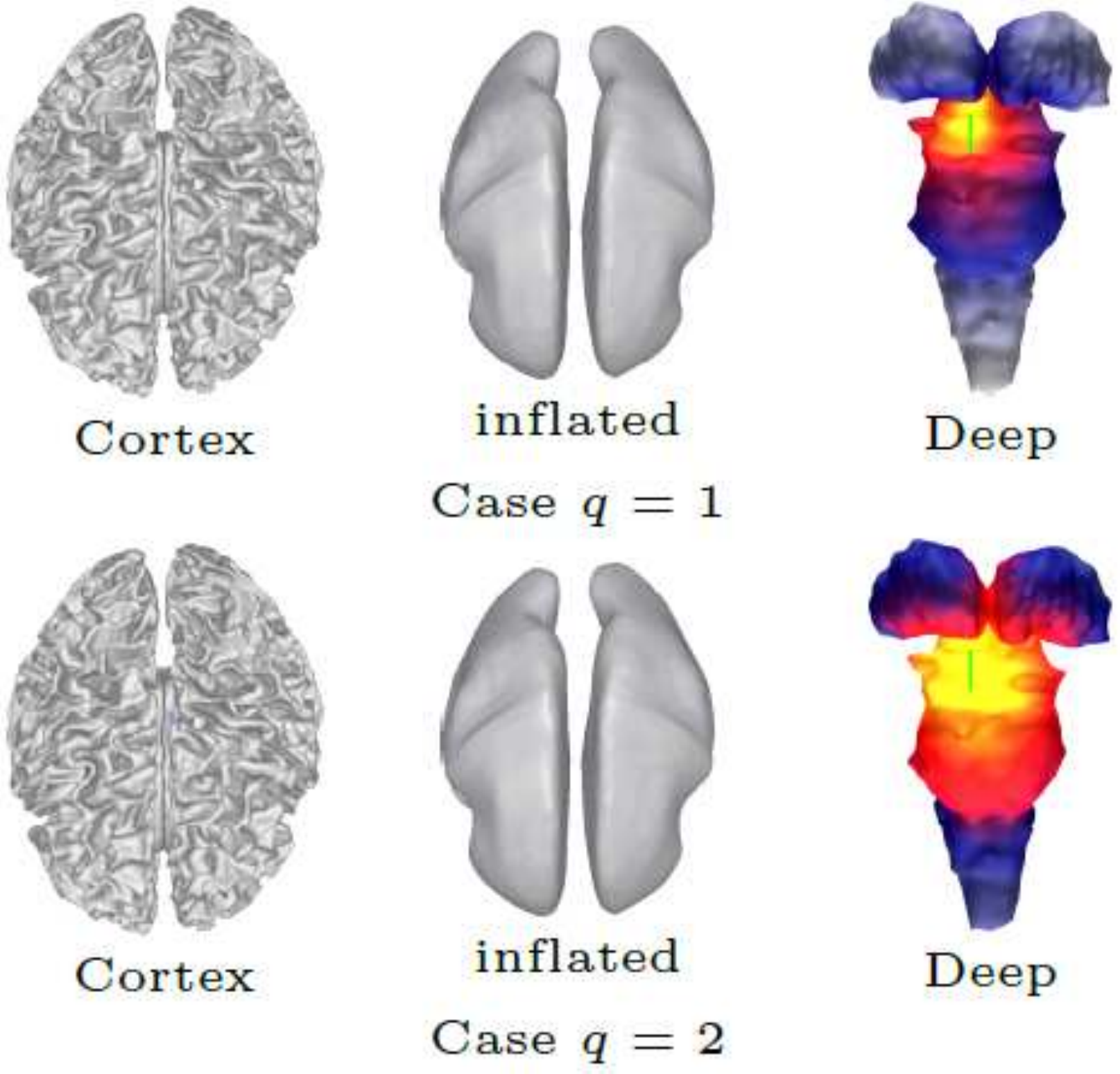} \\ \vskip0.07cm 
      {\bf Configuration (II): Iterative alternating sequential } \\ \vskip0.07cm
   \centering
    \includegraphics[height=3.55cm]{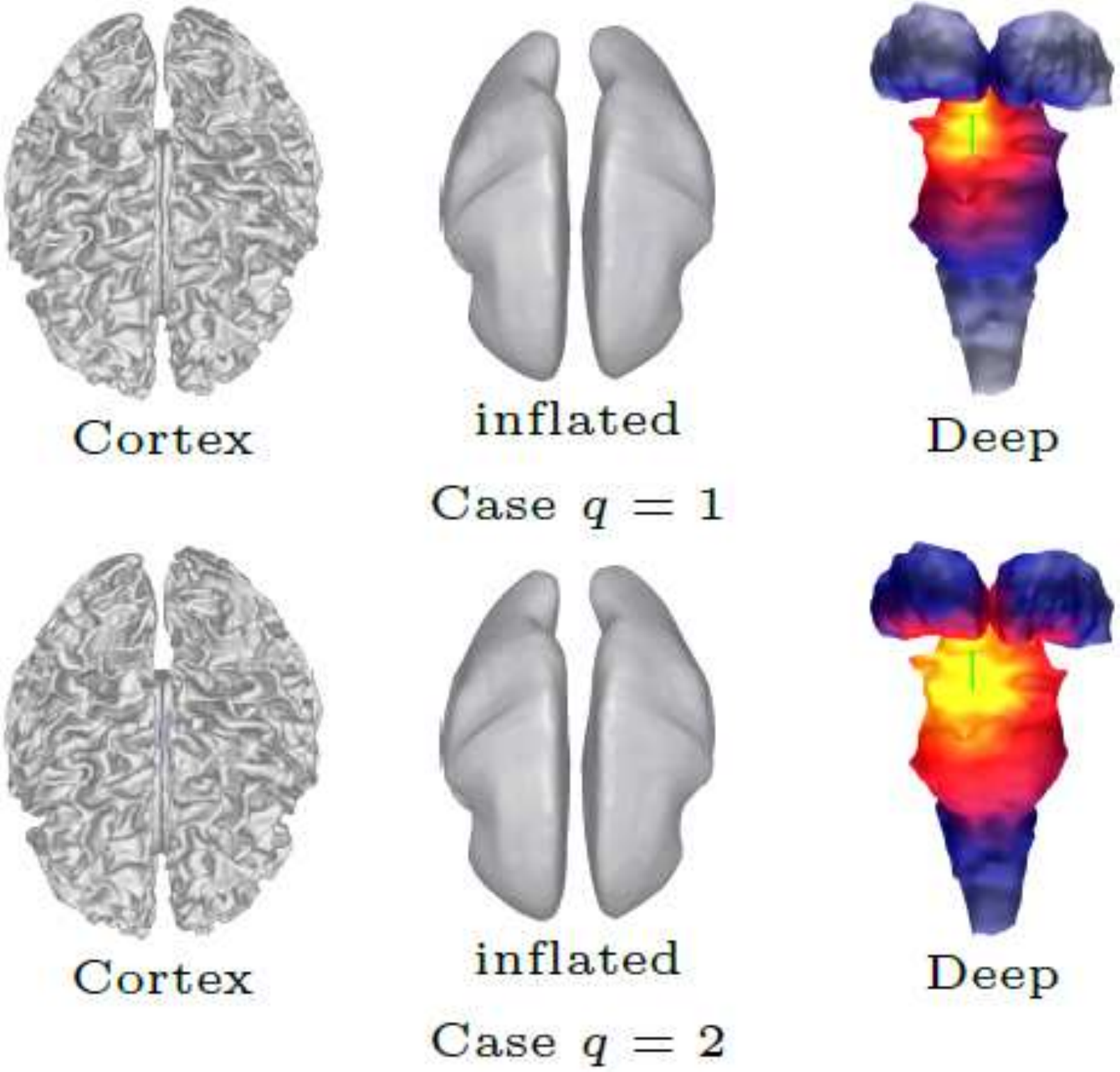} \\ \vskip0.07cm
     {\bf Configuration (II): SESAME } \\ \vskip0.07cm
   \centering
    \includegraphics[height=1.5cm]{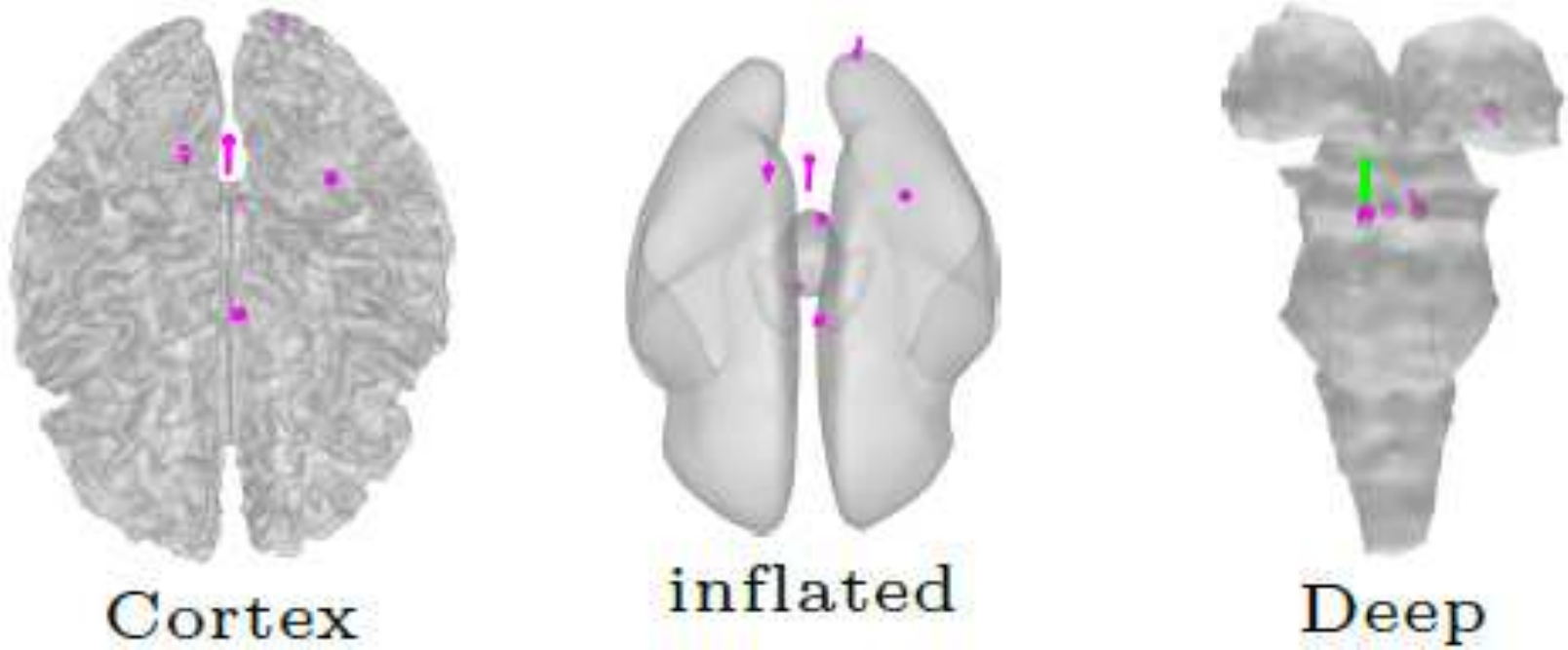} \\ \vskip0.07cm
    \end{tiny}
    \caption{The reconstructions obtained with the MRI-based head multi-compartment model with 3 \% noise for source   configuration ({\bf II}), including a single deep source placed in the brainstem (green pointer). The maximum found corresponds to the actual position,  a similar amplitude is obtained with both $q = 1$ and  $q = 2$, the focality being greater with $q = 1$. The EM and IAS reconstruction methods were found to perform essentially similarly for both source configurations with the most significant differences in the deep component of the configuration ({\bf I}). Magenta pointers in SESAME section represents its dipole estimations. }\label{fig:mri_based_2}
     \end{figure}
           \begin{figure}[h!]
    \centering
    \begin{tiny}

           {\bf  Configuration (III): Expectation maximization} \\ \vskip0.07cm
     \centering
     \includegraphics[height=3.55cm]{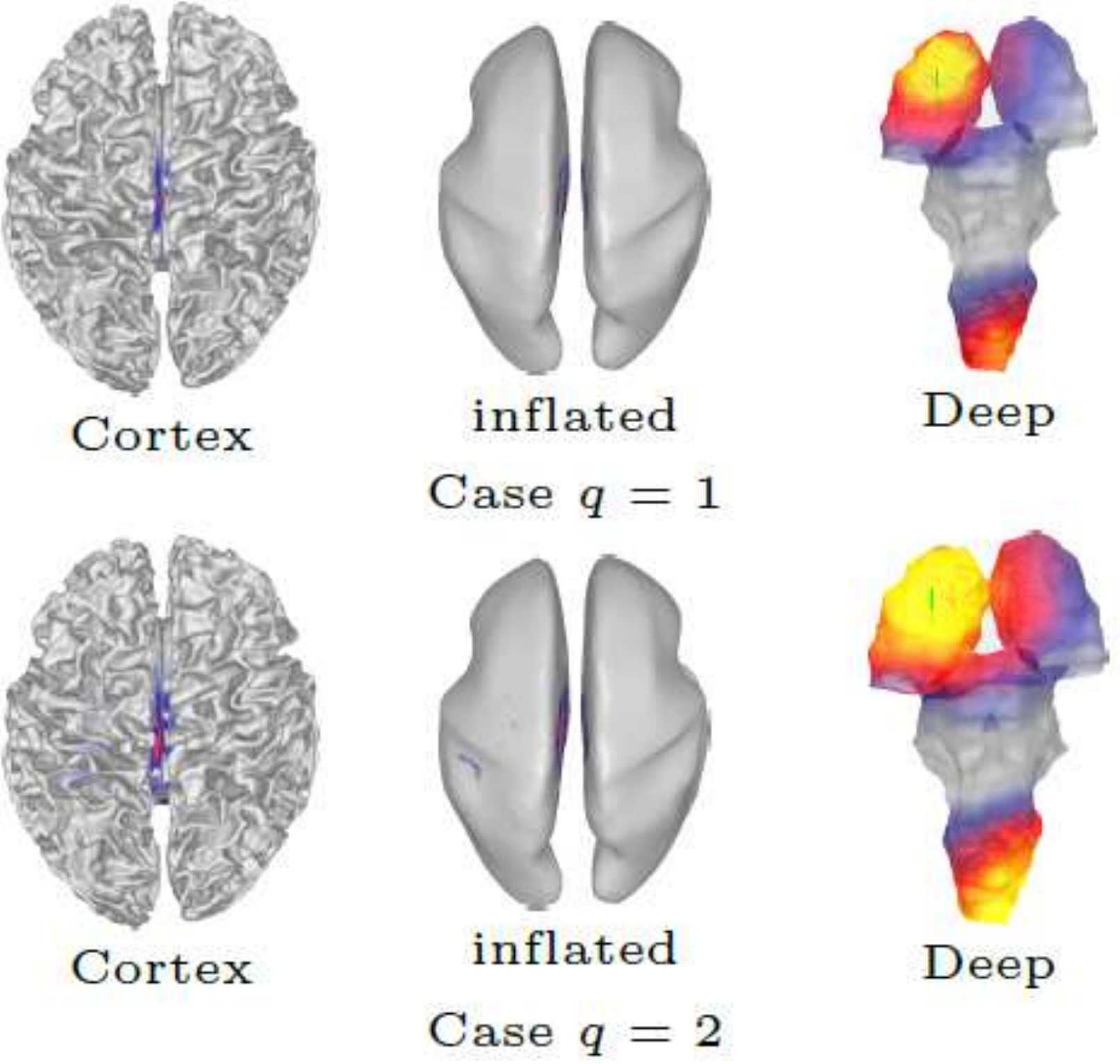}\\ \vskip0.07cm 
      {\bf Configuration (III): Iterative alternating sequential } \\ \vskip0.07cm
   \centering
    \includegraphics[height=3.55cm]{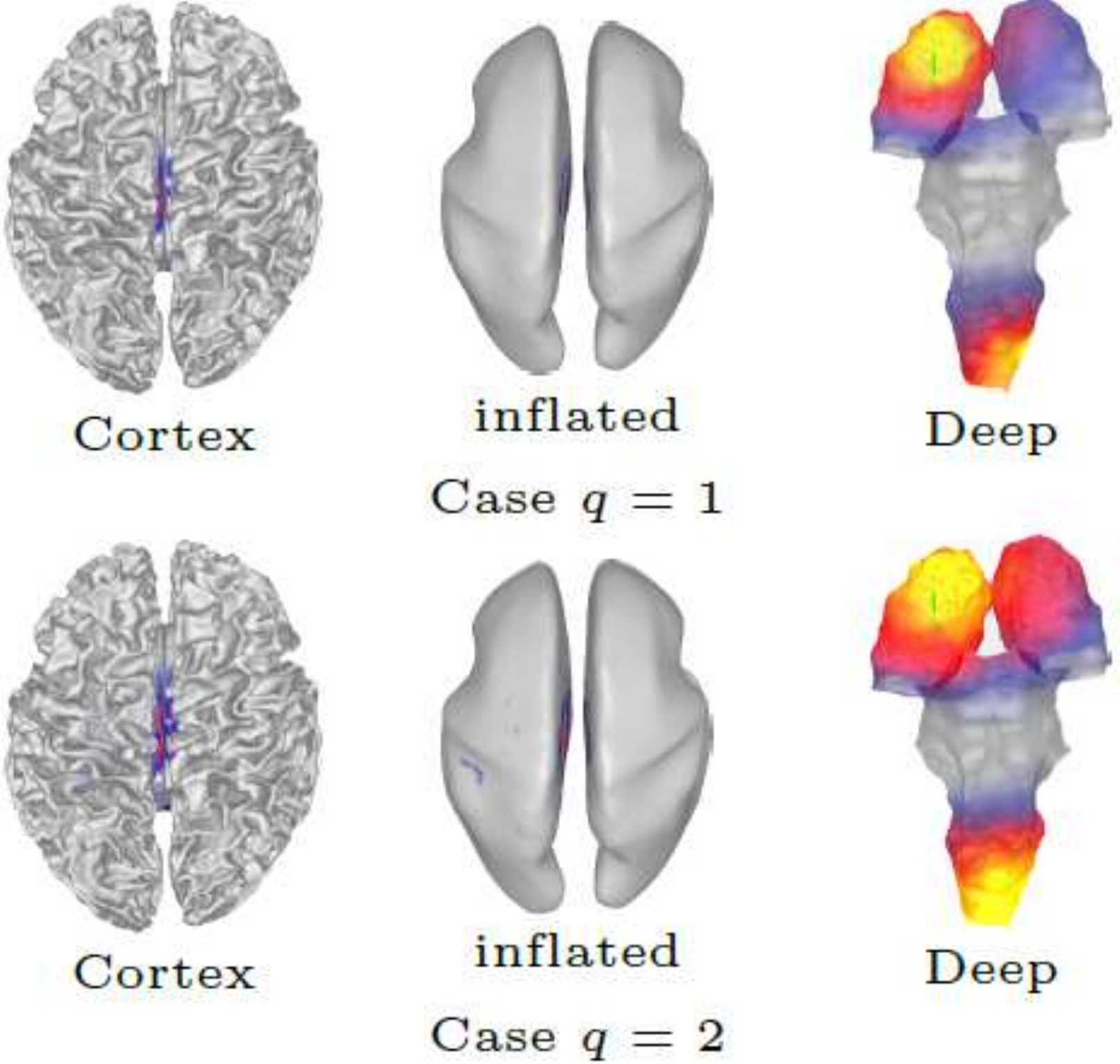} \\ \vskip0.07cm
     {\bf Configuration (III): SESAME } \\ \vskip0.07cm
   \centering
     \includegraphics[height=1.5cm]{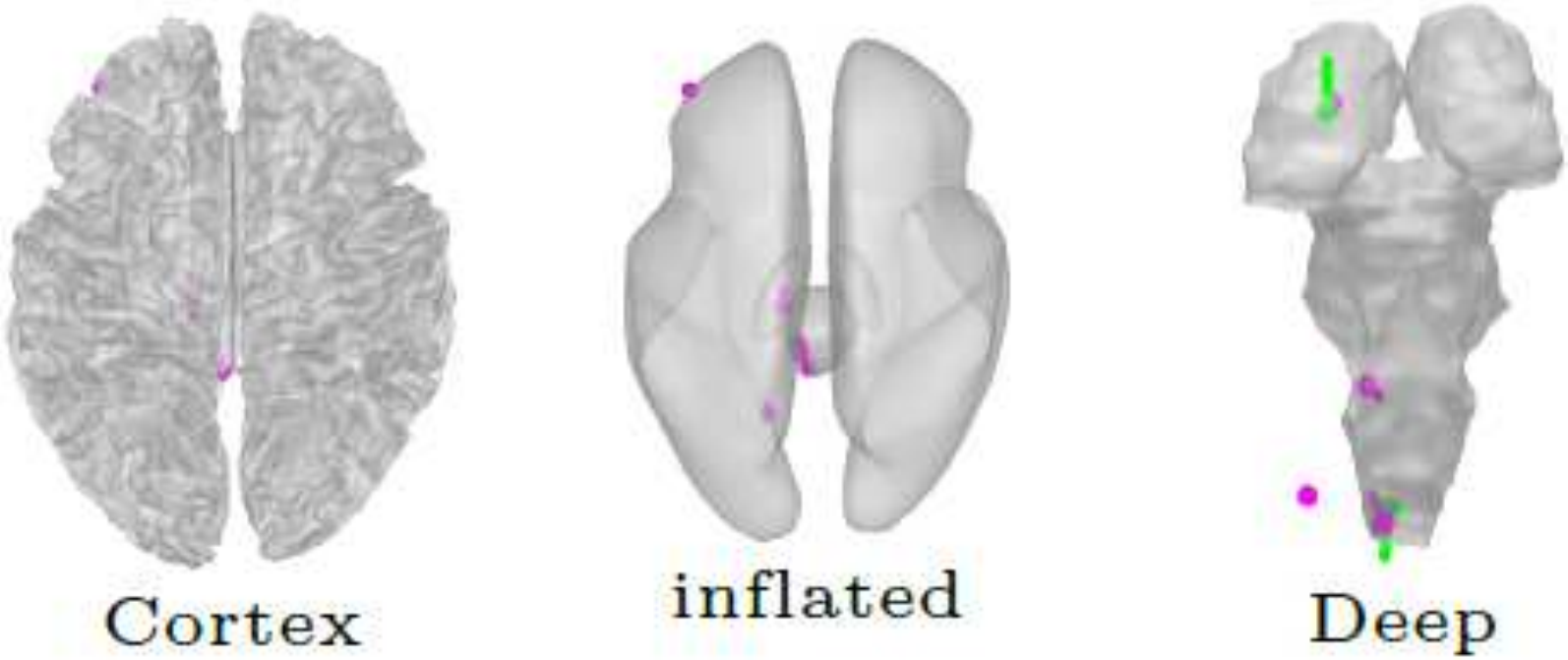} \\ \vskip0.07cm
    \end{tiny}
    \caption{The reconstructions obtained with the MRI-based head multi-compartment model with 3 \% noise for source  configuration ({\bf III}) with a quadrupolar deep source configuration formed by two dipolar components, one  in the ventral part of the left thalamus \cite{buchner1995origin} and another one  in lower medulla of the brainstem \cite{hsieh1995interaction}  (green arrows). The quadrupolar configuration is detected with both $q = 1$ and $q = 2$, and also with SESAME (magenta pointers). The reconstruction obtained is more focal in the former case. }
    \label{fig:mri_based_3}
     \end{figure} 
     

               \begin{figure*}[h!]
    \centering
    \begin{tiny}
     \vskip0.07cm \hrulefill \vskip0.07cm  
   {\bf Configuration (I)} \\ \vskip0.07cm   
\centering
\includegraphics[height=6.60cm]{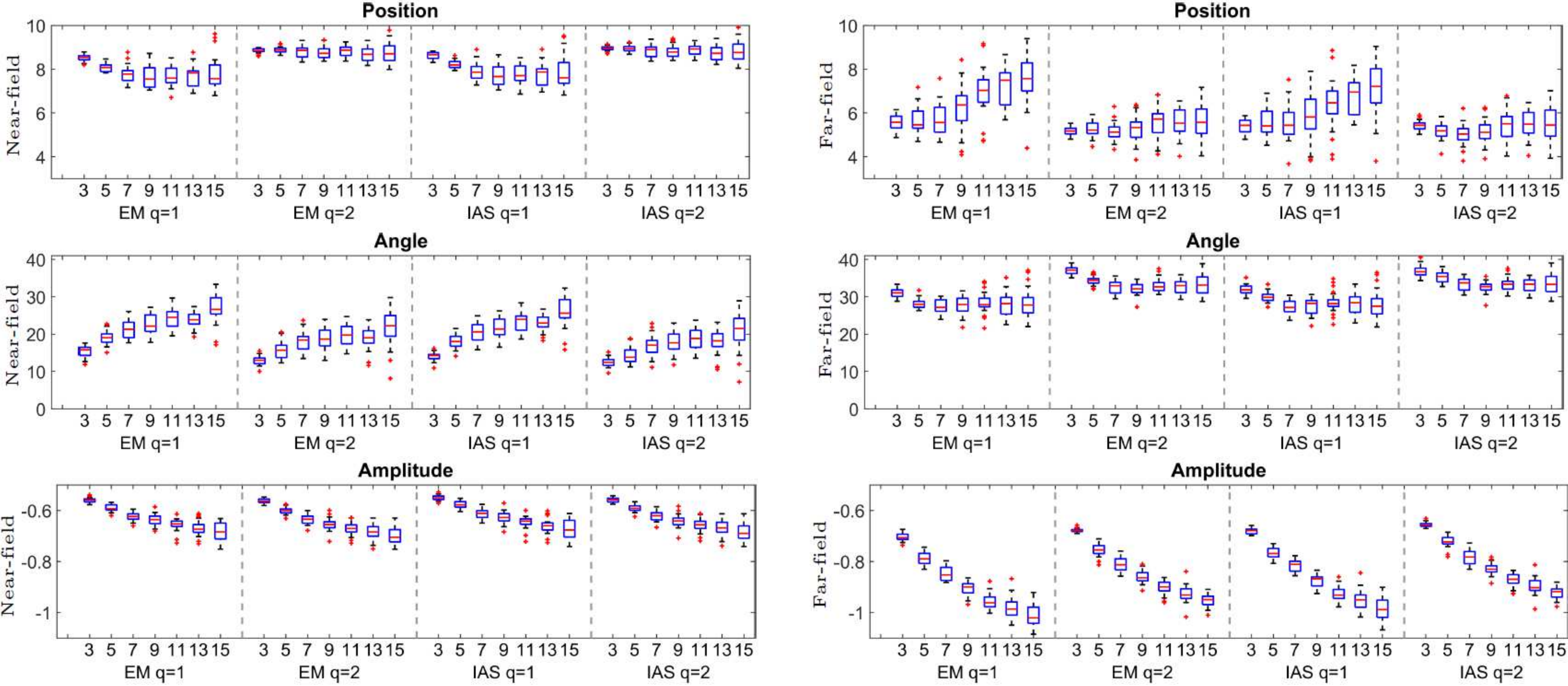} 
\\  \vskip0.07cm \hrulefill \vskip0.07cm  
{\bf Configuration (II)} \\ \vskip0.07cm
\centering
\includegraphics[height=6.66cm]{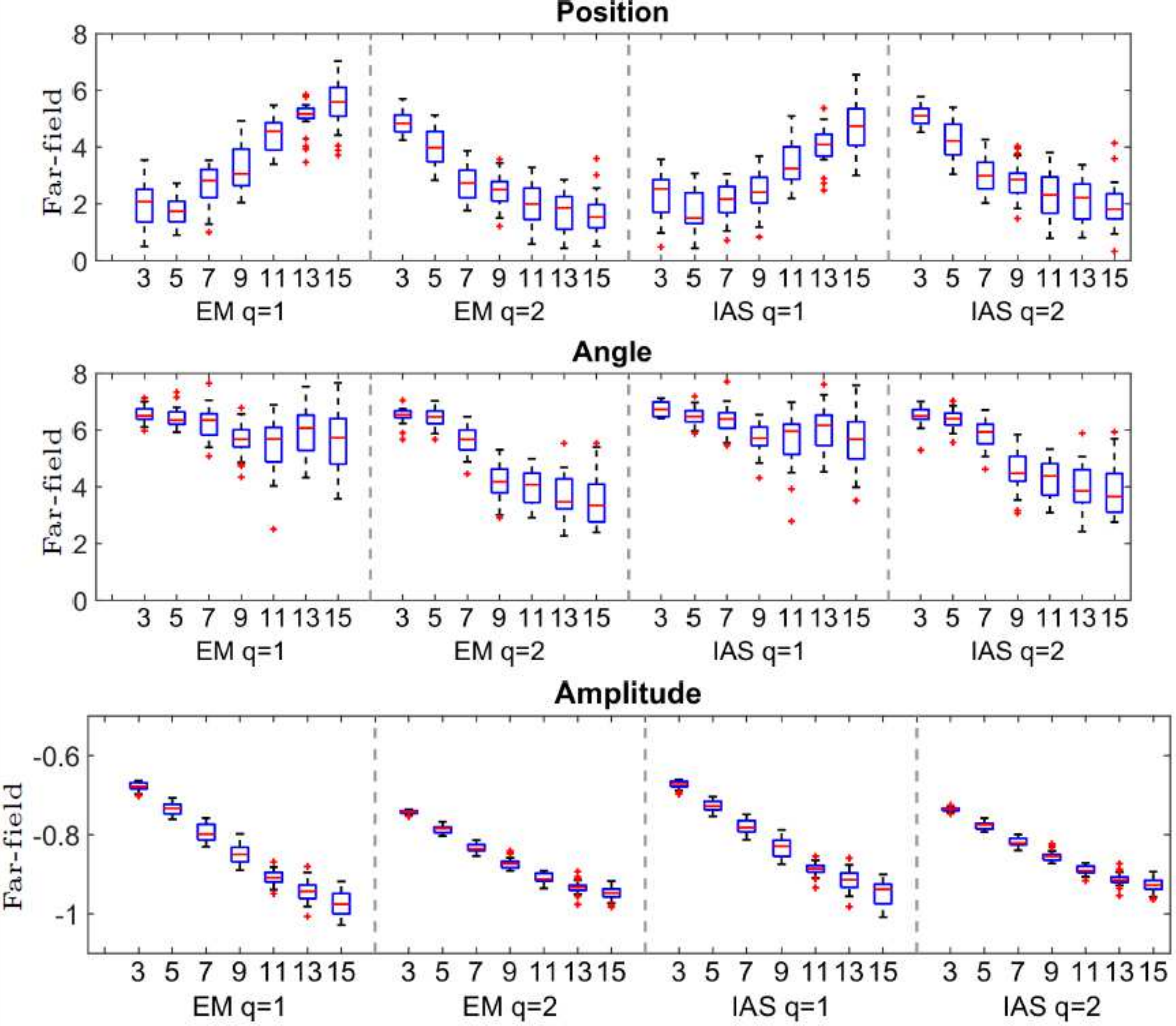}\\  \vskip0.07cm \hrulefill \vskip0.07cm  
        \end{tiny}
    \caption{The source localization accuracy measures  (position, angle and amplitude) for multiple noise levels. A total of 25 estimates have been calculated for each noise level. The error measures are evaluated in the spherical Ary model for source configurations ({\bf I}) and ({\bf II}) as shown in Section \ref{Sect:Sythetic_data} by applying CEP with the EM and IAS algorithm and both prior degrees 1 and 2.}
    \label{fig:boxplot_EXP}
\end{figure*}

\begin{figure*}[h!]
    \centering
    \begin{tiny}
     \vskip0.07cm \hrulefill \vskip0.07cm  
   {\bf Configuration (I): SESAME} \\ \vskip0.07cm 
\centering
 \includegraphics[height=4.60cm]{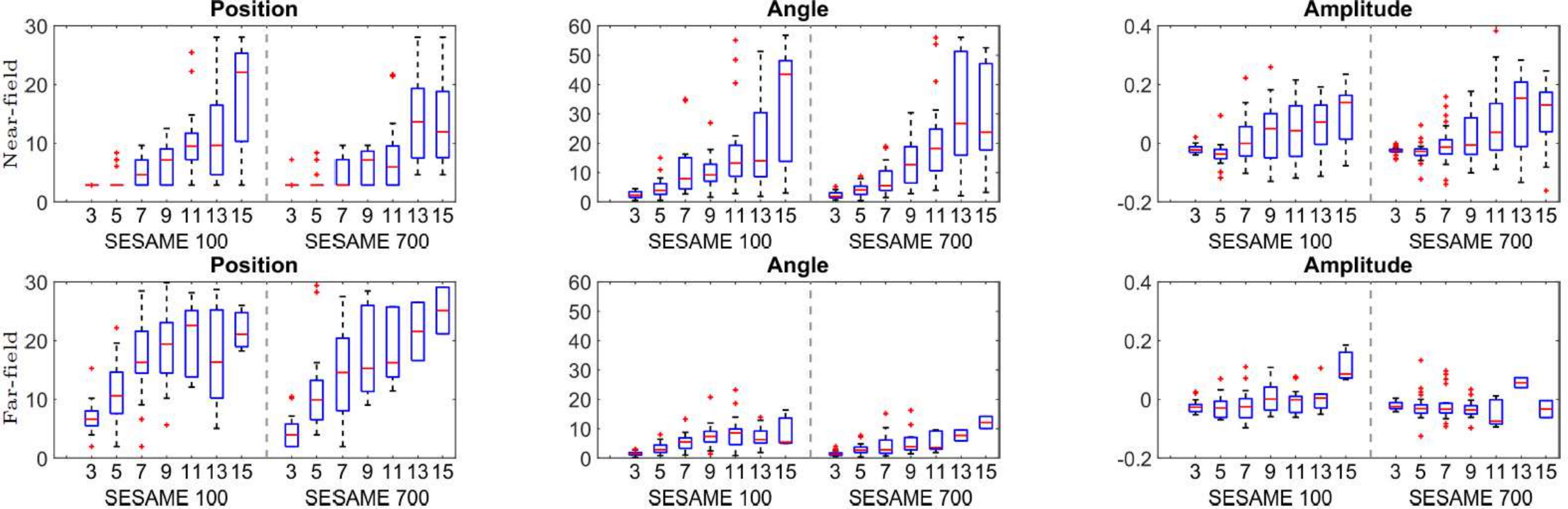} 
\\  \vskip0.07cm \hrulefill \vskip0.07cm  
  {\bf Configuration (II): SESAME} \\ \vskip0.07cm
\centering
 \includegraphics[height=2.30cm]{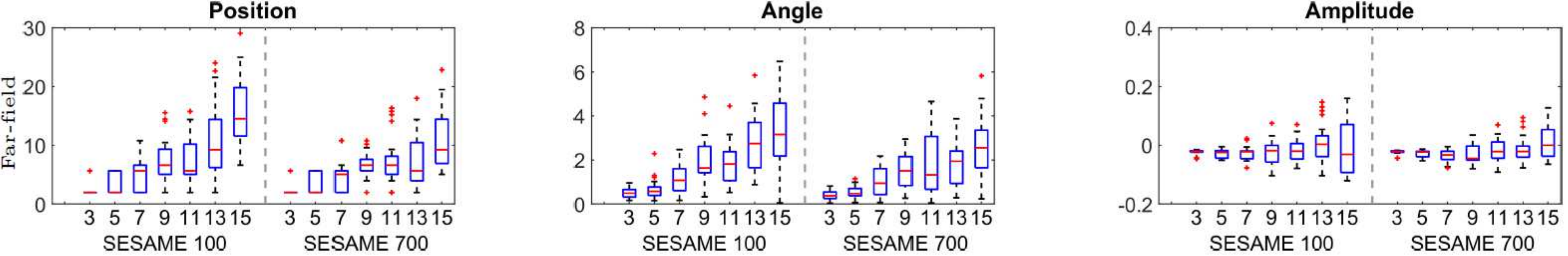} \\
    \vskip0.07cm \hrulefill  \vskip0.07cm
        \end{tiny}
    \caption{Dipole  reconstruction errors for SESAME-based source localization estimates calculated for multiple noise levels. Each box-plot bar has been obtained based on 25 reconstructions. }
    \label{fig:boxplot_SESAME}
\end{figure*}

               \begin{figure}[h!]
    \centering
    \begin{tiny}
     \vskip0.07cm \hrulefill \vskip0.07cm  
   {\bf Configuration (I)} \\ \vskip0.07cm   
\begin{minipage}{7.0cm}
 \includegraphics[height=2.10cm]{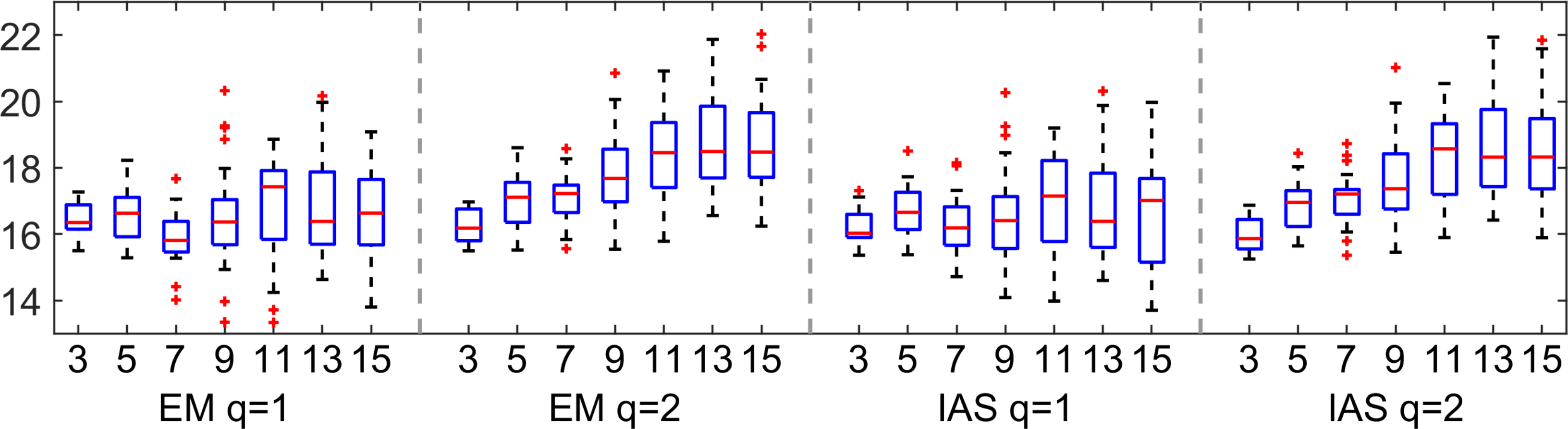} 
    \end{minipage}
      \\
      \vskip0.07cm \hrulefill  \vskip0.07cm 
  {\bf Configuration (II)} \\ \vskip0.07cm
  \begin{minipage}{7.0cm}
 \includegraphics[height=2.10cm]{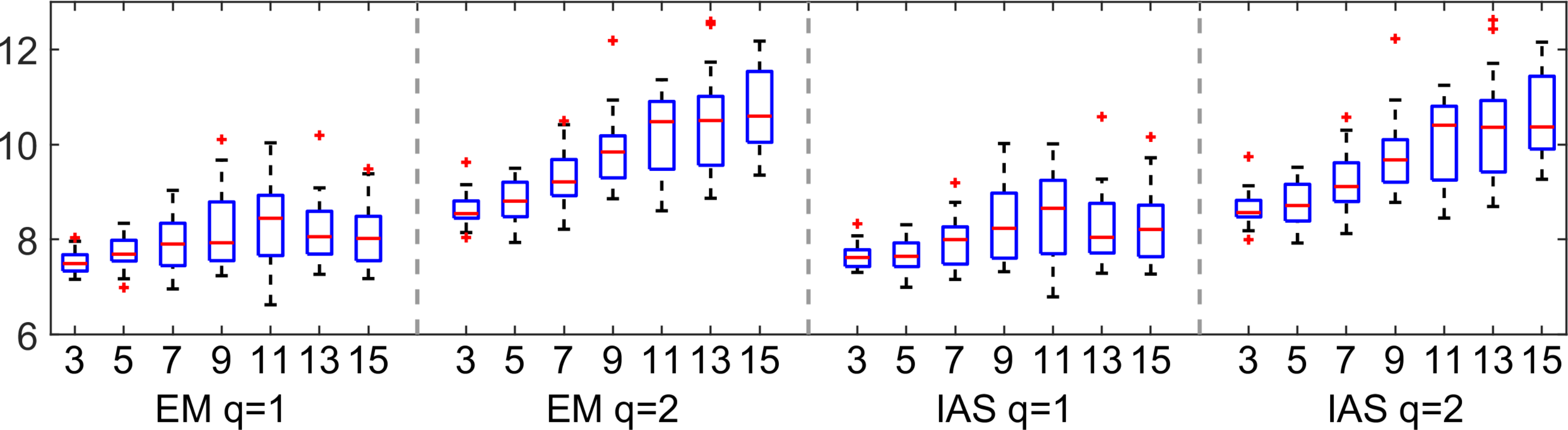} 
    \end{minipage}
      \\
      \vskip0.07cm \hrulefill  \vskip0.07cm 
  {\bf Configuration (I): SESAME} \\ \vskip0.07cm
\begin{minipage}{7.0cm}
\centering
    \includegraphics[height=2.10cm]{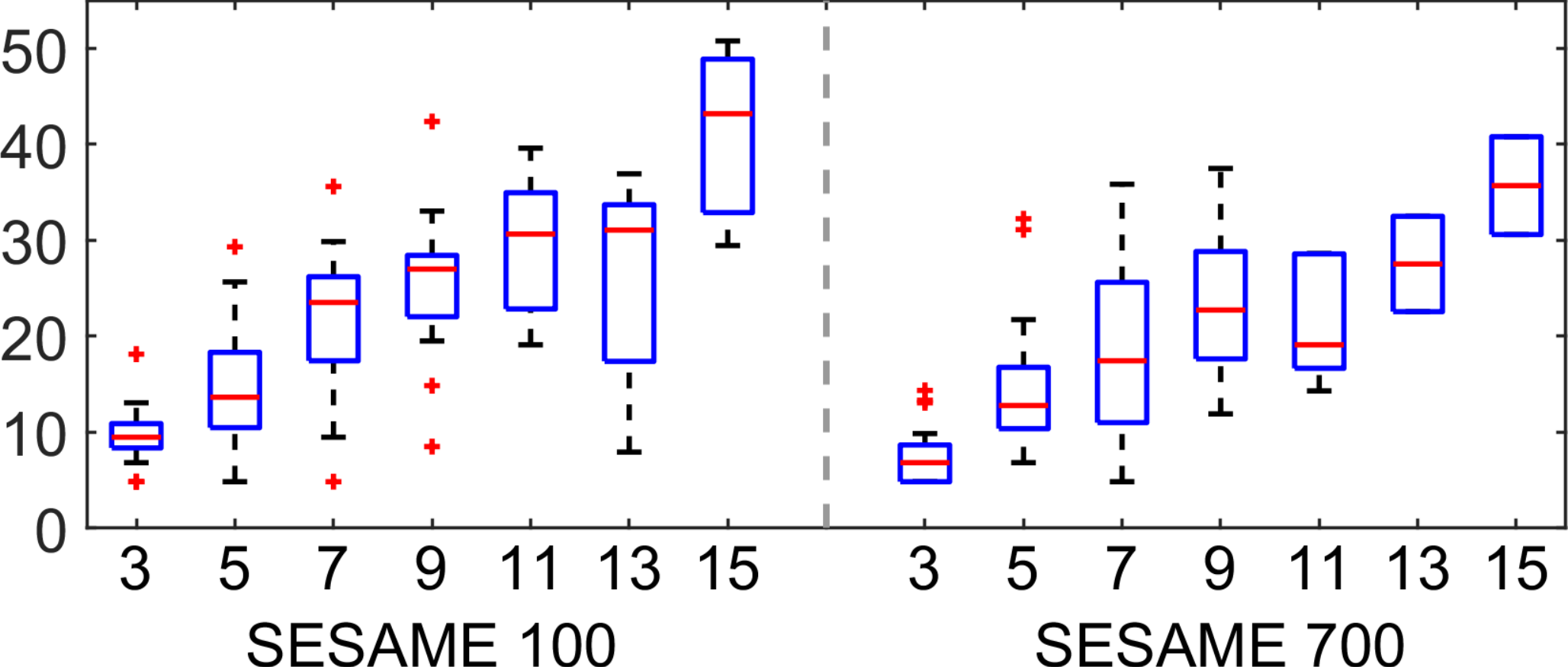} 
    \end{minipage}
      \\
       \vskip0.07cm \hrulefill  \vskip0.07cm 
{\bf Configuration (II): SESAME} \\ \vskip0.07cm
\begin{minipage}{7.0cm}
\centering
    \includegraphics[height=2.10cm]{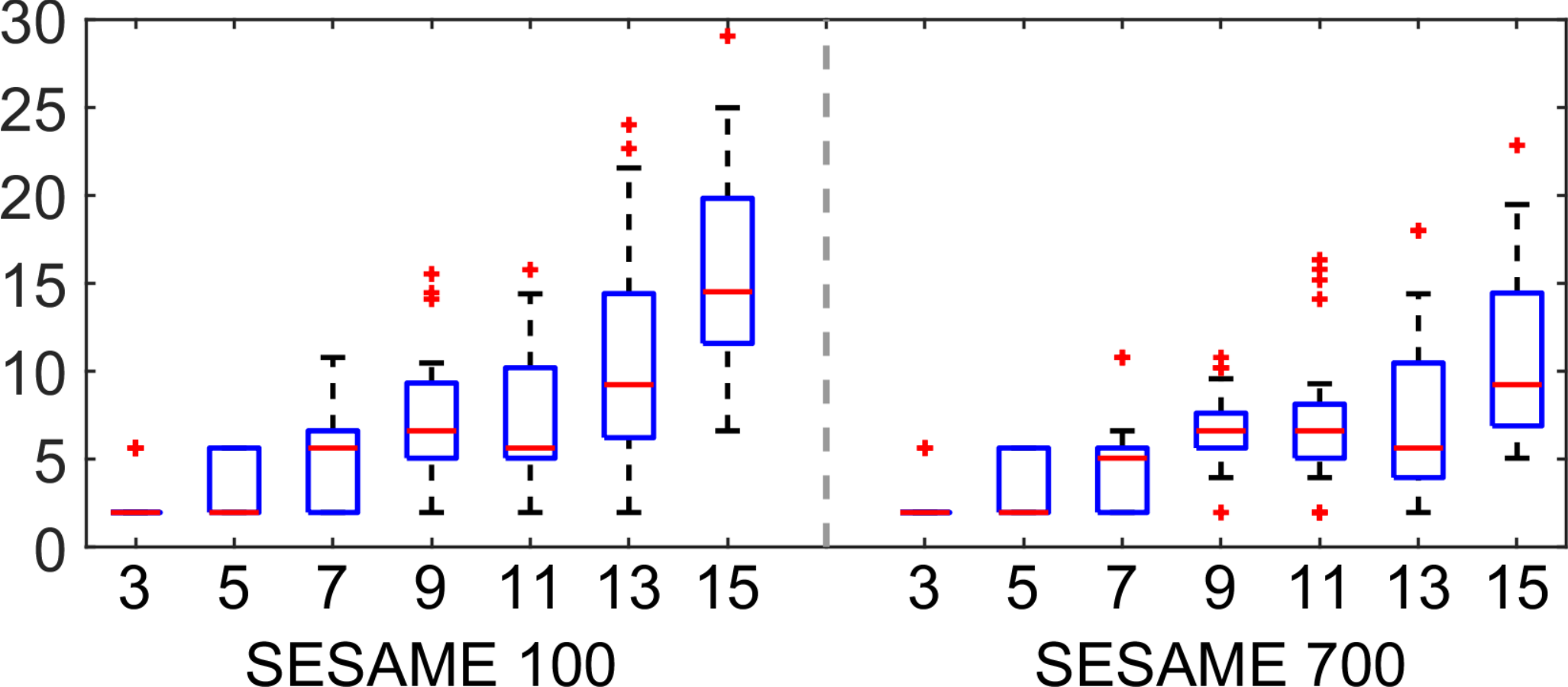} 
    \end{minipage}
      \\
      \vskip0.07cm \hrulefill  \vskip0.07cm 
        \end{tiny}
    \caption{The earth mover's distance (EMD) evaluated in a 45 mm radius sphere co-centered with the true source for noise levels from 3 \% to 15 \%. For  SESAME, the  dipoles with a distance greater than 45 mm to the true source location were excluded from the EMD evaluation. Each box-plot bar has been obtained using 25 different reconstructions.}
    \label{fig:boxplot_EMD}
\end{figure}
    \begin{figure}[h!]
    \centering
    \begin{tiny}

     {\bf  Configuration (I): Expectation maximization} \\ \vskip0.07cm
     \centering
    \includegraphics[height=3.55cm]{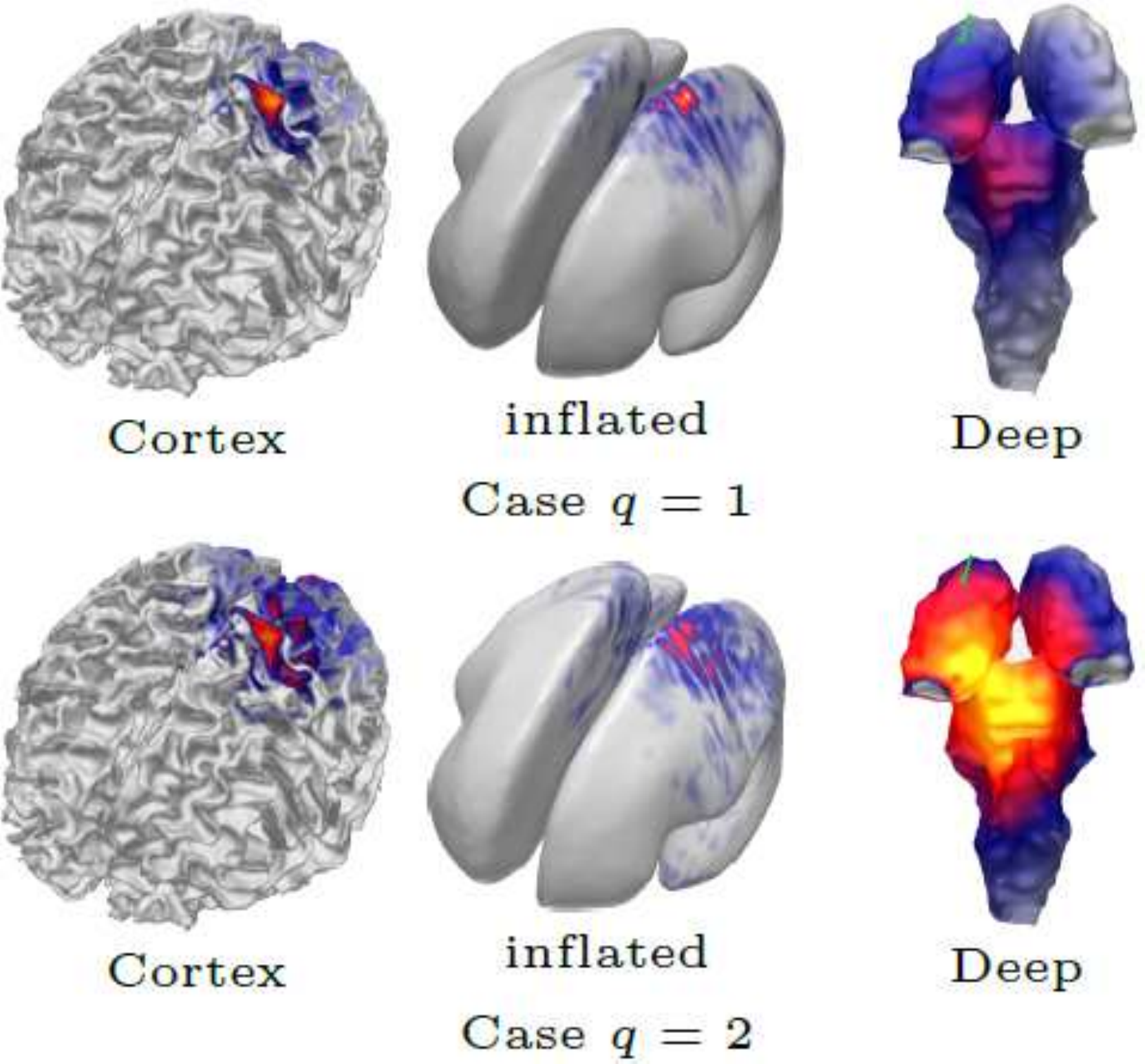} \\
 \vskip0.07cm      
      {\bf Configuration (I): Iterative alternating sequential } \\ \vskip0.07cm
   \centering
    \includegraphics[height=3.55cm]{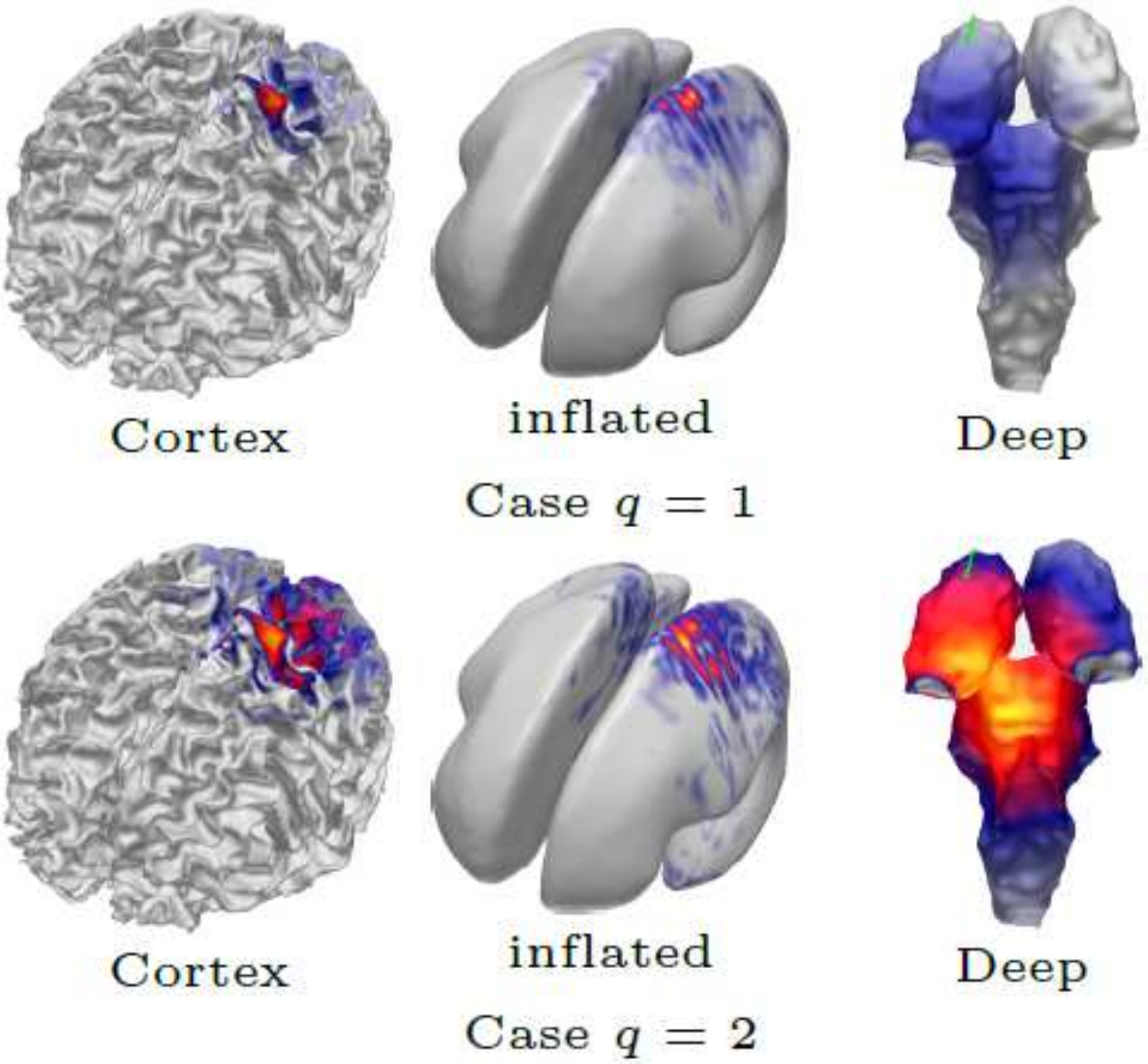} \\
        \vskip0.07cm   
     {\bf Configuration (I): SESAME } \\ \vskip0.07cm
   \centering
     \includegraphics[height=1.5cm]{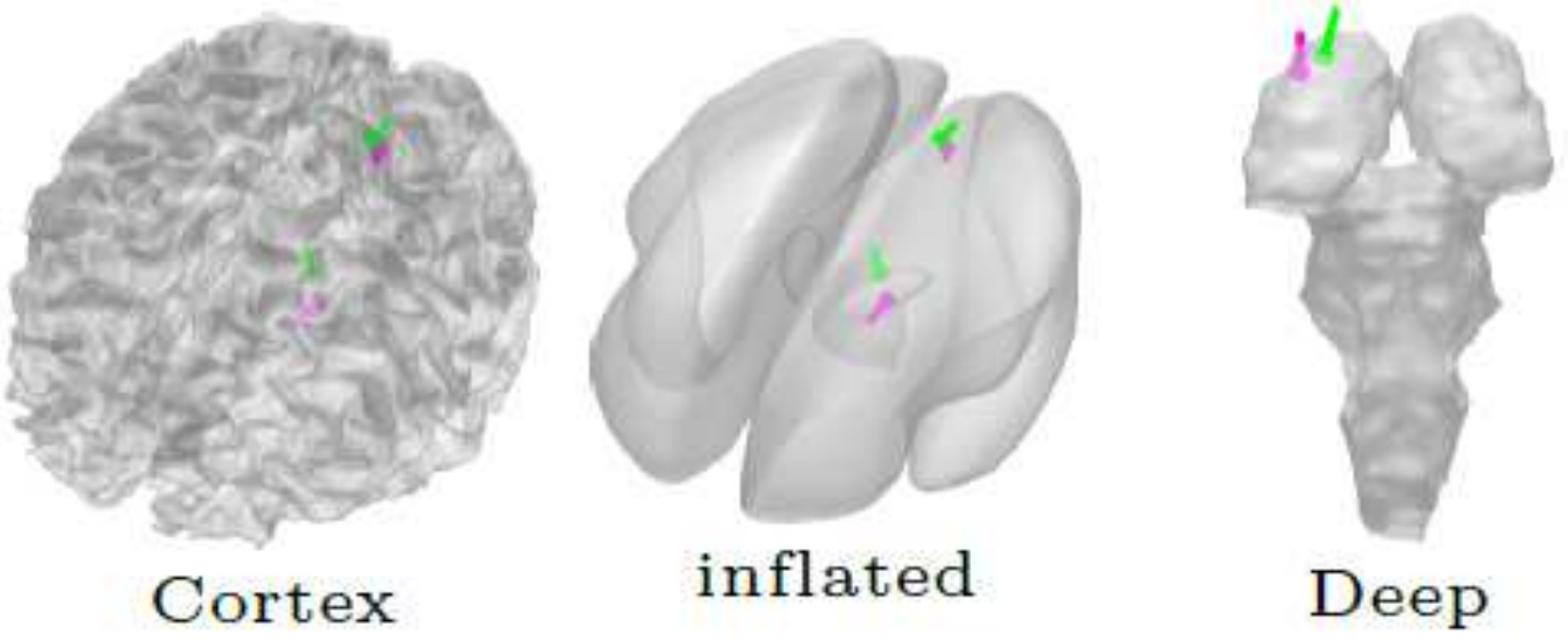} \\  \vskip0.07cm
    \end{tiny}
    \caption{The reconstructions obtained with 5 \% noise and the MRI-based head multi-compartment model for source  configuration ({\bf I})  including two sources, one placed in the left 3b Brodmann area of the central sulcus, pointing inwards in the direction of the local surface normal vector,  and a vertical source placed in the ventral posterolateral part of the left thalamus (Figure \ref{fig:source_configurations}).  The actual source position is visualized by a green pointer in each image. Magenta pointers are dipole realizations of SESAME.}\label{fig:mri_based_1_5}
     \end{figure}
  
      \begin{figure}[h!]
    \centering
    \begin{tiny}
    {\bf  Configuration (II): Expectation maximization} \\ \vskip0.07cm
     \centering
    \includegraphics[height=3.55cm]{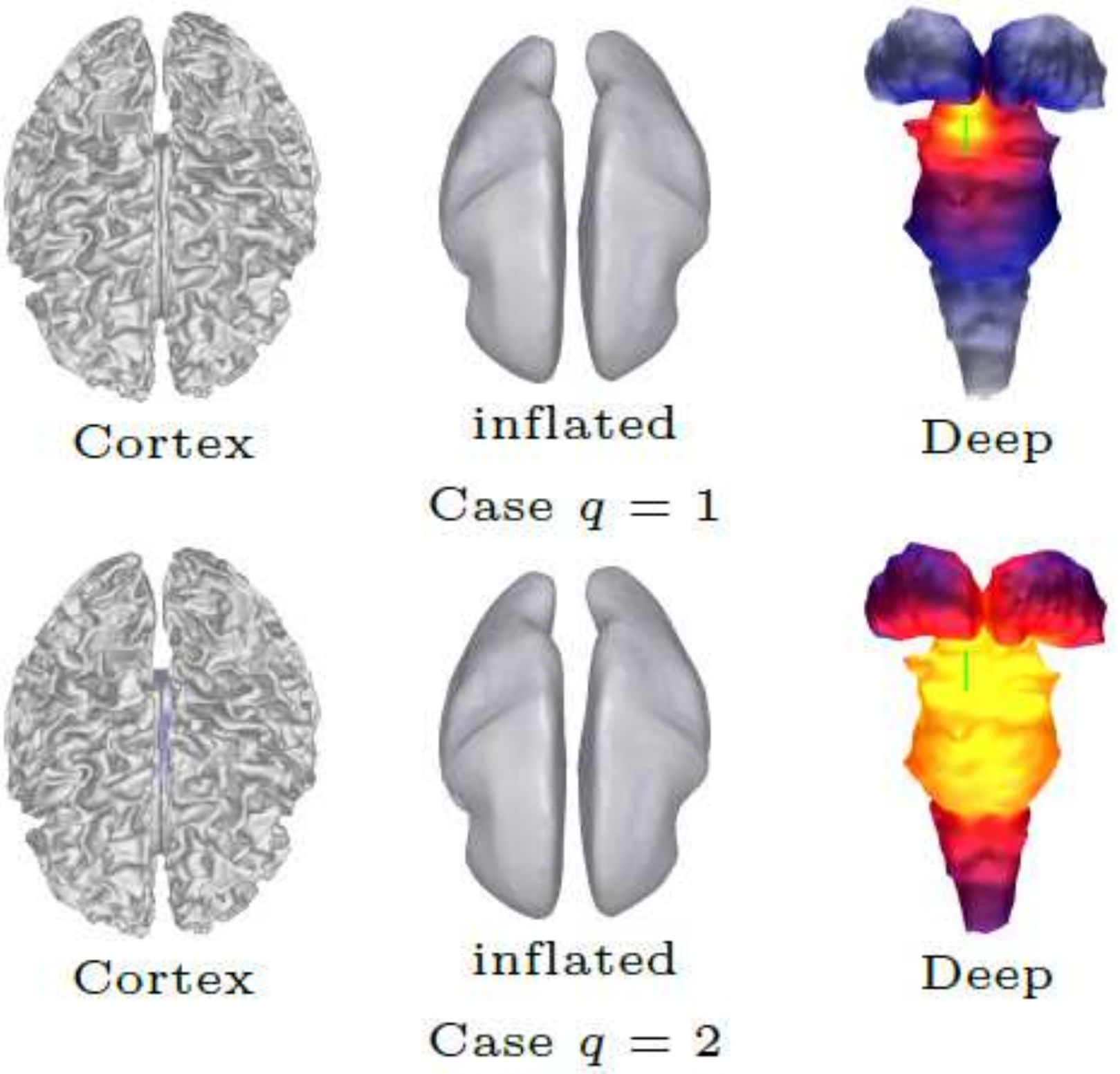} \\ \vskip0.07cm 
      {\bf Configuration (II): Iterative alternating sequential } \\ \vskip0.07cm
   \centering
    \includegraphics[height=3.55cm]{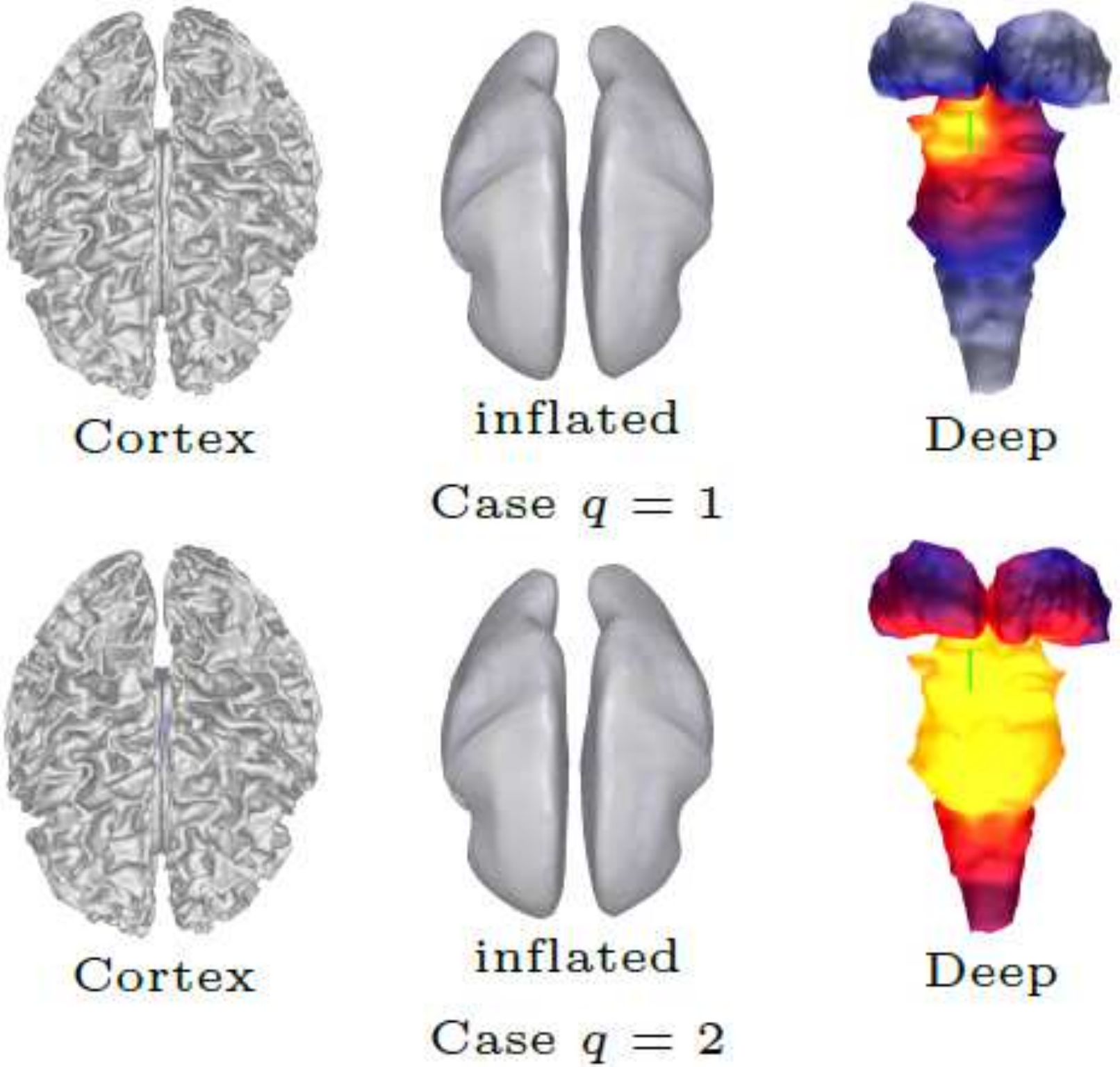} \\ \vskip0.07cm
     {\bf Configuration (II): SESAME } \\ \vskip0.07cm
   \centering
    \includegraphics[height=1.5cm]{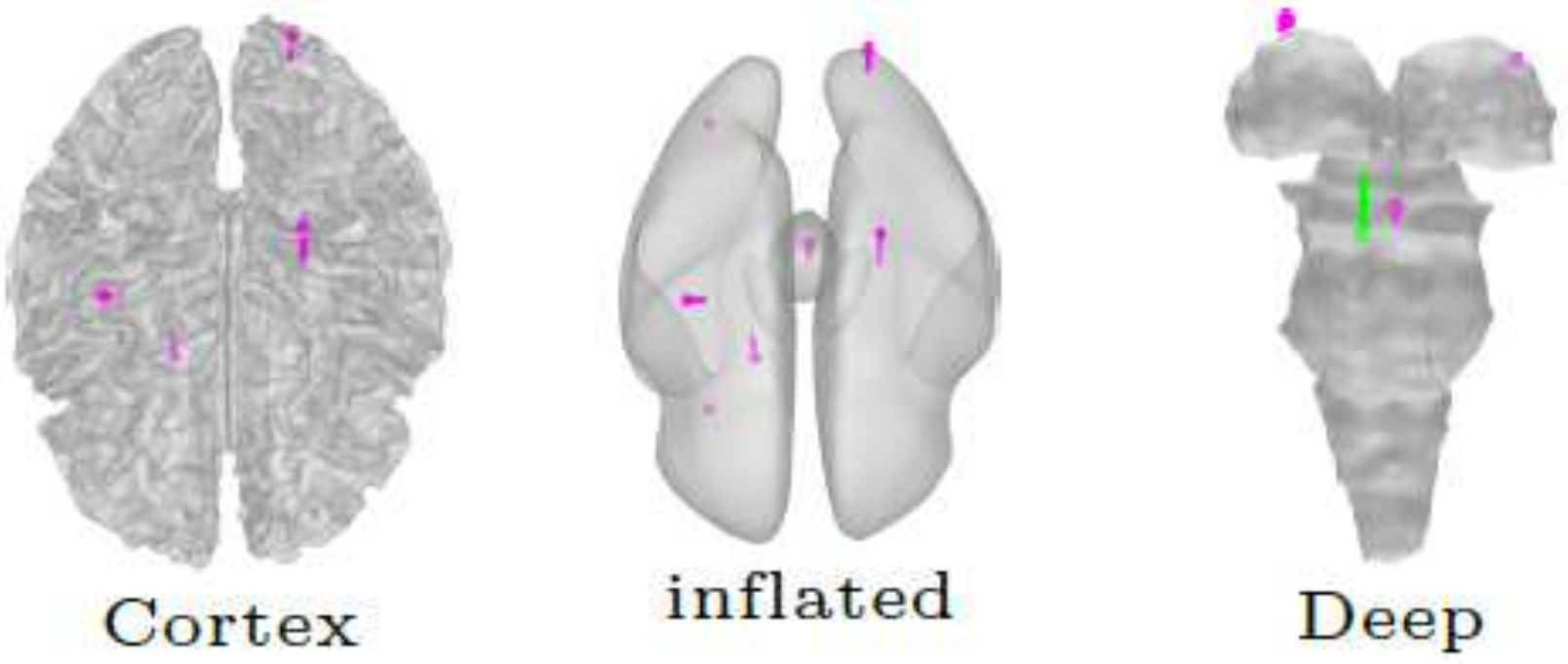}\\ \vskip0.07cm
    \end{tiny}
    \caption{The reconstructions obtained with 5 \% noise and the MRI-based head multi-compartment model for source   configuration ({\bf II}), including a single deep source placed in the brainstem (green pointer). The dipole estimations of SESAME are presented by magenta pointers.}\label{fig:mri_based_2_5}
     \end{figure}
     \begin{figure}[h!]
    \centering
    \begin{tiny}
           {\bf  Configuration (III): Expectation maximization} \\ \vskip0.07cm
    \includegraphics[height=3.55cm]{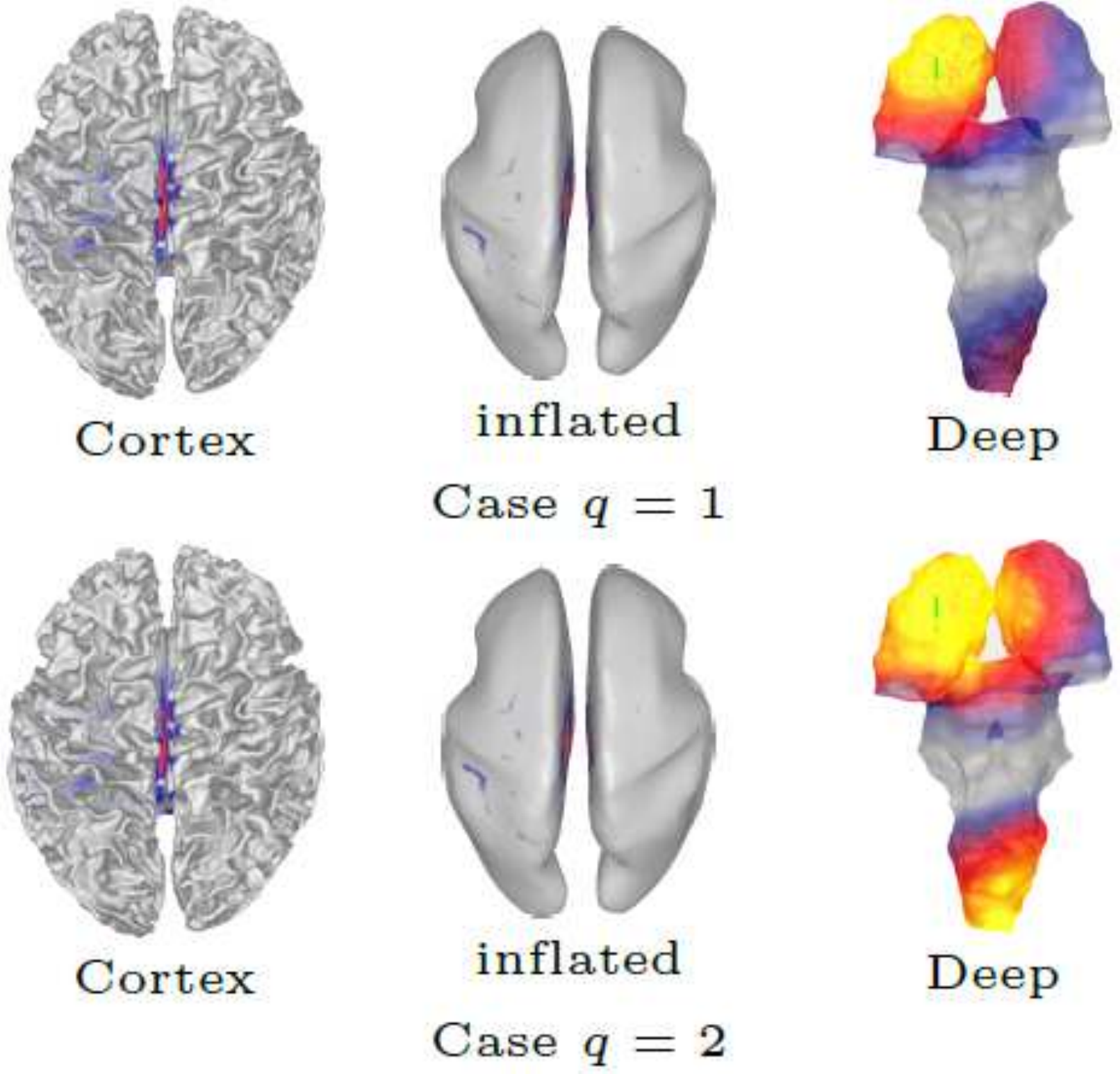} \\ \vskip0.07cm 
      {\bf Configuration (III): Iterative alternating sequential } \\ \vskip0.07cm
   \centering
    \includegraphics[height=3.55cm]{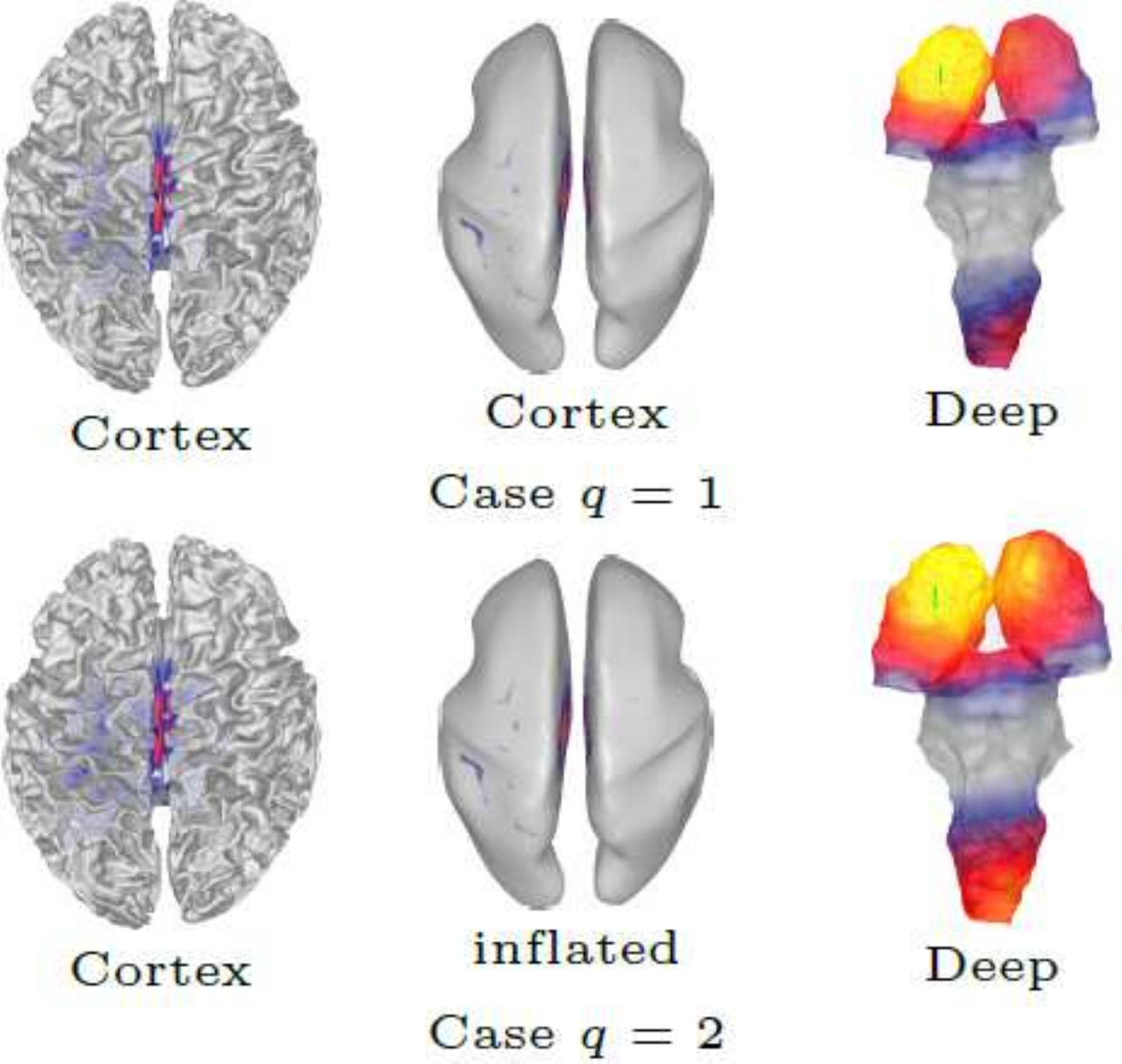} \\ \vskip0.07cm
     {\bf Configuration (III): SESAME } \\ \vskip0.07cm
    \includegraphics[height=1.5cm]{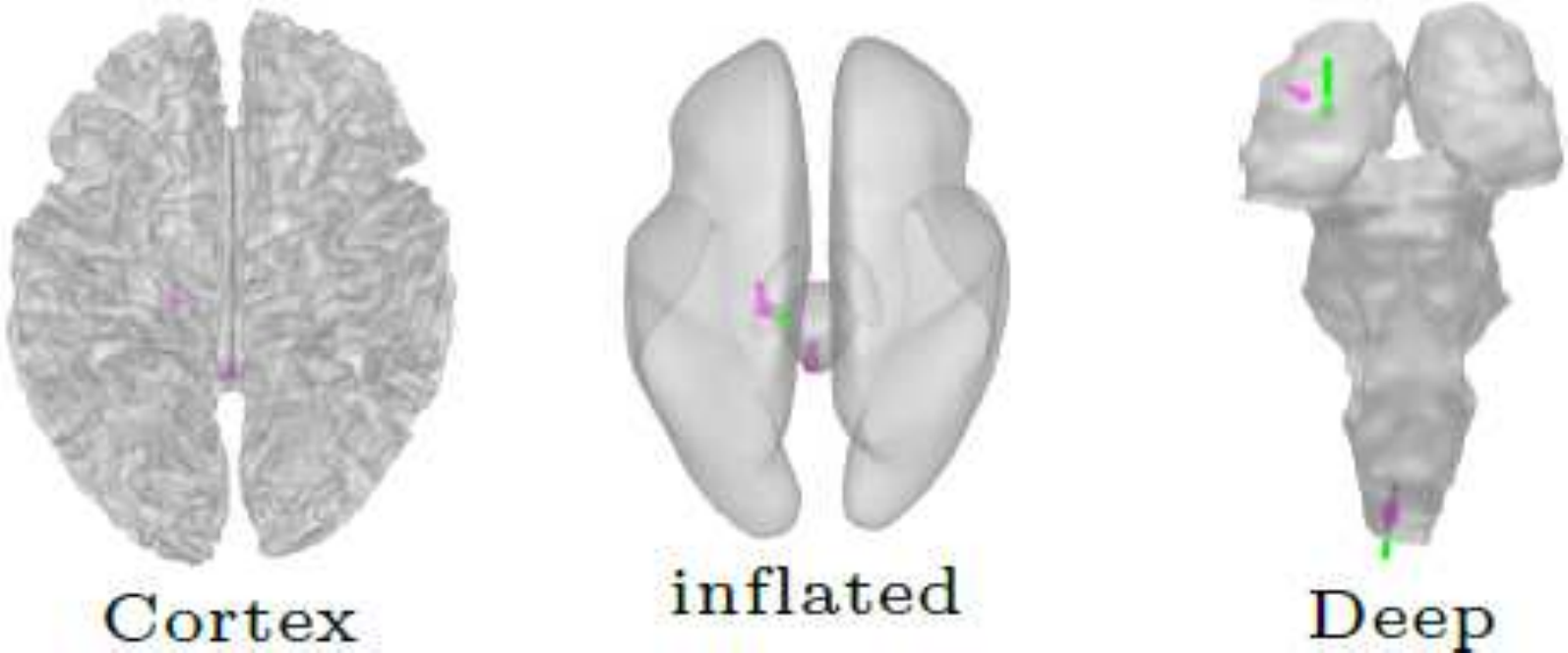}\\ \vskip0.07cm
    \end{tiny}
    \caption{The reconstructions obtained with 5 \% noise and the MRI-based head multi-compartment model for source  configuration ({\bf III}) with a quadrupolar deep source formed by two dipolar components, one  in the ventral part of the left thalamus \cite{buchner1995origin} and another one  in lower medulla of the brainstem \cite{hsieh1995interaction}  (green arrows). SESAME has found both of the dipole locations with high accuracy (Magenta pointers).}\label{fig:mri_based_3_5}
     \end{figure}

\section{Results}
\label{sec:results} 

The accuracy and focality results obtained with the spherical Ary model are shown as histograms in  Figures \ref{fig:results_accuracy5},   \ref{fig:results_focality5}, \ref{fig:results_EarthMover5} and \ref{fig:results_SESAME_accuracy} for 5 \% noise.  The CEP prior model reconstructs the near- and far-field activity with both prior degrees $q = 1$ and $q = 2$ and reconstruction techniques EM and IAS. While the degree of the prior does not directly appear to affect the position and orientation accuracy, the reconstructions were overall more focal in the case  $q=1$ compared to $q = 2$. With $q = 1$ the reconstruction of the near-field source was more focal compared to that of the deep one, while with $q=2$ such a tendency was less obvious or absent.  EMD and hard threshold, measuring the dynamical structure of the reconstructed distribution (Section \ref{sec:focality_measures}), show a similar tendency as the hard threshold. Of these, EMD reveals smaller relative differences, suggesting that the  overall variation of the reconstruction is rather similar for $q=1$ and $q = 2$ or the amount of artifacts is roughly the same, while the focality of the maximum peak varies more significantly. The mutual differences between the EM and IAS reconstruction techniques are less obvious than those following from the depth of the source,  the degree of the prior and the noise level. The results for SESAME show high accuracy for a near-field source, but the localization error for far-field source is diminished. However, SESAME gives an accurate estimation for dipole angle and amplitude regardless the estimated location. It can also be seen that in the case of the P14/N14 component, there is no significant difference in the accuracy obtained with 100 samples and 700 samples. EMD results for the configuration ({\bf II}) shows the effect of sample size in the location estimates obtained for SESAME: the median does not change as much as 90 \% credibility is shrunk.

The results obtained with the MRI-based head mo\-del are visualized in Figures \ref{fig:mri_based_1}--\ref{fig:mri_based_3} and in Figures \ref{fig:mri_based_1_5}--\ref{fig:mri_based_3_5} for 3 and 5 \% noise, respectively. The first-degree CEP, i.e., the case $q = 1$, leads to an overall more focal  reconstruction compared to $ q = 2$. Akin to the results obtained for the Ary model, the near-field component, i.e., the cortical pattern obtained with configuration ({\bf I}) is more focal with $q = 1$ than with $q = 2$; in the former case it is clearly restricted to the Brodmann 3b area, where cortical true source is located. In the latter one, it spreads more clearly towards the back areas of the brain: posterior cortex, e.g., Brodmann areas 5 and 7. Neither of the cases finds the accurate location of the simultaneous deep activity which appears more focal with $q =1$. On the contrary, the dipole estimating SESAME finds both sources within high accuracy. For the single-source configuration ({\bf II}), agreeing with the case of the Ary model, the maximum of the deep source found corresponds to the actual position and  a similar amplitude is obtained with both $q = 1$ and  $q = 2$, while in the former case,  the distribution  is more concentrated around the actual source position at the medial lemniscus pathway, distinguishing the upper brainstem as the area of activity.  The  EM and IAS reconstruction methods were found to perform essentially similarly for both source configurations with the most significant differences in the deep component of the configuration ({\bf I}). SESAME finds the dipole location roughly, but produces many false sources with significant dipole strength.
For ({\bf III}),  the thalamic and sub-thalamic components of the quadrupolar configuration are detected with both prior degrees $q = 1$ and $q = 2$, the results being more focal in the former case concerning, especially, the lateral localization accuracy; with  $q = 1$, the activity is more clearly limited to the left lobe of the thalamus, whereas with $q = 2$, both lobes show activity. SESAME detect activities also correctly, altough, for 3 \% noise we obtain false cortical activities similarly to ({\bf II}).

The result for multiple noise levels shown in Figures \ref{fig:boxplot_EXP} and \ref{fig:boxplot_SESAME} shows that, on the one hand, CEP localizes the near-field component slightly better with $q=1$ compared to $q=2$. On the other hand, CEP localizes the far-field component better with $q=2$. Between EM and IAS posterior maximizers, there are no significant differences considering the overall accuracy of the spherical Ary model and EEG. The dipole amplitude reconstruction has a natural  downward trend with respect to increasing noise due to  increasing dispersion.  This decreasing tendency  is weaker with prior degree 1. The reason for this can be explained by the lasso step, which acts as an "inhibitor" of the propagation of the activity estimate. Comparing the result with SESAME, one  can observe that, with low noise, CEP is not as accurate on angle and amplitude, yet it is more robust with respect to growing  noise. With configuration ({\bf I}), SESAME did not find the cortical source in two cases out of 25 at noise levels 13 and 15 \% with both sample sizes. The sub-cortical source was not detected with 100 samples in 7 cases at 7 \%, 6 cases at 9 \%, 11 cases at 11 \%, 16 cases at 13 \% and 22 cases at 15 \% noise level. For the case of 700 samples, there were 6 cases at 7 \%, 10 cases at 9 \%, 20 cases at 11 \%, 23 cases at 13 \% and 23 cases at 15 \% noise level, where the sub-cortical source was not found.

Our comparison demonstrates that CEP, when implemented in the context of the RAMUS technique, gives an advantage  to reconstruct far-field activities more accurately. Looking at the EMD measurements with multiple noise levels in  Figure \ref{fig:boxplot_EMD}, a clear difference between the prior degrees is observed: the median of the first-degree CEP does not change much, while the EMD of the second-degree CEP is increasing, indicating an increase in the reconstruction mass, i.e. the propagation of the estimated activity, as the localization accuracy does not deteriorate so strongly. For SESAME, the results show a progressive deterioration in both deep and cortical source localization with a significant number of outliers occuring in the case of 700 samples.

As shown by the box plot, the results obtained are generally of lower quality when noise is increased. Considering the MRI-based head mo\-del,   the focality of the estimates obtained with $q = 1$ maintains pronounced compared to the case of $q = 2$ with the 5 \% noise level. Moreover, EM yields a superior reconstruction of  the simultaneous cortical and deep activity of P20/N20 compared to the IAS algorithm. Since the histograms and box plots for elevated noise do not suggest significant differences between EM and IAS, we deem that these  observations might be due to the increased overall level of uncertainty.

\section{Discussion}
\label{sec:discussion}

In this study, we investigated the newly introduced hybrid of  hierarchical Bayesian modelling (HBM) and randomized multiresolution scanning (RAMUS) approach \cite{rezaei2020randomized} as a method to enhance the performance of focal depth localization in EEG. Technically, RAMUS constitutes  a frequentist hybrid solver that is applied to reduce the unknown modelling errors of the Bayesian source localization process. Namely, it relies on the frequentist principle that the actual distribution that the random variables (modelling errors) obey is unknown, but that it is possible to estimate those via sampling. A frequentist model is applied due to the necessity that the modelling accuracy is limited, i.e., there is no {\em a priori} information available beyond some limit, e.g., the resolution of the head model.  While hybrid methods have been introduced and their importance has been shown in various contexts \cite{YuanAo2011Bhmw,YuanAo2009BFHI}, there is no generally agreed approach to combine Bayesian and Frequentist methodologies and, therefore, we consider RAMUS here rather as a technique to improve the posterior optimization process than as a fully  independent statistical method. Future work considering RAMUS as a frequentist method is, however, well-motivated method development goal to find out, e.g., the expected loss of the proposed sample mean estimator \cite{Murphy2012}.

Here the  conditionally Gaussian prior (CGP) model \cite{Calvetti2009}, previously applied in the context of RAMUS in \cite{rezaei2020randomized,RezaeiAtena2021Rsac}, was interpreted as a special case of the conditionally exponential prior (CEP) to improve its focal reconstruction capability. These two priors were compared in numerical experiments using two different reconstruction techniques: the expectation maximization (EM) and iterative alternating sequential (IAS) algorithm, and as an alternative sampling-based technique the recently introduced sequential semi-analytic Monte Carlo estimation (SESAME) technique \cite{SorrentinoAlberto2014Bmmo,sommariva2014}. Whereas EM and IAS find the primary source current density of the brain as a distribution, SESAME reconstructs a finite set of dipolar sources.  

The performance of the CEP was analyzed numerically using both  the spherical three-compartment Ary model and an MRI-based  multi-compartment model. The first one of these was applied in quantitative accuracy and focality analysis considering multiple noise levels and focusing more carefully on the reconstructions within 5 \% noise, and the second one, MRI-based, in a  qualitative investigation, especially, to learn about the possible physiological relevance of the reconstruction difference. The EM and IAS method were implemented in the context of RAMUS to enable the simultaneous detection of both cortical and sub-cortical activities. In particular, as shown in Section \ref{sec:parameter_choice}, the combination of CEP and RAMUS  allows selecting the source-wise shape and scale parameters of the hyperprior through the physiological properties of the brain activity \cite{niedermeyer2004,hamalainen1993} relating the expected reconstruction noise and amplitude levels to those as suggested in \cite{rezaei2020parametrizing}. This combination omits the group effects following from a single focal activity being associated with multiple densely distributed sources. This is possible,  since in RAMUS, the initial reconstruction is first found for a coarse and randomly distributed set of sources (here a set of 10 source positions) and after that propagated towards a denser distribution of sources. If this is not the case, the scale parameters need to be adjusted based on the density of the source space: the greater the source space size, the smaller the scale \cite{rezaei2020parametrizing,he2019zeffiro}. 

To reconstruct the far-field components optimally regardless of the electrode positioning, a multiresolution method such as RAMUS might be required to apply the source sparsity technique  \cite{rezaei2020randomized,krishnaswamy2017sparsity}. Reflecting the observed results with the experimental localization errors obtained with a spherical head model  \cite{cuffin2001spherical}, RAMUS assisted reconstruction techniques are below ($<$ 8 mm including outliers) the experimental  limit found for a far-field source (8.8 -- 11.3 mm in \cite{cuffin2001spherical} with SNR between 25  and 15 dB). Thus, both CGP and CEP, when combined with RAMUS and physiology-based parameter selection, proved to be robust to an increased noise level up to 15 \% (16 dB SNR) of the entry-wise maximum measurement deviation.  This is an interesting finding given that methods have been developed to reduce the effect of noise \cite{LarsonEric2018RSNi,CaiChang2021Reon}. While SESAME outperforms  the other techniques numerically with  low noise, the spread (the interquartile range between 25 and 75 \% quantiles) of its  localization error is greater than the experimentally obtained limit, when the noise level is above 13 \%  (18 dB). In the case of SESAME, the superior dipole  estimation quality obtained  with low noise can at least partially follow from the usage of  dipoles in both generating the  synthetic dataset and reconstructing the activity. This qualitatively implies that SESAME's posterior distribution is relatively narrow and thus very vulnerable when modelling errors and measurement noise are present (similarly to overfitting problems). On the other hand, CEP which employs the distributed source model (in the likelihood) and sparsity priors, has more flexibility in the sense that it does not assume the activity to be limited in only few locations and a strict number of sources.

Based on the results obtained with low or moderate noise, it is evident that RAMUS tends to spread the far-field components, if this tendency is not taken into account {\em a priori}. CEP seems to provide a potential alternative for limiting the  spread that occurs with CGP, while CGP can still be regarded advantageous for distinguishing activity when multiple sources are simultaneously active.
When coupled with RAMUS, the CEP was found to localize deep activity with both prior degrees $q = 1$ and $q = 2$. Relying on the results, the source localization performance provided by the first-order prior can be seen crucial for the focal detection of both near- and  far-field activity. In the case $q = 1$, the near-field  fluctuations (cortical patterns) in configuration ({\bf I}) as well as the far-field (thalamic and sub-thalamic) components in  ({\bf II}) and  ({\bf III}) were observed to be more concentrated to their actual positions  which is potentially significant regarding the identification of the activity in the corresponding brain regions.  For the configuration ({\bf II}) involving simultaneous near- and far-field components,  the deep activity obtained with $q = 1$,  was slightly deviated with respect to the actual position, following clearly from the difficulty to  recover a weakly detectable deep source in the presence of the near-field  activity.  This difference is, nevertheless, minor in the case of the MRI-based head model which distinguishes the sub-cortical areas as disentangled  compartments. Overall, the EM and IAS methods were found to provide qualitatively similar results.

  The combination of first-degree CEP and RAMUS constitutes a potential technique for localizing focal sources with limited distinguishability, e.g., SEP-origin\-at\-ors. Those are classically  observed via invasive depth electrodes, i.e., stereo-EEG, to improve the visibility of the weak components. The present results suggest overall that non-invasive measurements might be successfully used to detect focal deep sources via an efficient source localization strategy, e.g., the CEP. The feasibility of detecting deep activity non-invasively has been shown recently \cite{seeber2019subcortical,pizzo2019deep}. Studies concerning the originators of median nerve SEPs \cite{buchner1994preoperative,buchner1995origin,buchner1997,hsieh1995interaction} associate a significant uncertainty on the early far-field components. Of the median nerve SEP components, the clearest visibility has been obtained for P14/N14 considered as an example case here. P14/N14 is observed when the median nerve SEP propagates within the brainstem. The components P16/N16, modelled in this study,  and P18/N18 occurring  16 and 18 ms post-stimulus, respectively, are likely to involve more than one deep originator and there is yet no exact knowledge on the actual location of those. The P20/N20 component, also considered here, is the first one involving cortical activity, the presence of which, if not taken into account appropriately, might hinder detecting the simultaneous deep activity. Based on our results, the  originators of the P14/N14, P16/N16  and P20/20 component might be non-invasively detectable. Furthermore, the  identification of these might  be enhanced with the  first-degree CEP ($q = 1$).  Surely, the detectability of different SEP originators will need to be carefully studied further both via numerical simulations and with experimental data in order to  gain a deeper insight into the practical applicability of the CEP.

The present source localization approach provides a potential solution for investigating the function and connectivity of the  neural activity networks with focal and weakly detectable far-field components, the analysis of the SEP originators being only one example of potential future  applications. To obtain an optimal performance in connectivity analysis,  RAMUS will need to be employed as a part of a dynamical process, where a reconstruction is obtained based on multiple time points. A  dynamical version of RAMUS can be obtained either via a  straightforward extension of the present single-time sampling (averaging) process or, alternatively, by applying data from multiple points in a single sampling run, where also the source correlations could be taken into account. Our present approach to use a single simulated time step as a basis of the reconstruction can be considered the most relevant, when the investigated processes, e.g., the early SEP components, have a short duration (a few milliseconds), i.e., when the number of available measurements per time point is likely to be small.

\section{Conclusion}
\label{sec:conclusion}
This paper introduced a hybrid of CEP and RAMUS  as a means to enhance the focality of the reconstructions in the localization of near- and far-field sources in EEG measurements.  This combination was shown to allow selecting the shape and scale parameters of the hyperprior relying on a typical brain activity and measurement amplitude, incorporating the noise level due to different uncertainty factors, e.g., measurement and forward modelling errors and thus differentiating it from iterative reweighted $\ell_1$ and $\ell_2$ methods. We observed that CEP, when combined with RAMUS, provides robust source localization estimates up to a 15 \% noise level and that it outperforms SESAME when the noise is high. The first-degree CEP was found to  improve the focality of the source localization estimate compared to the second-degree case, which corresponds to the previously introduced CGP. This improvement was found to be significant especially in the reconstruction of both near- and far-field sources, e.g., to distinguishing   activity in the thalamus  simultaneously with a source in the Brodmann 3b area or a dipolar or a  quadrupolar source configuration focally in the brainstem.  These findings might be crucial, for example, in reconstructing and analyzing SEPs, e.g., the originators of the P20/N20, P14/N14 and P16/N16 components of the median nerve SEP, respectively. 

\section*{Acknowledgments}
JL, AR, and SP were supported by the Academy of Finland Centre of Excellence in Inverse Modelling and Imaging 2018--2025 and ERA PerMed Project (AoF, 344712) "Personalised diagnosis and treatment for refractory focal paediatric and adult epilepsy" (PerEpi). AR was also supported by the Vilho, Yrjö and Kalle Väisälä Foundation and Alfred Kordelini Foundation. AK was supported by the Academy of Finland Postdoctoral Researcher grant number  316542. JL is supported by Väisälä Fund’s (Finnish Academy of Science and Letters) one-year young researcher grant admitted in 2021.

\appendix

\section{Iterative alternating sequential algorithm}
%
\label{sec:app_ias} 
Here we present the standard form of the iterative alternating sequential (IAS) algorithm  (as appeared in \cite{Calvetti2009} for a conjugate hyperprior). In particular, the IAS solves the MAP estimate $({\bf x}^{\hbox{\tiny (MAP)}},{\bm \lambda}^{{\tiny (MAP)}}) = \arg\max\{\pi({\bf x},{\bm \lambda} \mid {\bf y})\}$ as follows :
\begin{algorithm}
\caption{IAS MAP Estimation for Conditionally Exponential Model}\label{alg:IAS}
\begin{algorithmic}
\Procedure{IAS\_MAP}{${\bf x}_{0}$, ${\bf y}$, ${\bf L}$, $q$, $\kappa$, $\theta$, tol}
\State $\hat{\bf x} \gets {\bf x}_{0}$
\While{error $>$ tol}
\State $\hat{\bf x}_{\textnormal{old}} \gets \hat{\bf x}$
\State $\gamma_j \gets (\kappa+q^{-1}-1)/(\theta+\left|\hat{x}_j\right|^q)$, $ j=1,\cdots,3n$
\State $\hat{\bf x} \gets \underset{{\bf x}}{\hbox{arg min}} \left\lbrace \frac{1}{2\sigma^2}\|{\bf L}{\bf x}-{\bf y}\|_2^2+\sum_{i=1}^{3n}\gamma_i|x_i|^q \right\rbrace$
\State error $\gets \left\|\hat{\bf x}-\hat{\bf x}_{\textnormal{old}}\right\|$
\EndWhile
\EndProcedure
\end{algorithmic}

\end{algorithm}

\section{Randomized source space decomposition and RAMUS algorithm}
Here we describe how the source space decompositions are obtained and then we present the pseudoalgorithm of the RAMUS.
\begin{enumerate}
    \item {\em Initialization.} Choose the number of resolution level $R$ and the sparsity factor $s$. The number of sources in each $r=1,...,R$ resolution level is $n_r=\lfloor ns^{r-R} \rfloor$, where $\lfloor \cdot\rfloor$ is the floor function, i.e., rounding down function.
    \item {\em Sampling step.} Generate a desired number $D$ samples of randomized source spaces for each resolution level. These samples are called multiresolution decompositions $\lbrace \mathcal{D}_k^{(r)}\rbrace_{k=1}^D$. The multiresolution decompositions are formed by sampling $n_r$ uniform random center points within active brain tissue. Apply nearest point interpolation scheme between center points and source space such that every point set $B_k$ contains the source points nearest to the $k$th center point. Thereby, disjoint point sets $B_1,...,B_{K_m}$ covers the whole source space and represents the resolution: more there is center point, higher is the resolution. For every $\vec{p}\in B_k$, the unknown ${\bf x}(\vec{p})$ is set to be equal. This holds for every interpolation set respectively.
\end{enumerate}

\begin{algorithm}
\caption{RAMUS algorithm}\label{alg:RAMUS}
\begin{algorithmic}
\Procedure{RAMUS}{${\bf y}$, ${\bf L}$, $q$, $\kappa$, $\theta$, $s$, tol}
\State Generate the source space decompositions and form index sets $\mathcal{P}_{k,r}$, where $k=1,\cdots,N$, $r=1,\cdots,R$
\For{$k=1$ to $N$}\Comment{Loop over decompositions}
\State $\hat{\bf x}_{k,0} \gets {\bf 0}$
\For{$r=1$ to $R$}\Comment{Loop over resolution levels within a decomposition}
\State ${\bf A} \gets {\bf L}(:,\mathcal{P}_{k,r})$
\State Estimate $\hat{\bf x}_{k,r}$ using ${\bf A}$ as the leadfield
\State Interpolate $\hat{\bf x}_{k,r}$ to the finest source space
\State $\hat{\bf x}\gets \hat{\bf x}+\hat{\bf x}_{k,r}/N$
\EndFor
\EndFor
\State $\hat{\bf x}\gets \hat{\bf x}/\sum_{r=1}^Rs^r$\Comment{Scale the estimation appropriately such that $\left\|{\bf L}\hat{\bf x}\right\|/\left\|{\bf y}\right\|\approx 1$}
\EndProcedure
\end{algorithmic}

\end{algorithm}

\section{Solving the Lasso problem via Expectation Maximization}

Here we revisit the sparsity constraint problem (also referred to as Lasso problem) \cite{Tibshirani94regressionshrinkage} and we explain how it can be solved using
 Expectation Maximization (EM) \cite{Murphy2012,Figueiredo2003,Griffin2007}. The Lasso problem \cite{Tibshirani94regressionshrinkage} is to solve
\begin{equation}\label{eq:LassoMin}
\hat{\bf x}:=\underset{\bf x}{\text{arg min}}\left\lbrace\frac{1}{2\sigma^2}\|{\bf L}{\bf x}-{\bf y}\|_2^2+\sum_{i=1}^{3n} \gamma_{i}
|x_i|\right\rbrace,
\end{equation}
where $\gamma_{i}$ are fixed tuning parameters. Using a similar expression as in (\ref{eq:LassoMin}), we can conclude the following Bayesian framework. Given the Gaussian likelihood function of the form 
$\exp{\left(-\frac{1}{2\sigma^2}\|{\bf L} {\bf x}-{\bf y}\|_2^2\right)}$
and the Laplace prior $\pi({\bf x})=\prod_{i=1}^N\pi(x_i)$, where
\begin{equation}\label{eq:LaplacePrior}
\pi(x_i)=\mathrm{Lap}(x_i;0,1/\gamma_{i})=\frac{\gamma_{i}}{2}\exp{\left(-\gamma_{i}|x_i|\right)}
\end{equation}
Thereby, we have that the maximum a posteriori (MAP) estimation obtained via Bayes' rule is equivalent to the Lasso problem (\ref{eq:LassoMin}).
Using the previous decomposition and denoting ${\bf w}^2 = (w_1^2, \ldots, w_n^2)$, we can apply the Expectation Maximization (EM) algorithm \cite{Figueiredo2003,Caron2008} to solve the iteratively the optimization problem
{\scriptsize
\begin{equation}\label{eq:EMoptimization1}
\begin{split}
    {\bf x}^{(j+1)} :=
\underset{{\bf x}}{\mathrm{arg\, max}}\big\{\mathbb{E}_{\pi({\bf w}^2|{\bf x}^{(j)})}[\log
(\pi({\bf y} \mid {\bf x})\pi({\bf x} \mid {\bf w}^2)\pi({\bf w}^2))]\big\}.
\end{split}
\end{equation}}
The above optimization problem can be represented in the form
{\scriptsize
\begin{equation}\label{eq:EMoptimization3}
    {\bf x}^{(j+1)} :=\underset{{\bf x}}{\mathrm{arg\, min}}\left\lbrace\frac{1}{2\sigma^2}\|{\bf L}{\bf x}-{\bf y}\|_2^2+\sum_{i=1}^{3n}\frac{1}{2}x_i^2\mathbb{E}_{\pi(w_i^2|x_i^{(j)})}\left[w_i^2\right]\right\rbrace.
\end{equation}}
The expectation (E) step of EM is given by
$\mathbb{E}_{\pi(w_i^2|x_i^{(j)})}[w_i^2]$ $=$ $\int_{0}^{\infty}w_i^2 \pi(w_i^2|x_i^{(j)})\,\mathrm{d} w_i^2.$
 We can use the Bayes' rule for $\pi(w_i^2|x_i)$ to show that \\
$\mathbb{E}_{\pi(w_i^2|x_i^{(j)})}[w_i^2] =\frac{\gamma_{i}}{|x_i^{(j)}|}.$
That is why, we solve the minimization problem
\begin{equation}\label{eq:EMminimization}
{\bf x}^{(j+1)} :=\underset{{\bf x}}{\mathrm{arg\, min}}\left\lbrace\frac{1}{2\sigma^2}\|{\bf L}{\bf x}-{\bf y}\|_2^2+\sum_{i=1}^{3n}\frac{\gamma_{i}}{2|x_i^{(j)}|}x_i^2\right\rbrace
\end{equation}
which is the maximization (M) step (i.e., the point estimate for ${\bf x}$).
For unimodal posterior we set  $\gamma_{i}=\gamma_{i}/\sigma$ in
(\ref{eq:EMminimization}) \cite{Park2008}.

\end{document}